\makeatletter

\newdimen\paperwidth
\newdimen\paperheight

\def\papersize#1#2{\let\p@persize\relax\paperwidth#1\paperheight#2}

\def\Afour{\papersize{210truemm}{297truemm}}

\let\p@persize\Afour

\let\onesidestyle\@twosidefalse
\let\twosidestyle\@twosidetrue

\def\margins{\@ifnextchar[{\@margins}{\@margins[\z@]}}

\def\@margins[#1]#2#3{
  \p@persize\dimen0 #3\dimen0 .5\dimen0\normalsize%
  \oddsidemargin-1truein\advance\oddsidemargin#2%
  \evensidemargin-1truein\advance\evensidemargin#2%
  \topmargin-1truein\advance\topmargin\dimen0\headsep\dimen0\footskip\dimen0%
  \textwidth\paperwidth\advance\textwidth-#2\advance\textwidth-#2%
  \textheight\paperheight\advance\textheight-#3\advance\textheight-#3%
  \headheight\baselineskip\advance\topmargin-.5\baselineskip%
  \advance\headsep-.5\baselineskip%
  \footheight\baselineskip
  \advance\textwidth-#1\advance\oddsidemargin#1
  \if@twoside\def\@themargin%
    {\ifodd\count\z@\oddsidemargin\else\evensidemargin\fi}\fi}

\def\headlinesep#1{\advance\topmargin\headsep\advance\topmargin -#1
  \advance\topmargin.5\baselineskip\headsep #1\advance\headsep-.5\baselineskip}

\def\headline{\if@twoside\let\n@xt\h@dlin@\else\let\n@xt\h@@dlin@\fi\n@xt}
  
\def\h@dlin@#1#2{%
  \def\@oddhead{%
    {{\leftskip\z@\rightskip\z@\noindent\normalsize#1}}}
  \def\@evenhead{%
    {{\leftskip\z@\rightskip\z@\noindent\normalsize#2}}}}

\def\h@@dlin@#1{%
  \def\@oddhead{{{\leftskip\z@\rightskip\z@\noindent\normalsize#1}}}}

\def\footline{\if@twoside\let\n@xt\f@tlin@\else\let\n@xt\f@@tlin@\fi\n@xt}

\def\f@tlin@#1#2{%
  \def\@oddfoot{%
    {{\leftskip\z@\rightskip\z@\noindent\normalsize#1}}}
  \def\@evenfoot{%
    {{\leftskip\z@\rightskip\z@\noindent\normalsize#2}}}}

\def\f@@tlin@#1{%
  \def\@oddfoot{{{\leftskip\z@\rightskip\z@\noindent\normalsize#1}}}}

\def\normalpage{\global\@specialpagefalse}

\makeatother

\makeatletter

\def\ft{\@ifnextchar[{\ft@s}{\ft@}}
\def\ft@{\ft@@@s[\f@size]}
\def\ft@s[{\@ifnextchar{a}{\ft@sz[}{\ft@@s[}}
\def\ft@@s[{\@ifnextchar{s}{\ft@sz[}{\ft@@@s[}}
\def\ft@@@s[#1]{\ft@sz[at #1pt]}
\def\ft@sz[#1]#2{\font\fonttemp=#2 #1\fonttemp\ignorespaces}

\makeatother

\makeatletter

\input{epsf.sty}

\makeatother

\def\smallcircc{\mathop{\mkern3.5mu\hbox{\raise.58ex\hbox{\ft{lcircle10}a}}}}
\def\varemptyset{{\hbox{\raise.21ex\hbox{$\not$}}\mkern.15mu\mathrm{O}\mkern.15mu}}

  \let\epsilon\varepsilon
      \let\theta\vartheta
          \let\phi\varphi
   \let\emptyset\varemptyset

\documentstyle[12pt]{article}

\makeatletter

\let\Larg@\Large
\let\hug@\huge

\def\usepackage#1{\input{#1.sty}}

\input{geom.sty}

\def\r@adlabel#1#2{\global\@namedef{#1@\the\@key}{#2}}

\let\Large\Larg@
\let\huge\hug@

\def\smallskip{\vskip\smallskipamount}
\def\medskip{\vskip\medskipamount}
\def\bigskip{\vskip\bigskipamount}

\def\mytrivlist{\parsep\parskip\@nmbrlistfalse
  \my@trivlist \labelwidth\z@ \leftmargin\z@
  \itemindent\z@ \def\makelabel##1{##1}}

\def\my@trivlist{\global\@newlisttrue \@outerparskip\parskip}

\def\end#1{\csname end#1\endcsname\@checkend{#1}%
  \expandafter\endgroup\if@endpe\@doendpe\fi
  \if@ignore \global\@ignorefalse \ignorespaces\fi}

\def\put{\@ifnextchar[{\@put}{\@@rput[\z@,\z@][r]}}
\def\@put[#1]{\@ifnextchar[{\@@put[#1]}{\@@@@@put[#1]}}
\def\@@put[#1][{\@ifnextchar{l}{\@@lput[#1][}{\@@@put[#1][}}
\def\@@@put[#1][{\@ifnextchar{c}{\@@cput[#1][}{\@@@@put[#1][}}
\def\@@@@put[#1][{\@ifnextchar{r}{\@@rput[#1][}{\relax}}
\def\@@@@@put[{\@ifnextchar{l}{\@@lput[\z@,\z@][}{\@@@@@@put[}}
\def\@@@@@@put[{\@ifnextchar{c}{\@@cput[\z@,\z@][}{\@@@@@@@put[}}
\def\@@@@@@@put[{\@ifnextchar{r}{\@@rput[\z@,\z@][}{\@@@@@@@@put[}}
\def\@@@@@@@@put[#1]{\@@rput[#1][r]}

\let\hm@d@\leavevmode

\long\def\@@lput[#1,#2][l]#3{\setbox0\hbox{#3}\hm@d@\raise#2\hbox to\z@{\dimen0 #1%
  \advance\dimen0-\wd0\kern\dimen0\dp0\z@\ht0\z@\wd0\z@\box0\hss}\ignorespaces}
\long\def\@@cput[#1,#2][c]#3{\setbox0\hbox{#3}\hm@d@\raise#2\hbox to\z@{\dimen0 #1%
  \advance\dimen0-.5\wd0\kern\dimen0\dp0\z@\ht0\z@\wd0\z@\box0\hss}\ignorespaces}
\long\def\@@rput[#1,#2][r]#3{\setbox0\hbox{\kern#1\raise#2\hbox{#3}}%
  \dp0\z@\ht0\z@\wd0\z@\hm@d@\box0\ignorespaces}

\def\flbox{\@ifnextchar[{\@flbox}{\@@rflbox[\z@,\z@][r]}}
\def\@flbox[#1]{\@ifnextchar[{\@@flbox[#1]}{\@@@@@flbox[#1]}}
\def\@@flbox[#1][{\@ifnextchar{l}{\@@lflbox[#1][}{\@@@flbox[#1][}}
\def\@@@flbox[#1][{\@ifnextchar{c}{\@@cflbox[#1][}{\@@@@flbox[#1][}}
\def\@@@@flbox[#1][{\@ifnextchar{r}{\@@rflbox[#1][}{\relax}}
\def\@@@@@flbox[{\@ifnextchar{l}{\@@lflbox[\z@,\z@][}{\@@@@@@flbox[}}
\def\@@@@@@flbox[{\@ifnextchar{c}{\@@cflbox[\z@,\z@][}{\@@@@@@@flbox[}}
\def\@@@@@@@flbox[{\@ifnextchar{r}{\@@rflbox[\z@,\z@][}{\@@@@@@@@flbox[}}
\def\@@@@@@@@flbox[#1]{\@@rflbox[#1][r]}
\long\def\@@lflbox[#1,#2][l]#3{\@@lput[#1,#2][l]{%
  \vtop{\leftskip\z@\parindent\z@\raggedleft\hm@d@#3}}}
\long\def\@@cflbox[#1,#2][c]#3{\@@cput[#1,#2][c]{%
  \vtop{\leftskip\z@\parindent\z@\raggedcenter\hm@d@#3}}}
\long\def\@@rflbox[#1,#2][r]#3{\@@rput[#1,#2][r]{%
  \vtop{\leftskip\z@\parindent\z@\raggedright\hm@d@#3}}}

\def\maketitle{\par
 \begingroup
 \def\thefootnote{\fnsymbol{footnote}}
 \def\@makefnmark{\hbox to 0pt{$^{\@thefnmark}$\hss}} 
 \if@twocolumn 
 \twocolumn[\@maketitle] 
 \else 
 \global\@topnum\z@ \@maketitle \fi\thispagestyle{plain}\@thanks
 \endgroup
 \setcounter{footnote}{0}
 \let\maketitle\relax
 \let\@maketitle\relax
 \gdef\@thanks{}\gdef\@author{}\gdef\@title{}\let\thanks\relax}

\def\@maketitle{ 
 \null
 \vskip 2em \begin{center}
 {\LARGE \@title \par} \vskip 1.5em {\large \lineskip .5em
\begin{tabular}[t]{c}\@author 
 \end{tabular}\par} 
 \vskip 1em {\large \@date} \end{center}
 \par
 \vskip 1.5em}
 
\def\partbeforeskip#1{\def\p@rtbeforeskip{#1}}
\def\partstyle#1{\def\p@rtstyl@{#1}}
\def\partdot#1{\def\partd@t{#1}}
\def\partafterskip#1{\def\p@rtafterskip{#1}}
\def\partintrostyle#1{\def\partintr@styl@{#1}}
\def\partintrodot#1{\def\partintr@dot{#1}}
\long\def\partintrosep#1{\long\def\partintr@sep{#1}}
\def\partnewpagetrue{\def\p@rtnewp@ge{\newpage}}
\def\partnewpagefalse{\long\def\p@rtnewp@ge{\par}}

\partbeforeskip{4ex}
\partstyle{\centering\Large\bf}
\partdot{}
\partafterskip{3ex}
\partintrostyle{\large}
\partintrodot{}
\partintrosep{\par}
\partnewpagefalse

\def\partname{Part}
\def\part{\p@rtnewp@ge\addvspace\p@rtbeforeskip\@afterindentfalse\secdef\@part\@spart}

\def\@part[#1]#2{\ifnum \c@secnumdepth >-1\relax  
        \refstepcounter{part}                     
        \def\@tempa{\addcontentsline{toc}{part}}  %
        \expandafter\@tempa\expandafter{\thepart  
          \hspace{1em}#1}\else                    
        \addcontentsline{toc}{part}{#1}\fi        
   {\p@rtstyl@                       
    \ifnum \c@secnumdepth >-1\relax        
      {\partintr@styl@\partname\ \thepart  
       \partintr@dot}\partintr@sep\nobreak 
    \fi                                    
    #2\partd@t\markboth{}{}\par}
    \nobreak                       
    \vskip\p@rtafterskip           
   \@afterheading                  
    }                              

\def\@spart#1{{\p@rtcentering\p@rtstyl@                      
    #1\partd@t\par}                 
    \nobreak                        
    \vskip\p@rtafterskip            
    \@afterheading                  
  }                                 

\newif\ifsection@ftind
\newif\ifsection@ftpar

\def\sectionbeforeskip#1{\def\s@ctbeforeskip{#1}}
\def\sectionstyle#1{\def\s@ctstyl@{#1}}
\def\sectiondot#1{\def\sectiond@t{#1}}
\def\sectionafterskip#1{\def\s@ctafterskip{#1}}
\def\sectionintrostyle#1{\def\sectionintr@styl@{#1}}
\def\sectionintro#1{\def\sectionintr@{#1}}
\def\sectionintrodot#1{\def\sectionintr@dot{#1}}
\def\sectionintrosep#1{\def\sectionintr@sep{#1}}
\def\sectionindenttrue{\def\s@ctind{\parindent}}
\def\sectionindentfalse{\def\s@ctind{\z@}}
\def\sectionafterindenttrue{\section@ftindtrue}
\def\sectionafterindentfalse{\section@ftindfalse}
\def\sectionafternewlinetrue{\section@ftpartrue}
\def\sectionafternewlinefalse{\section@ftparfalse}

\newif\ifsubsection@ftind
\newif\ifsubsection@ftpar

\def\subsectionbeforeskip#1{\def\ss@ctbeforeskip{#1}}
\def\subsectionstyle#1{\def\ss@ctstyl@{#1}}
\def\subsectiondot#1{\def\subsectiond@t{#1}}
\def\subsectionafterskip#1{\def\ss@ctafterskip{#1}}
\def\subsectionintrostyle#1{\def\subsectionintr@styl@{#1}}
\def\subsectionintro#1{\def\subsectionintr@{#1}}
\def\subsectionintrodot#1{\def\subsectionintr@dot{#1}}
\def\subsectionintrosep#1{\def\subsectionintr@sep{#1}}
\def\subsectionindenttrue{\def\ss@ctind{\parindent}}
\def\subsectionindentfalse{\def\ss@ctind{\z@}}
\def\subsectionafterindenttrue{\subsection@ftindtrue}
\def\subsectionafterindentfalse{\subsection@ftindfalse}
\def\subsectionafternewlinetrue{\subsection@ftpartrue}
\def\subsectionafternewlinefalse{\subsection@ftparfalse}

\newif\ifsubsubsection@ftind
\newif\ifsubsubsection@ftpar

\def\subsubsectionbeforeskip#1{\def\sss@ctbeforeskip{#1}}
\def\subsubsectionstyle#1{\def\sss@ctstyl@{#1}}
\def\subsubsectiondot#1{\def\subsubsectiond@t{#1}}
\def\subsubsectionafterskip#1{\def\sss@ctafterskip{#1}}
\def\subsubsectionintrostyle#1{\def\subsubsectionintr@styl@{#1}}
\def\subsubsectionintro#1{\def\subsubsectionintr@{#1}}
\def\subsubsectionintrodot#1{\def\subsubsectionintr@dot{#1}}
\def\subsubsectionintrosep#1{\def\subsubsectionintr@sep{#1}}
\def\subsubsectionindenttrue{\def\sss@ctind{\parindent}}
\def\subsubsectionindentfalse{\def\sss@ctind{\z@}}
\def\subsubsectionafterindenttrue{\subsubsection@ftindtrue}
\def\subsubsectionafterindentfalse{\subsubsection@ftindfalse}
\def\subsubsectionafternewlinetrue{\subsubsection@ftpartrue}
\def\subsubsectionafternewlinefalse{\subsubsection@ftparfalse}

\newif\ifparagraph@ftind
\newif\ifparagraph@ftpar

\def\paragraphbeforeskip#1{\def\p@rbeforeskip{#1}}
\def\paragraphstyle#1{\def\p@rstyl@{#1}}
\def\paragraphdot#1{\def\paragraphd@t{#1}}
\def\paragraphafterskip#1{\def\p@rafterskip{#1}}
\def\paragraphintrostyle#1{\def\paragraphintr@styl@{#1}}
\def\paragraphintro#1{\def\paragraphintr@{#1}}
\def\paragraphintrodot#1{\def\paragraphintr@dot{#1}}
\def\paragraphintrosep#1{\def\paragraphintr@sep{#1}}
\def\paragraphindenttrue{\def\p@rind{\parindent}}
\def\paragraphindentfalse{\def\p@rind{\z@}}
\def\paragraphafterindenttrue{\paragraph@ftindtrue}
\def\paragraphafterindentfalse{\paragraph@ftindfalse}
\def\paragraphafternewlinetrue{\paragraph@ftpartrue}
\def\paragraphafternewlinefalse{\paragraph@ftparfalse}

\newif\ifsubparagraph@ftind
\newif\ifsubparagraph@ftpar

\def\subparagraphbeforeskip#1{\def\sp@rbeforeskip{#1}}
\def\subparagraphstyle#1{\def\sp@rstyl@{#1}}
\def\subparagraphdot#1{\def\subparagraphd@t{#1}}
\def\subparagraphafterskip#1{\def\sp@rafterskip{#1}}
\def\subparagraphintrostyle#1{\def\subparagraphintr@styl@{#1}}
\def\subparagraphintro#1{\def\subparagraphintr@{#1}}
\def\subparagraphintrodot#1{\def\subparagraphintr@dot{#1}}
\def\subparagraphintrosep#1{\def\subparagraphintr@sep{#1}}
\def\subparagraphindenttrue{\def\sp@rind{\parindent}}
\def\subparagraphindentfalse{\def\sp@rind{\z@}}
\def\subparagraphafterindenttrue{\subparagraph@ftindtrue}
\def\subparagraphafterindentfalse{\subparagraph@ftindfalse}
\def\subparagraphafternewlinetrue{\subparagraph@ftpartrue}
\def\subparagraphafternewlinefalse{\subparagraph@ftparfalse}

\sectionbeforeskip{\bigskipamount}
\sectionstyle{\large\bf}
\sectiondot{}
\sectionafterskip{.5\bigskipamount}
\sectionintrostyle{}
\sectionintro{}
\sectionintrodot{.}
\sectionintrosep{1.25ex}
\sectionindentfalse
\sectionafterindenttrue
\sectionafternewlinetrue

\subsectionbeforeskip{.8\bigskipamount}
\subsectionstyle{\normalsize\bf}
\subsectiondot{}
\subsectionafterskip{.4\bigskipamount}
\subsectionintrostyle{}
\subsectionintro{}
\subsectionintrodot{.}
\subsectionintrosep{1.25ex}
\subsectionindentfalse
\subsectionafterindenttrue
\subsectionafternewlinetrue

\subsubsectionbeforeskip{.6\bigskipamount}
\subsubsectionstyle{\normalsize\bf}
\subsubsectiondot{}
\subsubsectionafterskip{.3\bigskipamount}
\subsubsectionintrostyle{}
\subsubsectionintro{}
\subsubsectionintrodot{.}
\subsubsectionintrosep{1.25ex}
\subsubsectionindentfalse
\subsubsectionafterindenttrue
\subsubsectionafternewlinetrue

\paragraphbeforeskip{.5\bigskipamount}
\paragraphstyle{\normalsize\bf}
\paragraphdot{.}
\paragraphafterskip{1.25ex}
\paragraphintrostyle{}
\paragraphintro{}
\paragraphintrodot{.}
\paragraphintrosep{1.25ex}
\paragraphindentfalse
\paragraphafterindenttrue
\paragraphafternewlinefalse

\subparagraphbeforeskip{.5\bigskipamount}
\subparagraphstyle{\normalsize\bf}
\subparagraphdot{.}
\subparagraphafterskip{1.25ex}
\subparagraphintrostyle{}
\subparagraphintro{}
\subparagraphintrodot{.}
\subparagraphintrosep{1.25ex}
\subparagraphindenttrue
\subparagraphafterindenttrue
\subparagraphafternewlinefalse

\let\@partoken\par
\long\def\@@gobble#1{}
\def\ignorepar{\@ifnextchar\@partoken{\expandafter\ignorepar\@@gobble}{\ignorespaces}}

\def\@startsection#1#2#3#4#5#6{
   \@tempskipa #4\relax
   \csname if#1@ftind\endcsname\@afterindenttrue\else\@afterindentfalse\fi
   \advance\@tempskipa by\presection
   \if@nobreak \everypar{}\else
     \addpenalty{\@secpenalty}\addvspace{\@tempskipa}%
     \allowbreak\vskip -\presection \fi \@ifstar
     {\@ssect{#1}{#2}{#3}{#4}{#5}{#6}}{\@dblarg{\@sect{#1}{#2}{#3}{#4}{#5}{#6}}}}

\def\@sect#1#2#3#4#5#6[#7]#8{\def\object@type{#1}%
   \ifnum #2>\c@secnumdepth\def\@svsec{}\def\@tempb{}%
      \else\refstepcounter{#1}\def\@svsec{{\csname #1intr@styl@\endcsname%
        {\csname #1intr@\endcsname}\csname the#1\endcsname%
        \csname #1intr@dot\endcsname\kern\csname #1intr@sep\endcsname}}%
        \edef\@tempb{\noexpand\numberline{\csname the#1\endcsname}}\fi%
   \def\@tempa{\addcontentsline{toc}{#1}}%
   \csname if#1@ftpar\endcsname%
      \begingroup #6\relax%
        \@hangfrom{\hskip #3\relax\@svsec}{\interlinepenalty \@M{#8}%
        \csname #1d@t\endcsname\par}%
      \endgroup%
      \csname #1mark\endcsname{#7}%
      \expandafter\@tempa\expandafter{\@tempb #7}%
      \ifautolabel\label*{#8}\fi%
   \else%
      \def\@svsechd{#6\hskip #3\relax%
         \@svsec{#8}\csname #1mark\endcsname{#7}%
         \expandafter\@tempa\expandafter{\@tempb #7}%
         \ifautolabel\label*{#8}\fi}\fi%
   \@xsect{#1}{#5}\ignorepar}

\def\@ssect#1#2#3#4#5#6#7{%
   \ifnum #2>\c@secnumdepth\def\@tempb{}\else \def\@tempb{\numberline{}}\fi%
     \def\@tempa{\addcontentsline{toc}{s#1}}%
     \csname if#1@ftpar\endcsname
        \begingroup #6\relax
           \@hangfrom{\hskip #3}{\interlinepenalty \@M{#7}%
           \csname #1d@t\endcsname\par}%
        \endgroup
        \csname s#1mark\endcsname{#7}%
        \ifstarredcontents\expandafter\@tempa\expandafter{\@tempb #7}\fi%
        \ifautolabel\label*{#7}\fi%
     \else%
        \def\@svsechd{#6\hskip #3\relax{#7}\csname s#1mark\endcsname{#7}%
        \ifautolabel\label*{#7}\fi}\fi
   \@xsect{#1}{#5}\ignorepar}

\def\@xsect#1#2{
   \csname if#1@ftpar\endcsname 
       \par \nobreak \vskip #2\relax \@afterheading
    \else \global\@nobreakfalse \global\@noskipsectrue
       \everypar{\if@noskipsec \global\@noskipsecfalse
                   \clubpenalty\@M \hskip -\parindent
                   \begingroup \@svsechd \endgroup \unskip
                   \hskip #2\relax  
                  \else \clubpenalty \@clubpenalty
                    \everypar{}\fi}\fi\ignorespaces}

\def\section{\@startsection{section}{1}{\s@ctind}
  {\s@ctbeforeskip}{\s@ctafterskip}{\s@ctstyl@}}
\def\subsection{\@startsection{subsection}{2}{\ss@ctind}
  {\ss@ctbeforeskip}{\ss@ctafterskip}{\ss@ctstyl@}}
\def\subsubsection{\@startsection{subsubsection}{3}{\sss@ctind}
  {\sss@ctbeforeskip}{\sss@ctafterskip}{\sss@ctstyl@}}
\def\paragraph{\@startsection{paragraph}{4}{\p@rind}
  {\p@rbeforeskip}{\p@rafterskip}{\p@rstyl@}}
\def\subparagraph{\@startsection{subparagraph}{4}{\sp@rind}
  {\sp@rbeforeskip}{\sp@rafterskip}{\sp@rstyl@}}

\def\statementabove#1{\def\th@bove{#1}}
\def\statementstyle#1{\def\thstyl@{#1}}
\def\statementbelow#1{\def\thb@low{#1}}
\def\statementindentfalse{\let\thind@nt\relax}
\def\statementindenttrue{\let\thind@nt\indent}

\def\statementintrostyle#1{\def\thintr@style{#1}}
\def\statementintrodot#1{\def\thintr@dot{#1}}
\def\statementintrosep#1{\def\thintr@sep{#1}}
\def\statementintrobrackets#1#2{\def\thintr@left{#1}\def\thintr@right{#2}}

\statementabove{\medskip}
\statementstyle{\sl}
\statementbelow{\medskip}
\statementindenttrue

\statementintrostyle{\normalshape\bf}
\statementintrodot{.}
\statementintrosep{\kern1.25ex}
\statementintrobrackets{(}{)}

\def\@thskip{\dimen0\lastskip\vskip-\dimen0%
  \th@bove\dimen1\lastskip\vskip-\dimen1%
  \ifdim\dimen0>\dimen1\else\dimen0\dimen1\fi\vskip\dimen0}

\long\def\@@newtheorem#1#2#3{%
  \newenvironment{#3}%
    {\def\object@type{#3}\par\@thskip%
     \@ifnextchar[{\@enva{#3}{\thstyl@#1{#2}}}{\@envb{#3}{\thstyl@#1{#2}}}}%
    {\end{#3@}}%
  \@ifnextchar[{\@nothm{#3}}{\@nnthm{#3}}}

\def\@nothm#1[#2]#3{%
  \@ifundefined{c@#2}{\@latexerr{No theorem environment `#2' defined}\@eha}%
  {\expandafter\@ifdefinable\csname #1@\endcsname
  {\global\@namedef{the#1}{\@nameuse{the#2}}%
   \global\@namedef{c@#1}{\@nameuse{c@#2}}
   \global\@namedef{p@#1}{\@nameuse{p@#2}}
   \global\@namedef{#1@}{\@nnnthm{#2}{#3}}%
   \global\@namedef{end#1@}{\@endtheorem}}}}

\def\@nnnthm#1#2{\refstepcounter
    {#1}\@ifnextchar[{\@ynnnthm{#1}{#2}}{\@xnnnthm{#1}{#2}}}
 
\def\@xnnnthm#1#2{\@begintheorem{#2}{\csname the#1\endcsname}\ignorespaces}
\def\@ynnnthm#1#2[#3]{\@opargbegintheorem{#2}{\csname
       the#1@\endcsname}{#3}\ignorespaces}

\def\renewtheorem{\@ifnextchar[{\@renewtheorem}{\@renewtheorem[{}{}]}}

\long\def\@renewtheorem[#1]{\@@renewtheorem#1}

\long\def\@@renewtheorem#1#2#3{%
  \expandafter\let\csname#3@\endcsname\undefined
  \renewenvironment{#3}%
    {\def\object@type{#3}\par\@thskip%
     \@ifnextchar[{\@enva{#3}{\thstyl@#1{#2}}}{\@envb{#3}{\thstyl@#1{#2}}}}%
    {\end{#3@}}%
  \@ifnextchar[{\@nothm{#3}}{\@nnthm{#3}}}

\def\@begintheorem#1#2{\@opargbegintheorem{#1}{#2}{}}

\def\@opargbegintheorem#1#2#3{%
        \def\@tempx{#1}%
        \expandafter\let\expandafter\@tempy#2
        \def\@tempz{#3}%
        \mytrivlist\item[\thind@nt\hskip\labelsep%
        {\thintr@style%
          #1\if\@tempx\@empty\else\if\@tempy\relax\else\kern1ex\fi\fi#2%
          \ifx\@tempz\@empty%
            \if\@tempx\@empty\if\@tempy\relax%
            \else\thintr@dot\thintr@sep\fi\else\thintr@dot\thintr@sep\fi%
            \else%
            \if\@tempx\@empty\if\@tempy\relax\else\kern1ex\fi\else\kern1ex\fi%
           \thintr@left{#3}\thintr@right\thintr@dot\thintr@sep\fi}%
            \hskip-\labelsep]%
        \ifautolabel\label*{#3}\fi}

\def\@endtheorem{\strut\endtrivlist\thb@low}

\def\proofabove#1{\def\pf@bove{#1}}
\def\proofstyle#1{\def\pfstyl@{#1}}
\def\proofbelow#1{\def\pfb@low{#1}}
\def\proofindentfalse{\let\pfind@nt\relax}
\def\proofindenttrue{\let\pfind@nt\indent}

\def\proofintrostyle#1{\def\pfintr@style{#1}}
\def\proofintrodot#1{\def\pfintr@dot{#1}}
\def\proofintrosep#1{\def\pfintr@sep{#1}}
\def\proofintrobrackets#1#2{\def\pfintr@left{#1}\def\pfintr@right{#2}}

\proofabove{\medskip}
\proofstyle{}
\proofbelow{\medskip}
\proofindenttrue

\proofintrostyle{\sl}
\proofintrodot{.}
\proofintrosep{\kern1.25ex}
\proofintrobrackets{of\kern1ex}{}

\def\@pfskip{\dimen0\lastskip\vskip-\dimen0%
  \pf@bove\dimen1\lastskip\vskip-\dimen1%
  \ifdim\dimen0>\dimen1\else\dimen0\dimen1\fi\vskip\dimen0}

\renewenvironment{proof}%
  {\@pfskip\mytrivlist\item[\pfind@nt]\@ifnextchar[{\pro@f}{\pro@f[\prooftag]}}
  {\ifvoid\provedbox\else\hproved\fi\endtrivlist\pfb@low}

\def\pro@f[#1]{\setbox\provedbox\hbox{\provedboxcontents{#1}}\proofintro{#1}}

\def\proofintro#1{\expandafter\def\expandafter\@tempa\expandafter{#1}%
  {\pfintr@style{Proof\ifx\@tempa\empty\else\kern1ex\pfintr@left{#1}%
  \pfintr@right\fi}\pfintr@dot\pfintr@sep}\pfstyl@\ignorespaces}

\def\provedmark#1{\def\prm@rk{#1}}
\def\provedsep#1{\def\prs@p{#1}}

\provedmark{$\square$}
\provedsep{\kern1.25ex}

\def\provedtexttrue{\def\prb@x##1{\fbox{\small##1}}}
\def\provedtextfalse{\def\prb@x##1{\prm@rk}}
\def\provedmarkrighttrue{\let\prhf@l\hfill}
\def\provedmarkrightfalse{\let\prhf@l\relax}

\provedtextfalse
\provedmarkrightfalse

\def\provedboxcontents#1{\expandafter\def\expandafter\@tempa\expandafter{#1}%
  \ifx\@tempa\empty\prm@rk\else\prb@x{#1}\fi}

\def\proved{\ifmmode\eqno{\box\provedbox}\else\hproved\fi}

\def\hproved{\unskip\nobreak\prhf@l\penalty50\prs@p\hbox{}\nobreak\prhf@l
  \box\provedbox{\finalhyphendemerits=0\par}}

\def\captionstyle#1{\def\c@ptstyl@{#1}}
\def\captionintrostyle#1{\def\c@pintr@style{#1}}
\def\captionintrodot#1{\def\c@pintr@dot{#1}}
\def\captionintrosep#1{\def\c@pintr@sep{#1}}

\captionstyle{\small\sf}
\captionintrostyle{\bf}
\captionintrodot{.}
\captionintrosep{\hskip1.25ex}

\long\def\@makecaption#1#2{%
    \vskip\captionskip
    \setbox\@tempboxa\hbox{%
      \ifproofing\@ifundefined{the@label}{}
        {\hbox to 0pt{\vbox to 0pt{\vss\hbox{\tiny\the@label}\bigskip}\hss}}\fi
      \c@ptstyl@{\c@pintr@style #1\c@pintr@dot}\ignorespaces #2}%
    \@captionwidth=\hsize \advance\@captionwidth-2\@captionmargin
    \ifdim \wd\@tempboxa >\@captionwidth {%
        \rightskip=\@captionmargin\leftskip=\@captionmargin
        \unhbox\@tempboxa\par}%
      \else
        \hbox to\hsize{\hfil\box\@tempboxa\hfil}%
    \fi}

\def\end@Float#1{%
  \expandafter\caption\expandafter[\the@title]{%
   {\c@pintr@style%
   \ifx\the@caption\empty\ifx\the@title\empty
   \else\c@pintr@sep\fi\else\c@pintr@sep\fi
    \the@title\ifx\the@caption\empty%
     \expandafter\label\expandafter*\expandafter{\the@label}%
    \else\ifx\the@title\empty%
     \expandafter\label\expandafter*\expandafter{\the@label}%
    \else\c@pintr@dot\c@pintr@sep%
     \expandafter\label\expandafter*\expandafter{\the@label}\fi\fi}%
   \ignorespaces\the@caption}%
  \end{#1}}

\renewenvironment{Figure}{\@ifnextchar[%
  {\@myFloat{figure}}{\@myFloat{figure}[htbp]}}{\end@Float{figure}}

\newenvironment{Diagram}{\@ifnextchar[%
  {\@myFloat{diagram}}{\@myFloat{diagram}[htbp]}}{\end{diagram}}

\def\@myFloat#1[#2]#3{%
  \begin{#1}[#2]\def\the@label{#3}}

\def\showfig{\showfigurestrue\fig}
\def\fig#1#2{\@ifnextchar[{\@fig{#1}{#2}}{\@fig{#1}{#2}[0pt]}}

\def\@fig#1#2[#3]#4#5{%
  \def\the@title{#4}\def\the@caption{#5}\centerline{\fig@{#1}{#2}}\vskip#3}

\def\fig@@#1#2{\leavevmode{\figstyl@\vrule width 0pt height 1.8ex%
 \smash{\framebox{\strut\def\@temp{#1}\ifx\@temp\@empty{ #2 }\else{ #1 }\fi}}}}
\def\fig@@@#1#2{\epsfbox{#2}}

\def\figstyle#1{\def\figstyl@{#1}}

\figstyle{\normalshape\mediumseries\normalsize}

\def\oldFigure{%
  \renewenvironment{Figure}{\@ifnextchar[%
    {\@Float{figure}}{\@Float{figure}[htbp]}}{\end@Float{figure}}
    \let\showfig\@ldshowfig \let\fig\@ldfig
    \let\showfigurestrue\@ldshowfigurestrue
    \let\showfiguresfalse\@ldshowfiguresfalse
    \showfiguresfalse}

\def\@ldshowfig#1#2{\epsfbox{#2}}
\def\@ldfig@#1#2{\leavevmode{\framebox{\figstyl@\strut{ #1 }}}}
\def\@ldshowfigurestrue{\let\fig\showfig}
\def\@ldshowfiguresfalse{\let\fig\fig@}

\newcounter{diagram}

\let\thediagram\theequation
\def\ftype@diagram{2}
\def\ext@diagram{lod}
\def\diagram{\@float{diagram}}
\let\enddiagram\end@float
\newif\if@diagnum

\def\diag#1#2{\@ifnextchar[{\@diag{#1}{#2}}{\@diag{#1}{#2}[0pt]}}

\def\@diag#1#2[#3]#4{%
  \def\the@tag{#4}\@eqnswtrue%
  \centerline{\setbox0\hbox{\diag@{#1}{#2}}
  \dimen0 -0.5\wd0\dimen1 0.5\ht0\box0%
  \advance\dimen0 0.5\hsize\advance\dimen0 -\rightskip\advance\dimen1 #3%
  \let\@currentlabel\the@tag%
  \setbox0\hbox to 0pt{\hss\rm\normalshape\mediumseries\normalsize(\the@tag)}%
  \ifx\the@tag\@empty\refstepcounter{equation}\let\@currentlabel\theequation%
    \setbox0\hbox to 0pt{\hss\rm\normalshape\mediumseries\normalsize(\thediagram)}\fi%
  \if@eqnsw\else\let\@currentlabel\relax\setbox0\hbox to 0pt{}\fi%
  \advance\dimen1 -0.5\ht0%
  \put[\dimen0,\dimen1][l]{%
    \box0\expandafter\label\expandafter*\expandafter{\the@label}\kern0.15em}}}

\def\diag@@#1#2{\leavevmode{\diagstyl@\vrule width 0pt height 1.8ex%
 \smash{\framebox{\strut\def\@temp{#1}\ifx\@temp\@empty{ #2 }\else{ #1 }\fi}}}}
\def\diag@@@#1#2{\epsfbox{#2}}

\def\diagstyle#1{\def\diagstyl@{#1}}

\diagstyle{\normalshape\mediumseries\normalsize}

\def\showfiguresfalse{\let\fig@\fig@@}
\def\showfigurestrue{\let\fig@\fig@@@}

\def\showdiagramsfalse{\let\diag@\diag@@}
\def\showdiagramstrue{\let\diag@\diag@@@}

\def\n@number{\@eqnswfalse\let\@currentlabel\relax\let\the@tag\relax}

\def\equation{$$
  \@eqnswtrue\def\object@type{equation}\let\nonumber\n@number%
  \advance\c@equation1\edef\@currentlabel{\theequation}\advance\c@equation-1%
  \def\the@tag{\refstepcounter{equation}\eqno\hbox{\@eqnnum}}}

\def\tag#1{\edef\@currentlabel{#1}\def\the@tag{\eqno\hbox{\reset@font\rm(#1)}}}

\def\endequation{\the@tag$$
  \global\@ignoretrue}

\@definecounter{bibenumi}

\def\thebibliography#1{%
 \section*{\refname}\vskip-\lastskip%
 \list{[\arabic{bibenumi}]}{\topsep0pt\settowidth\labelwidth{[#1]}%
 \leftmargin\labelwidth\advance\leftmargin\labelsep\usecounter{bibenumi}}%
 \def\newblock{\hskip .11em plus .33em minus .07em}%
 \sloppy\clubpenalty4000\widowpenalty4000\sfcode`\.=1000\relax}

\makeatother

\proofingfalse
\autolabelfalse

\showfigurestrue
\showdiagramstrue

\newtheorem{stat}{\statname}  \unnumbered{stat}

\newtheorem{nstat}{\nstatname}[section]

\newtheorem{definition}[nstat]{Definition}
\newtheorem{lemma}[nstat]{Lemma}
\newtheorem{proposition}[nstat]{Proposition}
\newtheorem{theorem}[nstat]{Theorem}
\newtheorem{corollary}[nstat]{Corollary}

\newtheorem[{\ns}{}]{remark}[nstat]{Remark}

\let\ns\normalshape

\papersize{215.9truemm}{279.4truemm}
\margins{3.295cm}{2.12cm}

\headline{\hfill}
\footline{\small\hfill--\kern1ex\thepage\kern1ex--\hfill}
\makeatletter
  \textfloatsep\floatsep
  \intextsep\floatsep
\makeatother

\sectionbeforeskip{1.5\bigskipamount}
\sectionstyle{\centering\normalsize\bf}
\sectionafterskip{\bigskipamount}
\sectionafterindenttrue

\subsectionbeforeskip{\bigskipamount}
\subsectionstyle{\centering\normalsize\sl}
\subsectionafterskip{.5\bigskipamount}
\subsectionafterindenttrue

\paragraphbeforeskip{\bigskipamount}
\paragraphstyle{\centering\normalsize\sl}
\paragraphafterskip{.5\bigskipamount}
\paragraphafternewlinetrue
\paragraphafterindenttrue
\paragraphdot{}

\statementintrostyle{\sc}
\renewtheorem[{\ns}{}]{remark}[nstat]{Remark}
\renewtheorem[{\ns}{}]{definition}[nstat]{Definition}

\captionstyle{\small}
\captionintrostyle{\sc}

\newcommand{\id}{\mathop{\mathrm{id}}\nolimits} 
\newcommand{\Cl}{\mathop{\mathrm{Cl}}\nolimits}

\newcommand{\Fix}{\mathop{\mathrm{Fix}}\nolimits}
\renewcommand{\Im}{\mathop{\mathrm{Im}}\nolimits}
\newcommand{\Ker}{\mathop{\mathrm{Ker}}\nolimits}

\newcommand{\Lk}{\mathop{\mathrm{Lk}}\nolimits}
\newcommand{\St}{\mathop{\mathrm{St}}\nolimits}

\newcommand{\SO}{\mathop{\mathrm{SO}}\nolimits}

\newcommand{\subgr}{\mathrel{\leqslant}}

\newcommand{\nsubgr}{\mathrel{%
  \leqslant\mkern-2.18mu\vrule width .088ex height 1.45 ex depth -0.18 ex}}

\let\emptyset\varemptyset

\def\(#1){$(${\sl #1}\/$)$}

\def\wtM{{\mkern5mu\widetilde{\mkern-5mu{M}\mkern-5mu \mkern5mu}}}

\def\cal#1{{\text{\ft{eusm10}#1}}}

\def\G{{\text{\ft{eusm10}G}}}
\def\M{{\text{\ft{eusm10}M}}}
\def\H{{\text{\ft{eusm10}H}}}

\def\varemptyset{%
 {\text{\raise.21ex\hbox{$\not$}}\mkern.15mu\mathrm{O}\mkern.15mu}}
\def\widebar#1{%
 {\mathchoice
  {{\text{\setbox0\hbox{$#1$}\dimen0\ht0\dimen0.25\dimen0
   \hbox to 0pt{$\kern\dimen0\overline{\kern-\dimen0\phantom{#1}}$\hss}}#1}}
  {{\text{\setbox0\hbox{$#1$}\dimen0\ht0\dimen0.25\dimen0
   \hbox to 0pt{$\kern\dimen0\overline{\kern-\dimen0\phantom{#1}}$\hss}}#1}}
  {{\text{\setbox0\hbox{$\scriptstyle#1$}\dimen0\ht0\dimen0.25\dimen0
   \hbox to 0pt{$\kern\dimen0\overline{\kern-\dimen0\phantom{\scriptstyle#1}}$\hss}}#1}}
  {{\text{\setbox0\hbox{$\scriptscriptstyle#1$}\dimen0\ht0\dimen0.25\dimen0
   \hbox to 0pt{$\kern\dimen0\overline{\kern-\dimen0
	\phantom{\scriptscriptstyle#1}}$\hss}}#1}}}}
\font\ftt cmtt10 at 11pt
\font\fsc cmcsc10 at 12pt
\font\fsl cmsl12 at 12pt

\begin{document}\frenchspacing

\title{\large\bf BRANCHFOLDS AND RATIONAL CONIFOLDS%
       \label{Version 1.00 / \today}}
\author{
\fsc Riccardo Piergallini\\
\fsl Dipartimento di Matematica\\[-3pt]
\fsl Universit\`a di Camerino -- Italy\\
\ftt riccardo.piergallini@unicam.it
\and
\fsc Giacomo Tomassoni\\
\fsl Dipartimento di Matematica\\[-3pt]
\fsl Universit\`a di Camerino -- Italy\\
\ftt giacomo.tomassoni@unicam.it
}
\date{}

\maketitle

\begin{abstract}
\baselineskip13.5pt
\smallskip

\noindent
We extend the concept of orbifold to that of branchfold, in order to allow any cone
singularities with rational angles, and show why branchfolds naturally fit in the theory
of branched coverings. Then, we obtain a geometric goodness theorem for branchfolds and
apply it to prove that a conifold can be endowed with branchfold structure if and only if
it has locally finite holonomy.



\end{abstract}

\section*{Introduction}

This paper is aimed to introduce a class of spaces, which provide an algebraic point of
view for studying conifolds whose codimension two cone singularities have rational angle
of $2k\pi/h$ radians, analogously to what orbifolds do only for angles of $2\pi/h$
radians.

We decided to call these spaces branchfolds, since they naturally fit in the theory of
branched coverings. Roughly speaking, an $m$-dimensional branchfold is a space covered by
open sets $U$ admitting two regular branched coverings $V \leftarrow P \rightarrow U$,
with $P$ a polyhedron and $V$ an open subset of $R^m$ (see the basic example in Figure
\ref{example1/fig} at page \pageref{example1/fig}), rather than only one $V \to U$ as
in the orbifold case.

Since in dimension 3 conifolds with rational angles are dense in the space of all
conifolds, at least in the case of cone manifolds with link singularities (cf.
\cite{HK98}, \cite{Por98} and \cite{CHK00}), in principle branchfold could allow an
algebraic approach to the deformation theory of conifolds, which is a crucial analytic
aspect of the proof of the orbifold geometrization theorem \cite{BLP05}.

On the other hand, branchfolds could also be useful to shed some light on the
Cheeger-Simons problem of whether the volume of a compact spherical conifold with rational
angles is a rational multiple of the volume of the $m$-sphere (cf. \cite{CS85} and
\cite{DS00}).

Here, we limit ourselves to set up the basic theory of branchfolds and to establish the
relation between branchfolds and conifolds, while the above mentioned possible
applications will be considered in future papers.

In Section \ref{prelim/sec} we review the Fox theory of branched covering and prove some
preliminary result. The general setting of branchfold spaces and maps is given in Section
\ref{branchfold/sec}. In particular, Propositions \ref{branchfoldcover/thm},
\ref{coverbranchfold/thm} and \ref{quotientbranchfold/thm} relate branchfold coverings to
branched coverings.

The main results are contained Section \ref{geometry/sec}, where geometric branchfold are
introduced and a geometric goodness theorem (Theorem \ref{goodness/thm}) is proved. As a
consequence of such theorem, we get the above mentioned characterization of the conifolds
which admit a branchfold structure (Theorem \ref{conifoldtobranchfold/thm}), as the those
ones which have locally finite holonomy.

\section{Preliminaries\label{prelim/sec}}

In this section, after providing the needed general setting, we review some stand\-ard
theory of branched coverings between polyhedra in the sense of Fox \cite{Fox57},
reformulating definitions and basic facts in modern language. Then, we consider the notion
of good action, related to that of regular branched covering, and some elementary
properties of pseudo-manifolds, which will be used in the next section to define of
branchfold.

\medskip

By a {\sl polyhedron} we mean a topological space $P$ endowed with a {\sl polyhedral
structure}, that means a PL equivalence class of locally finite countable triangulations.
Such a {\sl triangulation} of $P$ is a locally finite countable simplicial complex $K$
whose underlying topological space is $P = |K|$. Moreover, two triangulations of $P$ are
called {\sl PL equivalent} if and only if they have a common linear subdivision.

\smallskip

We say that a map $f:P \to Q$ between polyhedra is a {\sl PL map} (resp. a {\sl locally PL
map}) if there are triangulations $K$ of $P$ and $L$ of $Q$, such that $f$ sends each
simplex of $K$ linearly {\sl onto} (resp. {\sl into}) a simplex of $L$. In both cases $f$
is called {\sl non-degenerate} if is preserves the dimensions of the simplexes of $K$. 

This terminology is convenient for the present purposes, but {\sl it is not standard
for non-compact polyhedra}, for which PL usually means locally PL in our terms. In
particular, it is worth emphasizing that according our definitions the class of PL maps is
not closed under composition, while such is the class of locally PL maps. However, it can
be shown that locally PL coincides with PL for proper maps.

\smallskip

A subspace $S \subset P$ of a polyhedron is a {\sl subpolyhedron} of $P$ if there exists a
triangulation $K$ of $P$ and a subcomplex $L \subset K$ such that $S = |L|$. In other
words, $S$ itself has a polyhedral structure making the inclusion $S \to P$ a PL map.
On the other hand, any open subspace of a polyhedron $P$ is understood to be a polyhedron,
with the unique polyhedral structure making the inclusion a locally PL map.

\medskip

Given a polyhedron $P$ and a point $x \in P$, we denote by $\St_x P$ (resp. $\Lk_x P$) any
{\sl star} (resp. {\sl link}) of $x$ in $P$, that is the underlying space of the
simplicial star $\St(x,K)$ (resp. link $\Lk(x,K)$) of $x$ in any triangulation $K$ of $P$
having $x$ as a vertex. Notice that $\St_x P$ and $\Lk_x P$ are subpolyhedra of $P$ well
defined up to pseudo-radial PL homeomorphisms centered at $x$.

\smallskip

Polyhedra have local conical structure. In fact, the stars $\St_x P$ gives a basis of {\sl
conical neighborhoods} of $x$ in $P$, being $\St_x P$ the cone of $\Lk_x P$ with apex $x$.
A similar local conical structure is exhibited by non-degenerate (locally) PL maps.
Namely, for any non-degenerate (locally) PL map $f: P \to Q$ and any $x \in P$, putting $y
= f(x)$ we have that $\St_x f = f_|: \St_x P \to \St_y Q$ is the cone of $\Lk_x f = f_|:
\Lk_x P \to \Lk_y Q$.\break
Also such restrictions are well defined up to pseudo-radial PL
homeomorphisms.\break By the {\sl local model} of $f$ at $x$ we mean the open cone
restriction $\St_xf - \Lk_xf$.

\subsection*{Good subpolyhedra\label{goodpoly/sub}}

Now we introduce and discuss the concept of good subpolyhedron. This will play a crucial
role in the definition of branched covering between polyhedra and is directly inspired by
the Fox's approach.

\medskip
\label{goodpoly/def}
Let $P$ be a polyhedron. Then a subpolyhedron $S \subset P$ is called a {\sl good
subpolyhedron} if it is nowhere dense in $P$ (equivalently it does not contain any top
simplex of any triangulation of $P$) and its complement $P - S$ is locally connected at
$S$, meaning that every point $x \in S$ has arbitrarily small neighborhoods $N$ such that
$N - S$ is connected.

\medskip

Here below we collect some elementary properties of good subpolyhedra, which will be
useful in dealing with branched coverings.

\begin{lemma}\label{goodstar/thm}
Let $P$ be a polyhedron and $S \subset P$ be a nowhere dense sub\-polyhedron. Then $S$ is
a good subpolyhedron of $P$ if and only if one of the following equivalent properties
holds:
\begin{itemize}\parskip0pt\itemsep0pt
\item[(1)] \vskip-\lastskip
$\St_x P - S$ is connected for every $x \in S$;
\item[(2)] 
$\Lk_x P - S$ is connected for every $x \in S$.
\end{itemize}\vskip-\lastskip\vskip-\baselineskip
\end{lemma}

\begin{proof} 
The equivalence between the goodness of $S$ and (1) immediately follows from the conical
structure of $\St_x P$, taking into account that $\St_x P \cap S = \St_x S$ when $\St_x P$
is sufficiently small. On the other hand, (1) and (2) are equivalent due to the fact that
$\Lk_x P - S$ is a deformation retract of $\St_x P - S$.
\end{proof}

\begin{proposition}\label{goodsub/thm} 
Let $P$ be a polyhedron and $S \subset P$ be a good subpolyhedron. Then any subpolyhedron
$R \subset S$ is a good subpolyhedron of $P$.
\end{proposition}

\begin{proof}
It suffices to observe that $\St_x P - S$ is dense in $\St_x P - R$ for every $x \in R$,
hence the connectedness of $\St_x P - S$ implies that of $\St_x P - R$.
\end{proof}

\begin{proposition}\label{goodconnect/thm}
If $P$ is a connected polyhedron and $S \subset P$ is a good subpolyhedron, then $P - S$
is connected.
\end{proposition}

\begin{proof}
By the contrary, let $P - S = A \cup B$ with $A$ and $B$ disjoint non-empty open sets.
Since $S$ is nowhere dense, we have $\Cl_P A \cup \Cl_P B = P$. Then, being $P$ connected,
there exists $x \in \Cl_P A \cap \Cl_P B$. Hence we have that $x \in S$ and $\St_x P - S$
meets both $A$ and $B$. So we get a disconnection of $\St_x P - S$, which cannot exist due
to the goodness of $S$.
\end{proof}

\begin{proposition}\label{goodlocal/thm}
Let $P$ be a polyhedron and $S \subset P$ be a nowhere dense sub\-polyhedron. Then $S$ is
a good subpolyhedron of $P$ if and only if one of the following equivalent properties
holds:
\begin{itemize}\parskip0pt\itemsep0pt
\item[(1)]\vskip-\lastskip
$\St_x S$ is a good subpolyhedron of $\St_x P$ for every $x \in S$;
\item[(2)]
$\Lk_x P$ is connected and $\Lk_x S$ is a good subpolyhedron of $\Lk_x P$
for every $x \in S$. \strut
\end{itemize}\vskip-\lastskip\vskip-\baselineskip
\end{proposition}

\begin{proof}
First of all we observe that, since $S$ is nowhere dense in $P$, also $\St_x S$ and $\Lk_x
S$ are nowhere dense respectively in $\St_x P$ and $\Lk_x P$ for every $x \in S$.

Now we prove that (1) holds when $S$ is good in $P$. Being the definition of good
subpolyhedron local in nature, for every $x \in S$ the goodness of the open\break star
 $\St_x S - \Lk_x S \subset \St_x P - \Lk_x P$ is inherited by that of $S \subset P$
(notice that $\St_x S - \Lk_x S = (\St_x P - \Lk_x P) \cap S$ for sufficiently small
stars). Then, the goodness of $\St_x S \subset \St_x P$ follows from the conical structure
of stars.

To prove that (1) implies (2), we fix a point $x \in S$ and assume that $\St_x S$ is
good in $\St_x P$. As above, the conical structure of stars allows us to see that $\Lk_x
S$ is good in $\Lk_x P$. Moreover, $\St_x P - \St_x S$ is connected by Proposition
\ref{goodconnect/thm} and by deformation retraction $\Lk_x P - \Lk_x S$ is connected too.
Therefore $\Lk_x P$ is connect\-ed, since $\Lk_x P - \Lk_x S$ is a dense subspace of it.

Finally, (2) implies that $S \subset P$ is good, by Lemma \ref{goodstar/thm} and
Proposition \ref{goodconnect/thm}.
\end{proof}

We remark that property (2) in the previous Proposition could be used to give an inductive
definition of good subpolyhedron. The induction would be on the (local) dimension of the
ambient polyhedron, starting from the case of dimension 0, where no good subpolyhedron
exists.

\begin{proposition}\label{goodunion/thm}
Let $P$ be a polyhedron and $S \subset P$ be a subpolyhedron.\break If $S$ is union of
good subpolyhedra of $P$, then $S$ is a good subpolyhedron of $P$.
\end{proposition}

\begin{proof}
Since goodness is a local property, it suffices to consider the case when $P$ is finite
dimensional and $S$ is a finite union of good subpolyhedra of $P$. By induction on the
number of such subpolyhedra, we can further reduce ourselves to the special case of $S =
S_1 \cup S_2$ with $S_1, S_2 \subset P$ good subpolyhedra.

The proof that $S$ is good in $P$ proceeds by induction on the dimension of $P$. The base
of the induction is trivially given by $\dim P = 0$. So, we assume $\dim P > 0$ and use
Proposition \ref{goodlocal/thm} to perform the inductive step. Clearly, $S$ is nowhere
dense in $P$. Now, given any point $x \in S$, we have that $\Lk_x P$ is connected since
$x$ is a point of the good subcomplex ($S_1$ or $S_2$) of $P$. Moreover, $\Lk_x S_1$ and
$\Lk_x S_2$ are good subpolyhedra of $\Lk_x P$ and hence $\Lk_x S = \Lk_x S_1 \cup \Lk_x
S_2$ is a good subpolyhedron in $\Lk_x P$ by the inductive hypothesis, being $\dim \Lk_x P
< \dim P$.
\end{proof}

\subsection*{Branched coverings\label{coverings/sub}}

In order to define the notion of branched covering, we need some terminology concerning a
non-degenerate PL map $f: P \to Q$ between polyhedra. We call $x \in P$ a {\sl regular
point} for $f$ if the restriction $f_|:\St_x P \to \St_{f(x)} Q$ is a homeomorphism.
Otherwise, we call $x$ a {\sl singular point} for $f$. \label{singularpoints/def} Then,
the {\sl singular set} $S_f \subset P$, consisting of all the singular points for $f$, is
a subpolyhedron of $P$. We also consider the subpolyhedron $B_f = f(S_f) \subset Q$, that
we call the {\sl branch set} of $f$, and the subpolyheron\break $S'_f = \Cl(f^{-1}(B_f) -
S_f) \subset P$, that we call the {\sl pseudo-singular set} of $f$. Furthermore, just for
notational convenience, we put $T_f = S_f \cup S'_f = f^{-1}(B_f) \subset P$.

\medskip

\label{covering/def}
By a {\sl branched covering} we mean a non-degenerate PL map $f: P \to Q$ between
non-empty polyhedra, which satisfies the following properties:
\begin{itemize}\parskip0pt\itemsep0pt
\item[(1)]\vskip-\lastskip
$S_f$ is a good subpolyhedron of $P$;
\item[(2)]
$B_f$ is a good subpolyhedron of $Q$;
\item[(3)]
$Q$ is connected.
\end{itemize}\vskip-\lastskip

\medskip

We point out that, contrary to the Fox's definition of branched covering, our definition
does not require the covering polyhedron $P$ to be connected. This choice turns out to be
more convenient for the present purposes. However, as Fox himself observes in the first
footnote at page 250 of \cite{Fox57}, this is not an essential point (being relevant only
in defining the universal ordinary covering). It is also worth observing that, according
to our definition, a PL map $f: P \to Q$ is a branched covering if and only if its
restriction $f_|: C \to Q$ is a branched covering for every connected component $C$ of
$P$. On the other hand, the connectedness of the base polyhedron $Q$ is needed for
describing the covering in terms of its monodromy.

Actually, in the Fox's definition $T_f$ is required to be a good subpolyhedron of $P$,
instead of $S_f$. But, in the light of Proposition \ref{covergood/thm} below, this does
not make a real difference.

\begin{proposition}\label{coverprop/thm}
Let $f: P \to Q$ be a branched covering. Then $f$ is a surjec\-tive open map and $T_f$ 
(resp. $B_f$) has local codimension $\geq 2$ in $P$ (resp. $Q$).
\end{proposition}

\begin{proof} Since $f$ is a closed map (as any PL map), while restriction $f_|: P - S_f
\to Q$ is an open map, by the very definition of $S_f$, we have that $f(P - T_f)$ is a
non-empty open and closed subset of $Q - B_f$. On the other hand, according to Proposition
\ref{goodconnect/thm}, properties (2) and (3) of branched coverings imply that $Q - B_f$
is connected.\break Then $f(P - T_f) = Q - B_f$ and hence $f$ is surjective, being $B_f$
nowhere dense.

The same argument above, with $\St_x P$ and $\St_{f(x)}$ respectively in place of $P$ and
$Q$, allows us to see that the restriction $f_|:\St_x P \to \St_{f(x)} Q$ is surjective
for every $x \in P$. This proves that $f$ is an open map.

Now, let $K$ and $L$ triangulations respectively of $P$ and $Q$ with respect to which $f$
is a simplicial map, and let $H \subset K$ the subcomplex such that $S_f = |H|$. 

Then, taking into account that $f$ is non-degenerate, the second part of the statement is
the same as saying that any top simplex $\sigma \in H$ of dimension $m$ is a face of an
$(m + 2)$-simplex of $K$. By the contrary, assume that this is not the case. Then, all
the top simplexes of $K$ containing $\sigma$ have dimension $m + 1$ (remember that $S_f$
is nowhere dense in $P$). Moreover, there is exactly one such top simplex of $K$,
otherwise $S_f$ locally disconnect $P$ (while it does not, being good in $P$ by
Proposition \ref{goodsub/thm}).\break But this easily implies that $\sigma \not\in H$.
\end{proof}

The next proposition says that for a non-degenerate PL map the property of being a
branched covering is a local one and it can be detected by looking at stars/links,
similarly to what happens for goodness of subpolyhedra (cf. Lemma \ref{goodstar/thm} and
Proposition \ref{goodlocal/thm}). We emphasize that the same is not true without assuming
from the beginning that the map is PL.

\begin{proposition}\label{coverlocal/thm}
Let $f: P \to Q$ be a non-degenerate PL map between polyhedra and assume that $Q$ is
connected. Then $f$ is a branched covering if and only if one of the following equivalent
properties holds:
\begin{itemize}\parskip0pt\itemsep0pt
\item[(1)]\vskip-\lastskip
$\St_x f = f_|: \St_x P \to \St_{f(x)} Q$ is a branched covering for every $x \in S_f$;
\item[(2)]
$\Lk_x f = f_|: \Lk_x P \to \Lk_{f(x)} Q$ is a connected branched coverings for every $x
\in S_f$.
\end{itemize}\vskip-\lastskip\vskip-\baselineskip
\end{proposition}

\begin{proof}
First of all we observe that, given $f$ as in the statement and $x \in P$, both the
restrictions $\St_x f$ and $\Lk_x f$ are non-degenerate PL maps. Then, the equivalence
between properties (1) and (2) in the statement can be deduced from the conical structure
of $f$ at $x$, by means of Proposition \ref{goodlocal/thm}. In particular, the
connectedness of $\Lk_x P$ (resp. $\Lk_{f(x)} Q$) is related to the goodness of
$S_{\,\St_x f}$ at $x$ (resp. $B_{\,\St_x f}$ at $f(x)$).

On the other hand, a direct inspection shows that $\St_x S_f = S_{\,\St_x f}$ for every $x
\in S_f$ and $\St_y B_f = \cup_{x \in f^{-1}(y)} B_{\,\St_x f}$ for every $y \in B_f$.
These equalities easily imply that $f$ is a branched covering if and only if it satisfies
property (1) in the statement, by Propositions \ref{goodsub/thm}, \ref{goodlocal/thm} and
\ref{goodunion/thm}.
\end{proof}

We notice that properties (1) and (2) in the previous Proposition concern only the
singular points of the map. However, the same properties hold by definition for the
regular points, except for the fact that the links are not necessarily connected at those
points. Moreover, property (2) could be used to give an inductive definition of branched
covering. The induction starts in dimension $\leq 1$, where singular points cannot exist
and branched coverings are the same as PL ordinary coverings, while the local models in
dimension $\geq 2$ are cones of PL homeomorphisms or branched coverings between connected
compact polyhedra. 

In particular, the first non-trivial local models appear in dimension $2$ and they are
given by the canonical projections $\pi_k: D^2 \to D^2/{{\Bbb Z}_k}$ with $k \geq 2$,
where the action of ${\Bbb Z}_k$ on the open disk $D^2$ is generated by the rotation of $2
\pi/k$ radians at the origin. These are the only local models for branched coverings
between surfaces. \label{flatbranchset} Actually, all the local models for branched
coverings between PL $m$-manifolds at the singular points where the singular set is a
locally flat PL $(m-2)$-submanifold are obtained by crossing the $\pi_k$'s with the
identity of $D^{m-2}$.

\smallskip

Now we want to provide the Fox's characterization of branched coverings as completions of
ordinary coverings. But first we need a technical result concerning the behavior of good
subpolyhedra with respect branched coverings, in order to recover the Fox's definition of
branched coverings, as discussed above.

\begin{proposition}\label{covergood/thm}
Let $f: P \to Q$ be a branched covering and $S \subset P$ be a subpolyhedron. Then $S$ is
a good subpolyhedron of $P$ if and only if $f(S)$ is a good subpolyhedron of $Q$. In
particular, $T_f$ is a good subpolyhedron of $P$.
\end{proposition}

\begin{proof}
Since both the notions of branched covering and good subpolyhedron are local in nature, it
is enough to deal with the case when $\dim P = \dim Q$ is finite. Then, we can proceed by
induction on such dimension, starting from the trivial case of dimension 0. 

To prove the inductive step for dimension $\geq 1$, we observe that $\Lk_x f: \Lk_x P \to
\Lk_{f(x)} Q$ is a branched covering between polyhedra of lower dimension for every $x \in
S$ (see Proposition \ref{coverlocal/thm} and the discussion following it). Thus $\Lk_x S$
is good in $\Lk_x P$ if and only if $f(\Lk_x S)$ is good in $\Lk_{f(x)} Q$, by the
inductive hypothesis. Moreover, $\Lk_x P$ is connected if and only if $\Lk_{f(x)} Q$ is
connected (both links are connected if $x \in S_f$, otherwise they are homeomorphic). Then
the thesis can be obtained by using Propositions \ref{coverlocal/thm}, \ref{goodlocal/thm}
and \ref{goodunion/thm} and taking into account that $\Lk_y f(S) = \cup_{x \in f^{-1}(y)}
f(\Lk_x S)$ for every $y \in f(S)$.
\end{proof}

Given a branched covering $f: P \to Q$, the restriction $g = f_|: P - T_f \to Q - B_f$ is
a PL ordinary covering. In fact, $\St_y Q$ is evenly covered for every $y \in Q - B_f$, by
the very definition of branch set.

As a consequence of the connectedness of $Q - B_f$, all the fibers of $g$ have the same
cardinality $n \leq \infty$. We define the {\sl degree} of $f$ by putting $d(f) = n$ and
call {\sl monodromy} of $f$ the usual monodromy homomorphism $\omega_f: \pi_1(Q - B_f) \to
\Sigma_{d(f)}$ of the covering $g$. In the light of Proposition \ref{coverlocal/thm}, we
also define the {\sl local degree} of $f$ at $x \in P$ as $d_x(f) = d(\St_x f)$ and the
monodromy $\omega_{f,y} : \pi_1(\St_y L - B_f) \to \Sigma_{d(f)}$ of $f$ at $y \in Q$ as
that of the restriction of $f$ over $\St_y Q$. Of course, $\omega_f$ and $\omega_{f,y}$
are determined only up conjugation in $\Sigma_{d(f)}$, depending on the numbering of the
sheets and on the choice of the base point.

The local degree $d_x(f)$, differently from the degree $d(f)$, is always finite for every
$x \in P$, due to the local compactness of polyhedra. Hence, the monodromy $\omega_{f,y}$
has finite orbits (i.e. $\Im \omega_{f,y} \subgr \Sigma_{d(f)}$ has finite orbits as a
group of permutations) for every $y \in Q$. In fact, such orbits correspond to the
restrictions $\St_x f$ with $x \in f^{-1}(y)$.

Actually, the covering space $P$ and the branched covering $f$ can be reconstructed
from the other data, namely from $Q$ and $g$ or equivalently from $Q$ and $\omega_f$. This
can be done by the following completion criterion provided by Fox in \cite{Fox57}.

\begin{proposition}\label{completion/thm}
Let $Q$ be a connected polyhedron, $B \subset Q$ be a good subpolyhedron and $g: R \to Q -
B$ be a PL ordinary covering, whose monodromy $\omega_{g,y}$ at $y$ has finite orbits for
every $y \in B$. Then there exist a simplicial complex $P$, a good subpolyhedron $T
\subset P$, a PL homeomorphism $h: P - T \to R$ and a branched covering $f: P \to Q$
uniquely determined up to PL homeomorphisms, such that $f_|: P - T \to Q - B$ coincides
with $g \circ h$.
\end{proposition}

\begin{proof}
This proposition is nothing more than a restatement in our context of the combination of
three theorems of \cite{Fox57}, namely the two theorems of Section 3, concerning existence
and uniqueness of the completion in the general context of spreads, and the first theorem
of Section 6, about the simplicial case.
\end{proof}

\label{completion/def}
According to Fox, we call the branched covering $f$ given by the above proposition the
{\sl completion} of the ordinary covering $g$ over the complex $Q$. The discussion above
tells us that any branched covering is the completion of the corresponding ordinary
covering.

\medskip

The following elementary results about compositions, factorizations and pullbacks of
branched coverings will be extensively used to define branchfolds and deal with them.

\begin{proposition}\label{covercomp/thm}
Let $f: P \to Q$, $g: Q \to R$ and $g \circ\mkern-2mu f: P \to R$ be PL maps between
polyhedra and assume that $Q$ is connected. If any two of the three maps $f$, $g$ and
$g\circ\mkern-2mu f$ are branched coverings, then also the third one is a branched
covering.
\end{proposition}

\begin{proof}
Taking into account that $S_{g \circ\mkern-2mu f} = S_f \cup f^{-1}(S_g)$ and $B_{g
\circ\mkern-2mu f} = g(B_f) \cup B_g$, the assertion can be easily proved by using
Propositions \ref{goodsub/thm}, \ref{goodunion/thm}, \ref{coverprop/thm} and
\ref{covergood/thm}.
\end{proof}

\begin{corollary}\label{finitecovercomp/thm}
If $f: P \to Q$ and $g: Q \to R$ are branched coverings and $g$ has finite degree, then
the composition $g \circ\mkern-2mu f: P \to R$ is a branched covering.
\end{corollary}

\begin{proof}
Since $g$ has finite degree, the fact that $f$ and $g$ are PL maps implies that also the
composition $g \circ\mkern-2mu f$ is a PL map (cf. \cite{RS72}). Then the corollary is an
immediate consequence of the previous proposition.
\end{proof}

Given two branched coverings $f_1: P_1\to Q$ and $f_2: P_2 \to Q$ with $P_1$ and
$P_2$\break connected, we define their pullback as follows. We put $B = B_{f_1} \cup
B_{f_2}$, $R_1 = P_1 - f_1^{-1}(B)$ and $R_2 = P_2 - f_2^{-1}(B)$. By Propositions
\ref{goodunion/thm} and \ref{covergood/thm}, these are good subpolyhedra of $Q$, $P_1$ and
$P_2$ respectively. Then, we consider the fiber product of $g_1 = f_{1|}: R_1 \to Q - B$
and $g_2 = f_{2|}: R_2 \to Q - B$, consisting of the polyhedron $R = \{(x_1,x_2) \in R_1
\times R_2\,|\,g_1(x_1) = g_2(x_2)\}$ together with the projections $\pi_1: R \to R_1$ and
$\pi_2: R \to R_2$. The maps $\pi_1$, $\pi_2$ and $g = g_1 \circ \pi_1 = g_2 \circ \pi_2:
R \to Q - B$ are ordinary coverings. Hence we can apply Proposition \ref{completion/thm}
to get the corresponding completions $p_1$, $p_2$ and $f$. In particular, $f$ is a
branched covering of degree $d(f) = d(f_1)\,d(f_2)$, with branch set $B_f = B_{f_1} \cup
B_{f_2}$ and monodromy $\omega_f =\break \omega_{f_1} \times \omega_{f_2}: \pi_1(Q - B_f)
\to \Sigma_{d(f_1)} \times \Sigma_{d(f_2)} \subset \Sigma_{d(f_1)\,d(f_2)}$. Now, by the
uniqueness of completions the branched coverings $p_1$, $p_2$ and $f$ can be assumed to
share the same covering space $P$. In this way, they fit into the following commutative
diagram of branched coverings, that is $f = f_1 \circ p_1 = f_2 \circ p_2$.

\vskip12pt
\begin{Diagram}[htb]{pullback1/dia}
\diag{}{pullback1.eps}{}
\end{Diagram}

\medskip

\label{pullback/def}
We call $f: P \to Q$ the {\sl pullback} of the connected branched coverings $f_1: P_1 \to
Q$ and $f_2: P_2 \to Q$. We emphasize that, in spite of the assump\-tion that $P_1$ and
$P_2$ are connected, $P$ is not necessarily connected. It is also worth remarking that
actually the diagram above is not a pullback in the category of PL maps. However,
according to next Proposition \ref{pullback/thm}, it is a pullback in the category of
branched coverings (cf. \cite{Hem90} and \cite{Pie95} for some very special cases).

\medskip

\begin{proposition}\label{pullback/thm}
Let $f_1: P_1 \to Q$ and $f_2: P_2 \to Q$ be connected branched coverings and $f: P \to Q$
be their pullback. Then any branched covering $f': P' \to Q$ which factors through $f_1$
and $f_2$ also factors through $f$. In other words, if there exist branched coverings
$p'_1: P' \to P_1$ and $p'_2: P' \to P_2$ such that $f' = f_1 \circ p'_1 = f_2 \circ
p'_2$, then there exists a PL map $p': P' \to P$ such that $f' = f \circ p'$. Hence, the
diagram\break below is commutative. Moreover, $p'$ restricts to a branched covering over
each connected component of $P$.
\end{proposition}

\begin{Diagram}[htb]{pullback2/dia}\vskip-3pt
\diag{}{pullback2.eps}{}\vskip6pt
\end{Diagram}

\begin{proof}
We consider the subpolyhedra $B = B_f \cup B_{f'} = B_{f_1} \cup B_{f_2} \cup B_{f'}
\subset Q$, $R_1 = P_1 - f_1^{-1}(B) \subset P_1$, $R_2 = P_2 - f_2^{-1}(B) \subset P_2$,
$R = P - f^{-1}(B) \subset P$ and $R' = P' - f'^{-1}(B) \subset P'$, which are good by
Propositions \ref{goodunion/thm} and \ref{covergood/thm}, and the ordinary coverings $g_1
= f_{1|}: R_1 \to Q - B$, $g_2 = f_{2|}: R_2 \to Q - B$, $g = f_|: R \to Q - B$ and $g' =
f'_|: R' \to Q - B$. Since $g'$ factorizes through $g_1$ and $g_2$, we have that
$g'_*(\pi_1(R',x')) \subgr g_{1*}(\pi_1(R_1,x_1)) \cap g_{2*}(\pi_1(R_2,x_2)) =
g_*(\pi_1(R,x)) \subgr \pi_1(Q - B,y)$,\break for any base points $x = (x_1,x_2) \in R$
and $x' \in R'$ such that $p'_1(x') = x_1$ and $p'_2(x') = x_2$ (hence $f'(x') = f(x) = y
\in Q - B$). This allows us to lift componentwise $g'$ through $g$ in order to get an
ordinary covering $h: R' \to R = P - f^{-1}(B)$ such that $g' = g \circ h$. Then the
wanted branched covering $p'$ can be obtained as the completion of $h$ over $P$, by
Proposition \ref{completion/thm}. The uniqueness of completions and liftings gives the
commutativity of the diagram, while the last part of the statement is true by
construction.
\end{proof}

\subsection*{Good actions\label{goodactions/sub}}

Here we recall some basic facts about regular branched coverings and consider the related
notion of good action. We also specialize to the regular case some of the results of the
previous subsection.

\medskip

\label{regular/def}
A branched covering $f: P \to Q$ is called {\sl regular} if there is a group $G$ a PL
automorphisms of $P$ and a PL homeomorphism $h: P/G \to Q$ such that $f = h \circ \pi_G$,
where $\pi_G: P \to P/G$ is the canonical projection.

\medskip

Before going on, we provide an intrinsic characterization of the PL actions on a
polyhedron whose canonical projection is a branched covering. We formally introduce such
actions in Definition \ref{goodaction/def} below, since they are crucial for the notion of
branchfold considered in the next section.

Given a polyhedron $P$ and a group $G$ of PL automorphisms of $P$, the projection $\pi_G:
P \to P/G$ is a non-degenerate PL map onto a polyhedron $P/G$ if and only if the action of
$G$ on $P$ is properly discontinuous, meaning that for any compact subpolyhedron $C
\subset P$ there are only finitely many $g \in G$ such that $g(C) \cap C \neq \emptyset$.
In fact, a standard argument shows that this is in turn equivalent to the existence of a
triangulation $P = |K_G|$ which makes the action simplicial (i.e. $g: K_G \to K_G$ is
simplicial for every $g \in G$). In this case, $\pi_G$ is simplicial with respect to
$K''_G$ (the second barycentric subdivision of $K_G$) and its singular set $S_{\pi_G}$ is
a $G$-invariant subpolyhedron of $P$ triangulated by a subcomplex of $K''_G$. By its very
definition, $S_{\pi_G}$ consists of all the points $x \in P$ whose stabilizer $G_x$ is
bigger than the stabilizer $G_{\St_x P}$ of $\St_x P$ (the\break set of elements of $G$
which fix $\St_x P$ pointwise).

In order to simplify the notation, in the following we will write $S_G$ and $B_G$ for
indicating the corresponding subpolyhedra $S_{\pi_G} = T_{\pi_G}$ (the equality is due
to the $G$-invariance of $S_{\pi_G}$) and $B_{\pi_G}$ associated to the canonical
projection $\pi_G$.

\begin{definition}\label{goodaction/def}
By a {\sl good action} of a group $G$ on a polyhedron $P$ we mean an effective properly
discontinuous PL action of $G$ on $P$, such that $S_G$ is a good subpolyhedron of $P$.
\end{definition}

\begin{proposition} \label{goodaction/thm}
Let $P$ be a polyhedron with an effective PL action of a group $G$ on it, such that $P/G$
is connected. Then, the action is good if and only if the canonical projection $\pi_G: P
\to P/G$ is a branched covering.
\end{proposition}

\begin{proof}
Taking into account the above observations, the only non-trivial fact to be proved is that
the branch set $\pi_G(S_G)$ of $\pi_G$ is a good subpolyhedron of $P/G$ when the action of
$G$ on $P$ is good. In fact, the $G$-invariance of $S_G$ implies that $\St_{\pi_G(x)}P/G -
\pi_G(S_G) = \pi_G(\St_x P - S_G)$ for every $x \in S_G$. Therefore the goodness of
$\pi_G(S_G) \subset P/G$ can be derived from that of $S_G \subset P$, by using Proposition
\ref{goodstar/thm}.
\end{proof}

We observe that, if $f: P \to Q$ is a regular branched covering, then any restric\-tion
$f_|: C \to Q$ to a connected component $C$ of $P$ is still a regular branched covering.
Namely, if $f \cong \pi_G: P \to P/G$, then $f_| \cong \pi_H: C \to C/H$, where $H$ is the
group of PL automorphisms of $C$ consisting of the restrictions of those $g \in G$ such
that $g(C) = C$. Moreover, up to PL homeomorphisms, the branched covering $f_|: C \to Q$
does not depend on the choice the component $C$, since the subgroups of $G$ leaving
invariant different components of $P$ are conjugate in $G$. So, it makes sense to call
$f_|: C \to Q$ the {\sl connected restriction} of the regular branched covering $f: P \to
Q$. \label{conntrestriction/def}

\medskip

In the light of the observation we have just made, the rest if this subsection is focused
on branched coverings whose covering space is connected. For the sake of brevity, we call
them {\sl connected branched coverings}. \label{connbranchedcover/def}

\medskip

By using completions, it can be easily shown that a connected branched covering $f: P \to
Q$ is regular if and only if the associated ordinary covering $g= f_|: P - T_f \to Q -
B_f$ is regular, or equivalently $g_*(\pi_1(P - T_f))$ is a normal subgroup of $\pi_1(Q -
B_f)$ (of course, here we could replace $B_f$ and $T_f$ respectively with $R$ and
$p^{-1}(R)$, where $R \subset Q$ is any good subpolyhedron of $Q$ containing $B_f$).
Moreover, if this is the case and $f$ is induced by the action of a group $G$ on $P$, then
the lifting properties of the connected ordinary covering $g$ imply that the singular set
$S_f = S_G$ consists of all the points $x \in P$ whose stabilizer $G_x$ is non-trivial.

\medskip

Now, we prove some properties of restrictions and factorizations of connected regular
branched coverings. In particular, we introduce the notion of regularization and
characterize those connected branched coverings which admit a regularization, meaning that
they fit into a factorization of a regular one.

\begin{proposition}\label{restaction/thm}
Let $f \cong \pi_G: P \to Q$ be the connected regular branched covering induced by a good
action of $G$ on $P$. Then the restriction $f_|: S \to T$ to connected open subspaces $S
\subset P$ and $T \subset Q$ is still a (connected regular) branched covering if and only
if $S$ is a connected component of $f^{-1}(T)$. In this case, $f_| \cong \pi_H$ is induced
by the good action of the subgroup $H = \{g \in G\;|\;g(S) = S\} \subgr G$ on $S$ given by
restriction.
\end{proposition}

\begin{proof}
Taking into account Propositions \ref{goodconnect/thm} and \ref{completion/thm}, the first
part of the state\-ment can be derived from the analogous property of connected ordinary
coverings by completion. On the other hand, the second part follows from Proposition
\ref{goodaction/thm}, once we observe that $f_|: S \to T$ represents the connected
restriction of $f_|: f^{-1}(T) \to T$ and that the action of $H$ on $S$ is effective, due
to the fact that a deck transformation of a connected ordinary covering is uniquely
determined by the image under it of any given base point.
\end{proof}

\begin{proposition} \label{factaction/thm}
Let $P$ be a connected polyhedron with a good action of a group $G$ on it. Then the
restriction of the action to any subgroup $H \subgr G$ is still good. Moreover the
canonical projection $\pi: P/H \to P/G$ is a branched covering, which is regular if and
only if $H$ is normal in $G$. On the other hand, every branched covering $f: P \to Q$
factorizing the canonical projection $\pi_G: P \to P/G$ is regular, being PL homeomorphic
to $\pi_H: P \to P/H$ for a subgroup $H \subgr G$.
\end{proposition}

\begin{proof}
Given a subgroup $H \subgr G$ as in the statement, we have $S_H \subset S_G$. Then the
restriction of the action of $G$ to $H$ is good by Proposition \ref{goodsub/thm}. Hence
$\pi$ is a branched covering by Propositions \ref{covercomp/thm} and \ref{goodaction/thm}.
Rest of the statement can be obtained by completion, using Propositions
\ref{completion/thm} and \ref{covercomp/thm}, after noticing that the claimed facts hold
for ordinary coverings.
\end{proof}
      
\label{regulariration/def}
By a {\sl regularization} of the connected branched covering $f: P \to Q$ we mean any
connected regular branched covering $r: R \to Q$ which factorizes as $r = f \circ s$ for
some (regular) branched covering $s: R \to P$. We call $r$ the {\sl minimal
regularization} of $f$ if it also satisfies the universal property that, for any other
connected regular branched covering $r' = f \circ s' : R' \to Q$, there exists a (regular)
branched covering $t: R' \to R$ fitting into the commutative diagram below (as a
consequence, if such minimal regularizations exists, then it is uniquely determined up to
PL homeomorphisms).

\begin{Diagram}[htb]{regularization/dia}
\diag{}{regularization.eps}{}
\end{Diagram}

\begin{proposition}\label{regularization/thm}
A connected branched covering $f: P \to Q$ has a regularization if and only if the
monodromy $\omega_{f,y}$ of $f$ at $y$ is finite (has finite image $\Im \omega_{f,y}
\subgr \Sigma_{d(f)}$) for every $y \in Q$. In this case, there also exists the minimal
regularization of $f$, which has finite degree when $f$ has finite degree.
\end{proposition}

\begin{proof}
Assume that $f$ has a regularization $r = f \circ s: R \to Q$ and that this is PL
homeomorphic to the canonical projection $\pi_G: R \to R/G$ induced by a good action of
$G$ on $R$. For any $y \in Q$, the restrictions $\St_x r$ with $x \in r^{-1}(y)$ are
equivalent under the action of $G$, hence their monodromies $\omega_{\St_x r}$ are all
conjugate in $\Sigma_{d(r)}$. As a consequence, since such monodromies have pairwise
disjoint supports, $\omega_{r,y}$ is isomorphic to any single one of them and therefore it
is finite. On the other hand, $\omega_{r,y}$ is invariant under the action of $G$ on $\{1,
\dots, d(r)\}$ induced by the numbering of the sheets. So, we can quotient it by the
action of the subgroup $H \subgr G$ which corresponds to the regular branched covering $s$
according to Proposition \ref{factaction/thm}. This quotient coincides with
$\omega_{f,y}$, up to identification of $\{1, \dots, d(f)\}$ with $\{1, \dots, d(r)\}/H$.
Then, $\omega_{f,y}$ is finite.

Viceversa, assume that the monodromy $\omega_{f,y}$ of $f$ at $y$ is finite for every $y
\in Q$. Then, Proposition \ref{completion/thm} guarantees the existence of the completion
over $Q$ of the ordinary covering of $Q - B_f$, whose fiber is $\Sigma_{d(f)}$ and whose
monodromy is $\omega_f$ acting on $\Sigma_{d(f)}$ by right translations. This completion
is regular just because the original ordinary covering is a regular. Namely, it is PL
homeomorphic to the canonical projection induced by the good action of $\Sigma_{d(f)}$ on
$R$ given by inverse left translations on the non-singular fibers. However it could not be
connected, so we consider its connected restriction $r: R \to P$ and the corresponding
good action of $G \subgr \Sigma_{d(f)}$ on $R$. This is a regularization of $f$. In fact,
a factorizing covering $s: R \to P$, such that $r = f \circ s$, can be realized by
restricting that good action to the subgroup $H \subgr G$ consisting of those elements of
$G$ which fix any given point in $R - T_r$.

To prove the existence of the minimal regularization $f$, we consider the regular ordinary
covering $o: O \to Q - B_f$ such that $M = \Im o_*$ is the largest normal subgroup of
$\pi_1(Q - B_f)$ contained in $f_{|*}(\pi_1(P - T_f))$. If $f$ admits a regularization $r$
as above, then $B_f \subset B_r$ and the inclusion $i: Q - B_r \to Q - B_f$ induces a
surjective homomorphism $i_*: \pi_1(Q - B_r) \to \pi_1(Q - B_f)$, due to the fact that
$B_r$ is a good subpolyhedron of $Q$. Therefore the restriction $r_|: R - T_r \to Q - B_r$
factorizes\break through $o$, being $(i \circ r_{|})_*(\pi_1(R - T_r)) \subgr M$ where $i:
Q - B_r \to Q - B_f$ is the inclusion. This implies that $o$ satisfies the monodromy
hypothesis of Theorem \ref{completion/thm} and hence it can be completed to a regular
branched covering $r_M: R_M \to Q$. The same argument applies to any regularization of $f$
to show that $r_M$ is the minimal one. Concerning the case when $f$ has finite degree $d$,
it suffices to observe that in this case also the regularization constructed above has
finite degree $ \leq d!$ and a fortiori the same holds for the minimal one.
\end{proof}

We conclude this subsection with some facts about pullbacks of connected regular branched
coverings and with the notion of connected pullback.

\begin{proposition}\label{regularpullback/thm}
Let $f_1: P_1 \to Q$ and $f_2: P_2 \to Q$ be connected branched coverings and $f: P \to Q$
be their pullback, which factorizes as $f = f_1 \circ p_1 = f_2 \circ p_2$ (cf.
commutative diagram \ref{pullback1/dia} at page \pageref{pullback1/dia}). If $f_1$ (resp.
$f_2$) is regular, being induced by a good action of a group $G_1$ (resp. $G_2$) on $P_1$
(resp. $P_2$), then also $p_2$ (resp. $p_1$) is regular, being induced by the lifting of
the given action to a good action of the same group on $P$. On the other hand, if both
$f_1$ and $f_2$ are regular as above, then also $f: P \to Q$ is regular, being induced by
the natural good action of $G_1 \times G_2$ on $P$.
\end{proposition}

\begin{proof}
Since we defined the pullback of $f_1$ and $f_2$ as completion of the standard pullback of
certain restrictions of them, the statement can be derived by completion, using
Proposition \ref{completion/thm}, from the usual pulling back group actions through maps.
\end{proof}

As a consequence of Proposition \ref{regularpullback/thm}, in the case when at least one
of the connected branched coverings $f_1: P_1 \to Q$ and $f_2: P_2 \to Q$ is regular, the
restrictions of their pullback $f: P \to Q$ to the connected components of $P$ are all
equivalent up to PL homeomor\-phisms (the same holds for the projections $p_1$ and $p_2$
in diagram \ref{pullback1/dia}).\break We call any such a restriction the {\sl connected
pullback} of $f_1$ and $f_2$. \label{connectedpullback/def}

\medskip

We notice that the connected pullback of $f_1$ and $f_2$ can be thought as the
com\-position of the connected restriction of the regular branched covering $p_2$ (resp.
$p_1$) in diagram \ref{pullback1/dia} with $f_2$ (resp. $f_1$), if $f_1$ (resp. $f_2$) is
regular, while it coincides with the connected restriction of the pullback of $f_1$ and
$f_2$ as branched coverings, if they are both regular.
From a different perspective, it is worth emphasizing that the connected pullback of $f_1$
and $f_2$ is the completion of the ordinary covering $g: R \to Q - B$ such that $B =
B_{f_1} \cup B_{f_2}$ and $g_*(\pi_1(R)) = f_{1|*}(\pi_1(P_1 - f_1^{-1}(B))) \cap
f_{2|*}(\pi_1(P_2 - f_2^{-1}(B)))$ (cf. proof of Proposition \ref{pullback/thm}).

\begin{proposition}\label{connectedpullback/thm}
Let $f_1: P_1 \to Q$ and $f_2: P_2 \to Q$ be connected branched coverings, at least one of
which is regular, and $f: P \to Q$ be the connected pullback of them. Then $f_1$, $f_2$
and $f$ fit into a commutative diagram of connected branched coverings like diagram
\ref{pullback1/dia} at page \pageref{pullback1/dia}, which satisfies the universal
pullback property in the category of the connected branched coverings (cf. Proposition
\ref{pullback/thm}).
\end{proposition}

\begin{proof}
This immediately follows from Propositions \ref{pullback/thm} and
\ref{regularpullback/thm}.
\end{proof}

\subsection*{Pseudo-manifolds\label{pseudoman/sub}}

Most of the polyhedra we will deal with in the next section, in particular branchfolds
themselves, belong to a particular class of polyhedra called pseudo-manifolds. 

\medskip

By an $m$-dimensional {\sl pseudo-manifold} we mean a polyhedron $P$ having a
triangulation $K$ with the following properties, which in this case actually hold for any
triangulation $K$ of $P$:
\begin{itemize}\parskip0pt\itemsep0pt
\item[(1)]\vskip-\lastskip
any top simplex of $K$ has dimension $m$ ($P$ is homogeneously $m$-dimensional);
\item[(2)]
any $(m-1)$-simplex of $K$ is a face of exactly two $m$-simplexes of $K$;
\item[(3)]
$P_{m-2} = |K_{m-2}|$ (the underlying space of the $(m-2)$-skeleton of $K$) is a good
subpolyhedron of $P$.
\end{itemize}\vskip-\lastskip

\medskip

Moreover, $P$ is said to be {\sl orientable} if the top simplexes of $K$ can be coherently
oriented (coherence being required as usual for any two top simplexes sharing an
$(m-1)$-face), while it is said to be {\sl locally orientable} if this can be done for the
top simplexes of the star $\St(x,K)$ for each vertex $x$ of $K$. Such a coherent
orientation of the top simplexes of $P$ (resp. of $\St(x,K)$) is called an {\sl
orientation} (resp. a {\sl local orientation}) of $P$.

\medskip

Of course, PL manifolds without boundary are special examples of locally orientable
pseudo-manifolds. Viceversa, pseudo-manifolds are PL manifolds with certain local
singularities in codimension $> 2$, being locally Euclidean at the points in the interior
of all the simplexes of codimension $\leq 2$. In particular, locally orientable orbifolds
without boundary are locally orientable pseudo-manifolds.

\medskip

We notice that our definition of pseudo-manifold is not completely standard. In fact,
instead of property (3) it is usually required that the connected components of $P$ are
strongly connected, meaning that top simplexes of $K$ in the same component can be joined
by a finite chain of top simplexes of $K$, where any two consecutive of them share an
$(m-1)$-face. By Proposition \ref{goodconnect/thm}, this last fact follows from property
(3), but the viceversa does not hold. So, in this respect our notion of pseudo-manifold is
more restrictive than the usual one.

On the other hand, property (3) presents the advantage of being local
in nature. Then, according to our definition, any open subspace of a pseudo-manifold is
still a pseudo-manifold of the same dimension. Moreover, we have the following
proposition.

\begin{proposition}\label{pm-cover/thm}
Let $f: P \to Q$ be a branched covering. Then $P$ is a pseudo-manifold of dimension $m$
if and only if $Q$ is a pseudo-manifold of the same dimension $m$. In this case,
orientability (resp. locally orientability) of $Q$ lifts to $P$.
\end{proposition}

\begin{proof}
This is an immediate consequence of Propositions \ref{coverprop/thm} and
\ref{covergood/thm}.
\end{proof}

We also notice that properties (1) to (3) of a pseudo-manifold $P$ holds for a
triangulation $K$ of $P$ if and only if they hold for any subdivision $K'$ of $K$ (in
particular, for property (3) this follows from the fact that no $(m-2)$-subpolyhedron can
disconnect an $m$-manifold). Hence, we can conclude that those properties also hold for
any triangulation of $P$.

As a consequence, we have that a subpolyhedron $S \subset P$ is good if and only if it has
local codimension $\geq 2$. In fact property (2) implies that a good subcomplex cannot
contain any codimension 1 simplex, viceversa any subpolyhedron of codimension $\geq 2$ is
good by property (3) and Proposition \ref{goodsub/thm}, being contained into a codimension
$2$ skeleton. 

In the light of this characterization of good subpolyhedra, all the previous subsections
would become simpler in the context of pseudo-manifolds. In particular, the following two
propositions just rewrite the very definition of branched covering and good action in
such context.

\begin{proposition}\label{cover-pm/thm}
Let $f: P \to Q$ be a non-degenerate PL map between pseudo-manifolds of dimension $m$ and
assume that $Q$ is connected. Then $f$ is a branched covering if and only if $\dim S_p
\leq m - 2$ and $\dim B_p \leq m - 2$.
\end{proposition}

\begin{proposition}\label{goodaction-pm/thm}
Let $P$ be pseudo-manifold of dimension $m$. Then a properly discontinuous action of a
group $G$ on $P$ is a good action if and only if $\dim S_G\leq m - 2$. In particular, any
properly discontinuous orientation preserving action on an orientable pseudo-manifold is
good.
\end{proposition}

\section{Branchfolds\label{branchfold/sec}}

This section is entirely devoted to introduce the notion of branchfold and to prove some
fundamental results about branchfolds spaces and related maps. We will consider only
branchfolds without boundary. However, the extension to the bounded case is
straightforward.

Branchfolds generalize locally orientable orbifolds without boundary, in that they admit a
much wider class of singularities. These include codimension two cone singularities with
any rational angle of $2k\pi/h$ radians, with $h$ and $k$ positive coprime integers,
instead of only those of $2\pi/h$ radians allowed for orbifolds.

The idea is just to define an $m$-dimensional branchfold as a spaces covered by open sets
$U$ modelled by two regular branched coverings $V \leftarrow P \rightarrow U$ with $V$
open in $R^m$, rather than only one $V \to U$ like in the orbifold case. 

\begin{Figure}[b]{example1/fig}\vskip6pt
\fig{}{example1.eps}{}{}
\end{Figure}

The basic example is when all the spaces $U$, $P$ and $V$ coincide with the open disk
$D^2$ and the two regular branched coverings are $p_h$ and $p_k$, respectively induced by
the cyclic actions on $D^2$ generated by the rotations of $2\pi/h$ and $2\pi/k$ radians,
with $h$ and $k$ coprime positive integers. This example is depicted in Figure
\ref{example1/fig}.

It is worth remarking that the two group actions modelling a branchfold chart are always
required to generate a single effective action on $P$ as in the above example, but in
general this is not the direct (or even a semidirect) product of them.

As we will see, the setup of these local data, the branchfold charts, requires a quite
long preparatory work, mainly in order to see how they can be glued together to give a
branchfold structure. In principle this is done in the same way as for orbifolds, but
details are more complicate, since we have to take into account of two group actions
instead of only one.

\subsection*{Charts and structures\label{branchfolds/sub}}

Let $X$ be a polyhedron. We want to endow $X$ of an
$m$-branchfold structure. The first step in this direction is to give a local
characterizations of such a structure in terms of branchfold charts.

\begin{definition}\label{chart/def}
An {\sl $m$-branchfold chart} on $X$ is a sextuple $(U, \phi, P, \psi, V, G)$, where: $U
\subset X$ and $V \subset R^m$ are both open; $P$ is a connected polyhedron; $G = H K$ is
a finite group generated by the subgroup $H \subgr G$ and the normal subgroup $K \nsubgr
G$; a good action of $G$ on $P$ is given, such that $\phi: U \to P/H$ is a PL
homeomorphism and $\psi : P/K \to V$ is a piecewise regular smooth homeomorphism letting
the\break induced action of $G/K$ on $P/K$ correspond to an orientation preserving smooth
action on $V$.
\end{definition}

\begin{Diagram}[htb]{chart1/dia}\vskip-15pt
\diag{}{chart1.eps}{}
\end{Diagram}

The commutative diagram shows the relevant actions and maps related to a branchfold chart
as in the definition. In particular, $\pi_G$, $\pi_H$ and $\pi_K$ respectively are the
canonical projections of the good action of $G$ on $P$ and of its restrictions to $H$ and
$K$. On the other hand, $\pi_{G/K}$ is the canonical projection of the induced action of
$G/K$ on $P/K$ (the existence of this last action is the only reason for requiring that
$K$ is normal in $G$). According to Proposition \ref{factaction/thm}, all these maps are
regular branched coverings, while $\pi$ is a (possibly irregular) branched covering.

For the sake of simplicity, in a branchfold chart as above we identify $P/H$ with $U
\subset X$ and $P/K$ with $V \subset R^m$ respectively through $\phi$ and $\psi$. We also
put $p_H = \phi^{-1} \circ \pi_H$, $p_{G/K} = \pi_{G/K} \circ \psi^{-1}$, $p_K = \psi
\circ \pi_K$ and $p = \pi \circ \phi$, to get the simplified commutative diagram below.
In this way, we can omit $\phi$ and $\psi$ from the notation and denote the branchfold
chart by the quadruple $(U, P, V, G = H K)$, where the subgroups $H \subgr G$ and $K
\nsubgr G$ are explicitly indicated.

\begin{Diagram}[htb]{chart2/dia}
\diag{}{chart2.eps}{}
\end{Diagram}

We remark that all the polyhedra in the above diagram are orientable pseudo-manifolds.
This is trivially true for $V$, while it follows from Proposition \ref{pm-cover/thm} for
$P$. Moreover, since the action of $G/K$ is orientation preserving by assumption, also the
actions of $G$ on $P$ is orientation preserving. So, we can apply Proposition
\ref{pm-cover/thm} once again to conclude that $P/G$ and $U$ are orientable
pseudo-manifolds. Actually, $P/G \cong V/(G/K)$ is a locally orientable orbifold.

\label{isochart/def}
Two branchfold charts $(U_1, P_1, V_1, G_1 = H_1 K_1)$ and $(U_2, P_2, V_2, G_2 = H_2
K_2)$ are\break called {\sl isomorphic} if there is an isomorphism $\eta: G_1 \to G_2$,
such that $\eta(H_1) = H_2$ and $\eta(K_1) = K_2$, and an $\eta$-equivariant PL
homeomorphism $h: P_1 \to P_2$ inducing a diffeomorphism $V_1 \cong V_2$. Moreover, those
charts are called {\sl strongly isomorphic} if $U_1 = U_2$ and the PL homeomorphism $U_1
\cong U_2$ induced by $h$ is the identity.

\begin{definition}\label{subchart/def}
An $m$-branchfold chart $(U',P',V',G' = H'K')$ is called a {\sl restriction} of the
$m$-branchfold chart $(U,P,V,G = HK)$, if $p_{H'}$ and $p_{K'}$ are restrictions of $p_H$
and $p_K$ respectively to the open subspaces $U' \subset U$, $P' \subset P$ and $V'
\subset V$.
\end{definition}

Propositions \ref{restaction/thm} and \ref{factaction/thm} tell us that $H' \subgr H$ and
$K' \nsubgr K$, and that the action of $G' \subgr G$ on $P'$ is the restriction of that of
$G$ on $P$. More precisely, we have $H' = \{h \in H \;|\; h(P') = P'\} \subgr H$ and $K' =
\{k \in K \;|\; k(P') = P'\} \subgr K$. However, we warn the reader that $G' \subgr \{g
\in G \;|\; g(P') = P'\}$, but in general these two groups do not coincide (cf. example in
Figure \ref{example2/fig} at page \pageref{example2/fig}).

The commutative diagram below, where $P'\!/G' \to P/G$ is the composition of the inclusion
$P'\!/G' \subset P/G'$ with the canonical projection $P/G' \to P/G$, summarizes how the
restriction chart is related to the original one.

\begin{Diagram}[htb]{restriction/dia}\vskip6pt
\diag{}{restriction.eps}{}\vskip6pt
\end{Diagram}

Of course, chart restriction is a transitive and anti-symmetric binary relation, hence
it induces a partial order on the set of all branchfold charts.

\medskip

\label{smallcharts/rem}
It is worth remarking that a chart restriction as above there always exists with $U'$ an
arbitrarily small open neighborhood of any given $x \in U$. When considering local models,
we will see that such a chart can be chosen to be very special, while here we limit
ourselves to derive its existence from the fact that the map $p: U \to P/G$ in diagram
\ref{chart2/dia} at page \pageref{chart2/dia} is a spread in the sense of Fox
\cite{Fox57}.\break This means that the connected components of the counterimages under
$p$ of the open subsets of $P/G$ form a basis for the topology of $U$. Thus, we can find a
connected open subset $W \subset P/G$ such that $p(x) \in W$ and the connected component
$U'$ of $p^{-1}(W)$ containing $x$ is arbitrarily small. Then, we define $(U',P',V',G' =
H'K')$, by putting $P'$ to be any connected component of $p_H^{-1}(U')$, $V' = p_K(P')$,%
\break $H' = \{h \in H \;|\; h(P') = P'\}$ and $K' = \{k \in K \;|\; k(P') = P'\}$. This
is a branchfold chart by Proposition \ref{restaction/thm}.

The reader should be aware that the restriction chart we have just constructed is not a
generic one, due to the special property that $\pi_{G|}: P' \to \pi_G(P') = W$ is\break a
regular branched covering, being $P'$ a connected component of $\pi^{-1}_G(W)$ (cf.
Proposition \ref{restaction/thm}). This implies that the map $P'\!/G' \to P/G$ can also be
thought as the composition of a regular branched covering $P'\!/G' \to \pi_G(P')$ with the
inclusion $\pi_G(P') \subset P/G$. We will refer to such a restriction as a {\sl special
restriction}. \label{specialsubchart/def}

\medskip

Now, we want to consider equivalent two $m$-branchfold charts when they give isomorphic
local singularities on the base space $X$, independently on the specific polyhedra and
good actions used to describe them. This equivalence relation is generated by that of
domination defined here below.

\begin{definition}\label{chartdomin/def}
An $m$-branchfold chart $(U, P', V, G' = H' K')$ is said to {\sl dominate} the
$m$-branchfold chart $(U, P, V, G = H K)$, which in turn is called a {\sl reduction} of
the first, if there exists a PL map $f: P' \to P$ such that $p_{H'} = p_H \circ f$ and
$p_{K'} = p_K \circ f$.
\end{definition}

In this case $f$ is a regular branched covering. Namely, we can apply Propositions
\ref{covercomp/thm} and \ref{factaction/thm} to the factorizations $p_{K'} = p_K \circ f$
and $p_{H'} = p_H \circ f$, in order to get a normal subgroup $N \nsubgr G'$ contained in
$H' \cap K'$, such that $f \cong \pi_N: P' \to P'\!/N \cong P$, $G'\!/N \cong G$ and, up
to these isomorphisms, we have the following commutative diagram.

\begin{Diagram}[htb]{domination/dia}\vskip-6pt
\diag{}{domination.eps}{}\vskip6pt
\end{Diagram}

Viceversa, any normal subgroup $N \nsubgr G'$ contained in $H' \cap K'$ gives raise in
this way to domination of $(U, P', V, G' = H' K')$ on $(U, P \cong P'/N, V, G = H K \cong
H'\!/N\,K'\!/N \cong G'\!/N)$ through the map $f \cong \pi_N$.

\medskip

Like restriction, also domination is a transitive and anti-symmetric binary relation,
hence it induces a partial order on the set of all branchfold charts. Moreover, any
special restriction of a domination (resp. reduction) chart is a domination (resp.
reduction) of a special restriction. This property does not hold in general, without
requiring the restrictions to be special. \label{specialsubchart/rem}

\begin{definition}\label{chartequiv/def}
Two $m$-branchfold charts $(U, P_1, V, G_1 = H_1 K_1)$ and $(U, P_2,\break V, G_2 = H_2
K_2)$ are called {\sl equivalent} when there exists a third $m$-branchfold chart $(U, P',
V, G' = H' K')$ which dominates both of them, as in the following commutative diagram.
\end{definition}

\begin{Diagram}[htb]{equivalence1/dia}\vskip-12pt
\diag{}{equivalence1.eps}{}\vskip6pt
\end{Diagram}

In order to see that chart equivalence is a true equivalence relation, the only
non-trivial property to be verified is transitivity. Since chart domination is a
transitive relation, it suffices to show that for any two charts $(U, P_1, V, G_1 = H_1
K_1)$ and $(U, P_2, V, G_2 = H_2 K_2)$ dominating the same chart $(U, P, V, G = H K)$
there exists a chart $(U, P', V, G' = H' K')$ which dominates both of them. This is done
in the commutative diagram \ref{equivalence2/dia}. We start with the regular branched
coverings $f_1: P_1 \to P$ and $f_2: P_2 \to P$ giving the assumed dominations. Then we
consider their connected pullback $q = f_1 \circ q_1 = f_2 \circ q_2: Q \to P$ and the
minimal regularization $r = \pi_G \circ q \circ s:\break P' \to P/G$ of the composition
$\pi_G \circ q$ (which is a possibly irregular branched covering by Corollary
\ref{finitecovercomp/thm}). The other PL maps $p_{H'}$, $p_{K'}$, $f'_1$ and $f'_2$ are
defined just by composition.
Proposition \ref{regularization/thm} allows us to think of $r$ as the canonical projection
$\pi_{G'}$ of a good action of a finite group $G'$ on $P'$. By Proposition
\ref{factaction/thm}, also $p_{H'}$ and $p_{K'}$ are regular branched coverings
corresponding to the restrictions of that action to certain subgroups $H', K' \subgr G'$.
Moreover $K'$ is normal in $G'$, since $p_{G/K}$ is regular. On the other hand, it is
clear from the diagram that $G' = H' K'$, hence we have a branchfold chart $(U, P', V, G'
= H' K')$. As desired, this chart dominates $(U, P_1, V, G_1 = H_1 K_1)$ and $(U, P_2, V,
G_2 = H_2 K_2)$ respectively through $f'_1: P' \to P_1$ and $f'_2: P' \to P_2$.

\begin{Diagram}[htb]{equivalence2/dia}\vskip6pt
\diag{}{equivalence2.eps}{}\vskip6pt
\end{Diagram}


At this point, we have all the ingredients needed to express the compatibility condition
for two branchfold charts on a polyhedron $X$, in order to belong to the same branchfold
structure on $X$. So, we can proceed with our main definitions.

\begin{definition}\label{atlas/def}
An {\sl $m$-branchfold atlas} on $X$ is a set $\cal A = \{(U_i, P_i, V_i, G_i = H_i
K_i)\}_{i\,\in\,I}$ of $m$-branchfold charts such that $\cal U = \{U_i\}_{i\,\in\,I}$ is
an open covering of $X$ and the following {\sl compatibility condition} holds: for any
$i,j \in I$, $x \in U_i \cap U_j$, $\widetilde x_i \in P_i$ and $\widetilde x_j \in P_j$
such that $p_{H_i}(\widetilde x_i) = p_{H_j}(\widetilde x_j) = x$, there exist two
restrictions $(U, P'_i, V'_i, G'_i = H'_i K'_i)$ and $(U, P'_j, V'_j, G'_j = H'_j K'_j)$
of the charts $(U_i, P_i, V_i, G_i = H_i K_i)$ and $(U_j, P_j, V_j, G_j = H_j K_j)$
respectively, with $x \in U$, $\widetilde x_i \in P'_i$ and $\widetilde x_j \in P'_j$,
which are equivalent up to strong isomorphisms of charts.
\end{definition}

The compatibility condition is summarized in the commutative diagram
\ref{compatibility/dia}, where $(U', P', V', G' = H' K')$ is the dominating branchfold
chart which gives the equivalence between the restrictions, $V'_i \cong V' \cong V'_j$
represents diffeomorphisms, $f_i$ and $f_j$ are the domination maps.

\begin{Diagram}[htb]{compatibility/dia}\vskip6pt
\diag{}{compatibility.eps}{}\vskip6pt
\end{Diagram}

\begin{definition}\label{branchfold/def}
An {\sl $m$-branchfold} (or $m$-dimensional branchfold) is a pair $X_{\cal B} = (X,\cal
B)$, where $X$ is a polyhedron and $\cal B$ is an {\sl $m$-branchfold structure}, meaning
a maximal $m$-branchfold atlas, on $X$. We will write $X$ in place of $X_{\cal B}$,
omitting the reference to the branchfold structure, when no confusion can arise.
\end{definition}

\label{atlas->structure}
A standard argument shows that any $m$-branchfold atlas determines a unique $m$-branchfold
structure containing it. In fact, the compatibility condition at a fixed $x \in X$ turns
out to be an equivalence relation on the charts $(U, P, V, G = HK)$ such that $x \in U$.
In particular, the main point for the transitivity is the possibility of choosing all the
restriction and domination charts in two consecutive compatibility conditions to be based
on the same open $U' \subset X$. This could be seen by using the special restrictions
considered at page \pageref{specialsubchart/def} and their properties with respect to
dominations and reductions (see page \pageref{specialsubchart/rem}). However, since the
claimed transitivity will become evident in the next subsection, after expressing the
compatibility condition in terms of local models (cf. observation following Definition
\ref{localmodel/def}), we do not go into the details of the proof here.

\medskip

We emphasize that $m$-branchfolds, like locally orientable $m$-orbifolds without boundary,
are particular $m$-dimensional pseudo-manifolds. This is a consequence of Proposition
\ref{pm-cover/thm} and of the local nature of our notion of pseudo-manifold. On the other
hand, locally orientable $m$-orbifolds without boundary coincide with the special
$m$-branchfolds which admit a branchfold atlas $\cal A = \{(U_i, P_i, V_i, G_i = H_i
K_i)\}_{i\,\in\,I}$ such that $K_i = 1$ for every $i \in I$. Indeed, such branchfold
charts trivially reduce to usual orientable orbifold charts and the same is true for the
compatibility condition between any two of them (cf. \cite{BS85}, \cite{MM91},
\cite{Sat56} or \cite{Thu79}).

Finally, we say that a branchfold is {\sl orientable} (resp. {\sl oriented}\/) referring
to the underlying pseudo-manifold. Moreover, given an oriented branchfold $X$, by an
oriented chart (resp. atlas) on $X$ we mean a chart which is (resp. an atlas whose all
charts are) oriented coherently with $X$.

\subsection*{Local models\label{localmodels/sub}}

The local model of a branchfold $X$ at a point $x \in X$ can be characterized in terms of
the conical branchfold charts centered at $x$ which are minimal with respect to the
domination order. In fact, such charts turn out to be all isomorphic. To prove this, we
first formalize the definition of conical branchfold chart and show that all the conical
restrictions at $x$ of a given branchfold chart are isomorphic. Then, we show that each
equivalence class of branchfold charts contains a unique representative that is minimal in
the above sense. Finally, we put these facts together to get the claimed unicity up to
isomorphism of the minimal conical charts at $x$.

\begin{definition}\label{conicalchart/def}
A branchfold chart $(U, P, V, G = HK)$ of an $m$-branchfold $X$ is called {\sl conical} if
$P$ is an open cone, the action of $G$ on $P$ is conical (meaning that it preserves the
cone structure of $P$) and the induced cone structure on $V$ is linear. Moreover, we say
that the chart is {\sl centered at $x$}, when $x \in X$ is the apex of the cone structure
induced on $U$.
\end{definition}

Up to chart isomorphisms, in a conical chart as above we can always assume $V = R^m$ with
apex at the origin, since any starlike open set is diffeomorphic to $R^m$, and $G/K <
\SO(m)$ acting on $R^m$ by Euclidean isometries, being $G/K$ identifiable with a finite
group of orientation preserving linear isomorphisms of $R^m$.

Given a branchfold chart $(U, P, V, G = HK)$ and a point $x \in U$, we can construct
arbitrarily small branchfold conical restrictions $(U', P', V', G' = H'K')$ of $(U, P, V,
G = HK)$ centered at $x$, in the following way.

Let $\widetilde x \in P$ be any point such that $p_H(\widetilde x) = x$ and let $\widebar
x = p_K(\widetilde x) \in V$. Consider (arbitrarily small) open stars $U' = \St_x U -
\Lk_x U \subset U$, $P' = \St_{\widetilde x} P - \Lk_{\widetilde x} P \subset P$\break and
$V' = \St_{\widebar x} V - \Lk_{\widebar x} V \subset V$, and the stabilizer
$G_{\widetilde x} \subgr G$. Moreover, put $H' = H_{\widetilde x} = H \cap G_{\widetilde
x} \subgr H$, $K' = K_{\widetilde x} = K \cap G_{\widetilde x} \subgr K$ and $G' = H' K' =
H_{\widetilde x} K_{\widetilde x} \subgr G_{\widetilde x} \subgr G$. Proposition
\ref{restaction/thm} tells us that $(U', P', V', G' = H'K')$ is a branchfold chart, hence
a conical restriction of $(U, P, V, G = HK)$.

We emphasize that in general $G'$ is a proper subgroup $G_{\widetilde x}$. The simplest
(conical) chart $(U, P, V, G = HK)$ where this happen is depicted in Figure
\ref{example2/fig}. Here, $H = \langle \sigma \rangle \cong {\Bbb Z}_2$, $K = \langle \rho
\rangle \cong {\Bbb Z}_2$ and $G = \langle \sigma,\rho \rangle \cong {\Bbb Z}_2 \times
{\Bbb Z}_2$, while the heavy lines represent cone singularities, whose angle is $2 \pi$
divided by the corresponding numeric label (this is the index that will be introduced in
Definition \ref{localmodel/def}). Then, for any $x \neq 0$ along the vertical axis of $U
\cong R^3$ and any $\widetilde x$ such that $p_H(\widetilde x) = x$, we have $G' \cong 1$
while $G_{\widetilde x} = \langle \sigma\mkern-1mu\rho \rangle \cong {\Bbb Z}_2$.

\begin{Figure}[htb]{example2/fig}
\fig{}{example2.eps}{}{}
\end{Figure}

As it is natural to expect, all the conical restrictions of $(U, P, V, G = HK)$ centered
at $x$ are isomorphic. In fact, any such restriction can be obtained by the above
construction, and this turns out to be independent on the choice of the lifting
$\widetilde x$, up to conjugation by an element of $H$ (which leaves invariant both the
subgroups $H$ and $K$, being the latter normal in $G$), and on the choice of the specific
realizations of the stars, up to pseudo-radial PL homeomorphisms.

\begin{definition}\label{reducedchart/def}
A branchfold chart $(U, P, V, G = HK)$ is called {\sl reduced} if it does not properly
dominate any other chart (in other words, it is minimal with respect to the domination
order).
\end{definition}

The discussion of domination following Definition \ref{chartdomin/def}, implies that $(U,
P, V, G = HK)$ is reduced if and only if $H \cap K$ does not contain any non-trivial
normal subgroup of $G$. Moreover, any branchfold chart $(U, P, V, G = HK)$ dominates a
unique\break (up to isomorphisms) reduced chart $(U, P', V, G' = H'K')$, which is given by
$P' = P/N$, $G' = G/N$, $H' = H/N$ and $K' = K/N$, where $N \nsubgr G$ is the maximal
normal subgroup of $G$ contained in $H \cap K$.

Consequently, any equivalence class of branchfold charts has a unique reduced
representative (up to isomorphisms), being any two equivalent reduced charts dominated by
the same chart.

\begin{proposition}\label{localmodel/thm}
Let $X$ be an $m$-branchfold. Then for any point $x \in X$ there exists a unique (up to
isomorphisms) reduced conical chart $(U_x, P_x, V_x, G_x = H_x K_x)$ of $X$ centered at
$x$.
\end{proposition}

\begin{proof}
We have already said that conical charts of $X$ centered at $x$ do exist. Then, the
existence of a reduced one follows from the trivial observation that any reduction of a
conical chart is still conical.

To prove the unicity up to isomorphisms, let $(U_1, P_1, V_1, G_1 = H_1 K_1)$ and $(U_2,
P_2, V_2, G_2 = H_2 K_2)$ any two reduced conical charts of $X$ centered at $x$. By the
compatibility condition there exists equivalent restrictions $(U, P'_1, V'_1, G'_1 = H'_1
K'_1)$ and $(U, P'_2, V'_2, G'_2 = H'_2 K'_2)$. Possibly passing to smaller ones, these
restrictions can be assumed to be conical. In this case, they are isomorphic to the
original charts, as conical restrictions of conical charts (cf. the above discussion on
conical charts). Then, they are also isomorphic to one another, being equivalent and
reduced (cf. the above discussion on reduced charts). So, we can conclude by using
transitivity.
\end{proof}

\begin{definition}\label{localmodel/def}
By a branchfold {\sl local model} we mean any reduced conical branchfold chart. In
particular, given a branchfold $X$ and a point $x \in X$, we define the {\sl local model}
of $X$ at $x$ to be the reduced conical chart $(U_x, P_x, V_x, G_x = H_x K_x)$, whose
existence and unicity are ensured by the previous proposition. Moreover, we call $G_x =
H_x K_x$ the {\sl isotropy group} of $X$ at $x$ (this notation should not be confused with
that for stabilizers, since $G_x$ acts on $P_x$ while $x \in U_x$) and $i_x =
|H_x|/|K_x|$, where $|\cdot|$ denotes the cardinality, the {\sl index} of $X$ at $x$.
\end{definition}

We observe that the two restrictions, as well as the intermediate chart dominating them,
in the compatibility condition for the charts $(U_i, P_i, V_i, G_i = H_i K_i)$ and $(U_j,
P_j, V_j, G_j = H_j K_j)$ in Definition \ref{atlas/def} can be assumed to be conical. As
an immediate consequence, those charts satisfy the compatibility condition at $x \in U_i
\cap U_j$ if and only if their conical restrictions centered at $x$ reduce to the same
local model $(U_x, P_x, V_x, G_x = H_x K_x)$ up to isomorphism.

\medskip

In terms of isotropy groups, locally orientable orbifolds without boundary can be
characterized as the branchfolds $X$ such that $K_x \cong 1$ for every $x \in X$. In this
case $G_x = H_x$ and $G_x/K_x \cong G_x$. Moreover, there are isomorphisms $p_{K_x}: P_x
\cong V_x \subset R^m$ and $p_x: U_x \cong P_x/G_x$ allowing us to identify $p_{H_x}$ with
$p_{G_x/K_x}$, which is always an orbifold local model, even for $K_x \ncong 1$.

At the opposite end of the branchfold spectrum, we call $X$ a {\sl pure branchfold} when
$H_x \subset K_x$ for every $x \in X$. In this case $G_x = K_x$ and $G_x/K_x \cong 1$,
hence we have an isomorphism $p_{G_x/K_x}: V_x \cong P_x/G_x$. Up to this isomorphism,
$p_{K_x} \cong p_x \circ p_{H_x}$ is the minimal regularization of the branched covering
$p_x: U_x \to P_x/G_x \cong V_x$. Therefore, $X$ is locally modelled on (conical) branched
coverings of $R^m$. The simplest non-trivial example of such a local model is depicted in
Figure \ref{example3/fig} (heavy lines and numerical labels have the same meaning as in
Figure \ref{example2/fig}), where $H = \langle \sigma \rangle \cong {\Bbb Z}_2$ and $G = K =
\langle \sigma,\rho \rangle \cong \Sigma_3$, while $p$ is the cone of the covering $S^2
\to S^2$ branched over three points with monodromies $(1\;2)$, $(2\;3)$ and $(1\;2\;3)$
respectively.

\begin{Figure}[htb]{example3/fig}
\fig{}{example3.eps}{}{}
\end{Figure}

\begin{definition}\label{singularset/def}
Given an $m$-branchfold $X$, we define the {\sl singular locus} of $X$ to be the good
subpolyhedron $\Sigma X = \{x \in X \;|\; G_x \neq 1\} \subset X$. We think of $\Sigma X$
as a stratified set of dimension $\leq m-2$, with the natural stratification such that the
local model is constant on the connected components of the strata.
\end{definition}

The fact that $\Sigma X$ is good in $X$ immediately follows from Propositions
\ref{goodsub/thm}, \ref{goodlocal/thm} and \ref{covergood/thm}, taking into account that
$\St_x \Sigma X \subset p_{H_x}(S_{G_x})$ for every $x \in \Sigma X$.

The natural stratification $\Sigma_0 X \subset \Sigma_1 X \subset \dots \subset
\Sigma_{m-2} X = \Sigma X$ can be obtained by defining $\Sigma_i X$ as the set of all $x
\in \Sigma X$ such that, for any triangulation of $P_x$ making the action of $G_x$
simplicial in the local model $(V_x, P_x, U_x, G_x = H_x K_x)$, the simplex $\sigma$ of
the induced triangulation of $U_x$ containing $x$ in its interior has $\dim \sigma \leq
i$. Taking into account the product structure of the chart over a neighborhood of the
interior of $\sigma$, it is straightforward to verify that in this way we get a
stratification. Moreover, such a product structure also implies that the local model is
the same for all points in the interior of $\sigma$, hence this is contained in $\Sigma_i
X$. For $x \in \Sigma_i X - \Sigma_{i-1}X$, there exists a simplex $\sigma$ as above with
$\dim \sigma = i$, which is a top simplex of $\Sigma_i X$. Then, the local model turns out
to be locally constant on the stratum $\Sigma_i X - \Sigma_{i-1}X$, so it is constant on
each connected component of it.

\medskip 

The next proposition focus on the $(m-2)$-stratum $\Sigma_{m-2}X - \Sigma_{m-3}X$ of
$\Sigma X$. Namely, it tells us that the local models at this stratum are particularly
simple, being the cartesian product of the basic example depicted in Figure
\ref{example1/fig} with $R^{m-2}$.

\begin{proposition}\label{singularset/thm}
Let $X$ be an $m$-branchfold and let $x \in \Sigma X$ be a singular point in the
$(m-2)$-stratum of $\,\Sigma X$. Then the local model $(V_x, P_x, U_x, G_x = H_x K_x)$ is
given by $U_x \cong P_x \cong V_x \cong R^m$, $H_x = \langle \rho_{2 \pi/h_x} \rangle
\cong {\Bbb Z}_{h_x}$, $K_x = \langle \rho_{2 \pi/k_x} \rangle \cong {\Bbb Z}_{k_x}$ and $G_x
= \langle \rho_{2 \pi/(h_x k_x)}\rangle \cong H_x \times K_x \cong {\Bbb Z}_{h_x k_x}$,
where $h_x$ and $k_x$ are the unique coprime positive integers such that $i_x = h_x/k_x$,
while $\rho_\alpha$ denotes the rotation of $\alpha$ radians around $R^{m-2} \subset R^m$.
\end{proposition}

\begin{proof}
Let us assume $V_x = R^m$ and $G/K \subgr \SO(m)$ acting by Euclidean isometry on $V_x$
(cf. observation following Definition \ref{conicalchart/def}). Since $x \notin
\Sigma_{m-3} X$, the branch set of the branched covering $p_{K_x}$ is the empty set or an
$(m-2)$-dimensional subspace of $R^m$. Then, as already observed at page
\pageref{flatbranchset}, $p_{K_x}$ is isomorphic to the canonical projection induced by
the action on $P_x \cong R^m$ generated by a rotation $\rho_{2 \pi/k_x}$ for some positive
integer $k_x$, that is $K_x = \langle \rho_{2 \pi/k_x}\rangle$. Now, also the singular set
of the branched covering $p_{H_x}$ is the empty set or an $(m-2)$-dimensional subspace of
$P_x \cong R^m$, and in this last case it must coincide with the singular set of
$p_{K_x}$. Therefore, also the action of $H_x$ is generated by a rotation $\rho_{2
\pi/h_x}$ for some positive integer $h_x$, that is $H_x = \langle \rho_{2 \pi/h_x}
\rangle$. As a consequence, we also have $U_x \cong R^m$. Finally, since $(V_x, P_x, U_x,
G_x = H_x K_x)$ is a reduced chart, we have that $h_x$ and $k_x$ are coprime and hence
$G_x = \langle \rho_{2 \pi/(h_x k_x)} \rangle \cong H_x \times K_x \cong {\Bbb Z}_{h_x
k_x}$.
\end{proof}

In the light of the previous proposition, for any connected component $C$ of the
$(m-2)$-stratum $\Sigma_{m-2} X - \Sigma_{m-3} X$ of $\Sigma X$, the branchfold structure
of a neighborhood of $C$ is completely determined by its index $i_C$, defined as the
common index\break $i_x$ at all the points $x \in C$. So, it makes sense to label each
such component $C$ with $i_C$, just as we have already done in Figures \ref{example2/fig}
and \ref{example3/fig}.

Of course, these labels have integer values if $X$ is an $m$-orbifold, and in this case
they coincide with the customary ones, while they have values of the type $1/n$ if $X$ is
a pure $m$-branchfold. The viceversa does not hold in general, but it trivially holds when
$\Sigma_{m-3} = \emptyset$, that is when the underlying space of $X$ is a PL $m$-manifold
and $\Sigma X$ is a PL locally flat $(m-2)$-submanifolds of $X$.

\subsection*{Maps and coverings\label{branchfoldmaps/sub}}

This subsection is entirely devoted to introduce branchfold maps and coverings and to
relate them to branched coverings. Branchfold coverings form a special class of maps,
which extends that of orbifold coverings to the branchfold context.

\begin{definition}\label{branchfoldmap/def}
Let $f: X \to Y$ be a map between branchfolds. We call $f$ a {\sl branchfold map} if for
every $x \in X$ there exist a branchfold chart $(U, P, V, G = HK)$ of $X$ with $x \in U$,
a branchfold chart $(U', P', V', G' = H'K')$ of $Y$ with $f(U) \subset U'$,\break a PL map
$\phi: P \to P'$ and a regular smooth map $\psi: V \to V'$, such that the following
diagram commutes. Moreover, we call $f$ a {\sl branchfold isomor\-phism} if in addition it
is a PL homeomorphism (then also $f^{-1}: Y \to X$ is a branchfold map).
\end{definition}

\begin{Diagram}[htb]{map1/dia}\vskip-3pt
\diag{}{map1.eps}[-3.5pt]{}\vskip9pt
\end{Diagram}

We notice that definition of branchfold map is local in nature. In particular, branchfold
maps are locally PL, since so is the map $f_|$ in the diagram, but they are not
necessarily PL. We also notice that the only reason why the composition of two branchfold
maps could not be a branchfold map, is the possible lack of regularity of the map $\psi$
in the diagram.

By considering conical restrictions of both the branchfold charts involved in diagram
\ref{map1/dia}, we can assume that $\psi: V \to V'$ is a linear map of maximum rank
between linear spaces and $f_|$ is a PL map. Actually $(U', P', V', G' = H'K')$ could be
assumed to be a local model, but in general one could not insist that also $(U, P, V, G =
HK)$ is a local model. 

\begin{Diagram}[b]{map2/dia}
\diag{}{map2.eps}[4pt]{}\vskip-9pt
\end{Diagram}

In particular, when $X$ and $Y$ have the same dimension, the map $\psi: V \to V'$ can be
assumed to be an isomorphism. In this case, Propositions \ref{covercomp/thm} and
\ref{factaction/thm} imply that $\phi$ is a regular branched covering, and then
Proposition \ref{covercomp/thm} and Corollary \ref{finitecovercomp/thm} imply that $f_|$
is a branched covering. Moreover, there exists a further branched covering $\rho: P/G \to
P'/G'$ which completes diagram \ref{map1/dia} to give the commutative
diagram~\ref{map2/dia}. Finally, by completion and the lifting properties of ordinary
coverings, we can show that there exists a group homomorphism $\eta: G \to G'$ such that
$\phi$ is $\eta$-equivariant, $\eta(H) \subset H'$ and $\eta(K) \subset K'$.

If $f: X \to Y$ is a branchfold isomorphism, then $\dim X = \dim Y$ and the map $f_|$ in
diagram \ref{map2/dia} is a PL homeomorphism. As a consequence, such diagram gives a
domination of $(U, P, V, G = HK)$ on $(U', P', V', G' = H'K')$, hence these branchfold
charts are equivalent. Thus, branchfold isomorphisms can be characterized as those PL
homeomorphisms that induce isomorphisms of local models at all points.

\medskip

Now, in order to give the definition of branchfold covering, we need the concept of
rectifiability of a path in a finite dimensional polyhedra.

Let $P$ be a finite dimensional polyhedron. Once a locally PL inclusion $P \subset R^n$ is
given, we can define the length $L(\alpha)$ of a path $\alpha: [0,1] \to P$ in the usual
way. Of course, $L(\alpha)$ depends on the inclusion. However, the intrinsic metrics
induced by different inclusions are locally Lipschitz equivalent, hence the compactness of
$\alpha([0,1])$ implies that the property $L(\alpha) < \infty$ does not depend on the
inclusion. 

\medskip

\label{rectifiable/def}
We say that the path $\alpha: [0,1] \to P$ is {\sl rectifiable}, or that it has {\sl
finite length}, when $L(\alpha) < \infty$, for any locally PL inclusion $P \subset R^n$.

\medskip

If $M$ is a smooth manifold, then the usual notion of rectifiability with respect to any
Riemannian metric on $M$, coincides with the above one, when $M$ is thought with the
polyhedral structure induced by any smooth triangulation of it.

\begin{definition}\label{branchcover/def}
A {\sl branchfold covering} is a branchfold map $f: X \to Y$ be-\break tween branchfolds
of the same dimension, which is {\sl complete} with respect to lifting of rectifiable
paths, meaning that any partial lifting $\widetilde \alpha_{|}: \left[0,1\right[ \to X$ of
a rectifiable path $\alpha: [0,1] \to Y$ extends to a complete lifting $\widetilde
\alpha:[0,1] \to X$.
\end{definition}

The completeness property of branchfold coverings specializes to locally PL maps between
polyhedra the notion of completeness given by Fox in the much more general context of
spreads (cf. \cite{Fox57}, pages 244 and 248).

According to the above remarks about branchfold maps, a branchfold covering $f: X \to Y$
is locally a branched covering at every point $x \in X$. Moreover, extending to the
locally PL map $f$ the notions of singular and branch set (as defined at page
\pageref{singularpoints/def} for PL maps), we have that $S_f \subset \Sigma X \cup
f^{-1}(\Sigma Y)$ and $B_f \subset f(\Sigma X) \cup \Sigma Y$.

We emphasize that $f$ it is not required to be a PL map, hence it is not necessarily a
branched covering. Actually, $S_f$ turns out to be always a good subcomplex of $X$, but
$B_f$ can be a quite pathological subset of $Y$, for example it can be dense in $Y$ (cf.
footnote~\ddag\ at page 250 of \cite{Fox57}). However, branchfold coverings behave better
than branched covering with respect to composition, being the composition of branchfold
coverings always a branchfold covering.

Arguing on diagram \ref{map1/dia} for the local models of a branchfold covering $f: X \to
Y$, we can easily realize that $Y$ is an orbifold if $X$ is and that $X$ is a pure
branchfold if $Y$ is. Moreover, in the orbifold case, our definition of branchfold
covering restricts to the usual one of orbifold covering, since the map $\phi$ in diagram
\ref{map1/dia} is forced to be a PL homeomorphism. We also note that in this case
$f(\Sigma X) \subset \Sigma Y$, hence there is no room for pathological branch sets.

\medskip

The next propositions establish some relations between PL branchfold coverings and
branched coverings. In particular, Proposition \ref{coverbranchfold/thm} tells us that any
branched covering $f: X \to Y$ of a given brachfold $Y$ can be interpreted as a branchfold
cov\-ering, by a suitable choice of the branchfold structure on $X$. An analogous result
holds in the context of orbifolds only for proper maps, if one allows both the orbifold
structures on $X$ and $Y$ to be suitably choosen. On the other hand, Proposition
\ref{quotientbranchfold/thm} extends the construction of good orbifolds as quotients of
smooth manifolds by the action of a discrete group of diffeomorphisms.

\begin{proposition}\label{branchfoldcover/thm}
Let $X$ and $Y$ be branchfolds of the same dimension, and assume that $Y$ is connected.
Then, any PL branchfold map $f: X \to Y$ is a branchfold covering and it is a branched
covering as a map between polyhedra.
\end{proposition}

\begin{proof}
The local properties of branchfold maps between branchfold of the same dimension ensure
that $f$ is non-degenerate. Then, the first part of the statement derives from the
observation that any non-degenerate PL map is complete, while the second part immediately
follows from Proposition \ref{coverlocal/thm}.
\end{proof}

\begin{proposition}\label{coverbranchfold/thm}
Let $f: X \to Y$ be a branched covering onto a branchfold $Y$. Then, the branchfold
structure of \ $Y$ can be lifted to a unique branchfold structure on $X$ making $f$ into 
a PL branchfold covering.
\end{proposition}

\begin{proof}
As the first step, we prove that any given branchfold chart $(U, P, V, G = HK)$ on $Y$ can
be lifted to a branchfold chart $(\widetilde U, \widetilde P, \widetilde V, \widetilde G =
\widetilde H \widetilde K)$ on $X$, such that the two charts are related by a commutative
diagram like diagram \ref{map1/dia}. This is done in diagram \ref{map3/dia}, which is
constructed in the following way. Start with the diagram of the chart $(U, P, V, G = HK)$
on the right side of diagram \ref{map3/dia}. Choose an arbitrary connected component
$\widetilde U$ of $f^{-1}(U)$ and consider the connected pullback $p_H \circ q_1 = f_|
\circ q_2: Q \to U$ of $p_H$ and $f_|: \widetilde U \to U$. Then, let $r = \pi_G \circ q_1
\circ s$ be the minimal regularization of the branched coverings $\pi_G \circ q_1$.
Finally, put $\widetilde V = V$ and define the remaining maps by composition.

\begin{Diagram}[htb]{map3/dia}
\diag{}{map3.eps}[12pt]{}
\end{Diagram}

We denote by $L$ the group of deck transformations of the covering $r$ and define
$\widetilde H \subgr L$ and $\widetilde K \nsubgr L$ to be the subgroups corresponding to
the coverings $p_{\widetilde H}$ and $p_{\widetilde K}$ respectively (see Proposition
\ref{factaction/thm}). Then, we put $\widetilde G = \widetilde H \widetilde K \subgr L$,
getting in this way a branchfold chart $(\widetilde U, \widetilde P, \widetilde V,
\widetilde G = \widetilde H \widetilde K)$ on $X$, which is related to the original chart
$(U, P, V, G = HK)$ on $Y$ by the maps $f_|$, $\phi$ and $\psi = \id$, as in diagram
\ref{map1/dia}.

Now, we have to verify that the lifted charts satisfy the compatibility condition of
Definition \ref{atlas/def}, in such a way that they form a branchfold atlas (they
obviously cover all of $X$). According to the observation following Definition
\ref{localmodel/def}, it suffices to show that the local model $(\widetilde U_x,
\widetilde P_x, \widetilde V_x, \widetilde G_x = \widetilde H_x \widetilde K_x)$, obtained
from the lifted chart $(\widetilde U, \widetilde P, \widetilde V, \widetilde G =
\widetilde H \widetilde K)$ by reduction of the conical restriction centered at any $x \in
\widetilde U$, is equivalent to the lifting $(\widetilde U_y, \widetilde P_y, \widetilde
V_y, \widetilde G_y = \widetilde H_y \widetilde K_y)$ of the local model $(U_y, P_y, V_y,
G_y = H_y K_y)$ of $Y$ at $y = f(x) \in U$.

To this end, we first consider the conical restriction $(U', P', V', G' = H'K')$ of the
chart $(U, P, V, G = HK)$ centered at $y$ and the lifting $(\widetilde U', \widetilde P',
\widetilde V', \widetilde G' = \widetilde H' \widetilde K')$ of it such that $x \in
\widetilde U'$. We want to prove that such lifting is chart equivalent to\break the
conical restriction $(\widetilde U', \widetilde P'', \widetilde V', \widetilde G'' =
\widetilde H'' \widetilde K'')$ of $(\widetilde U, \widetilde P, \widetilde V, \widetilde
G = \widetilde H \widetilde K)$, where $\widetilde P''$ is any connected component of
$\phi^{-1}(P')$ such that $p_{\widetilde H}(\widetilde P'') = \widetilde U'$.

The construction of $(\widetilde U', \widetilde P', \widetilde V', \widetilde G' =
\widetilde H' \widetilde K')$ is described by the diagram, obtained from diagram
\ref{map3/dia} just by putting a prime on spaces and maps. We refer to this new diagram as
diagram \ref{map3/dia}$'$, without drawing it. Comparing diagrams \ref{map3/dia} and
\ref{map3/dia}$'$, we see that there are open inclusions $U' \subset U$, $P' \subset P$,
$V' \subset V$, $\widetilde U' \subset \widetilde U$ and $Q' \subset Q$, the first three
by hypothesis and the last two by construction. Moreover, the maps $p_{H'}$, $p_{K'}$,
$f_|$, $q'_1$ and $q'_2$ in diagram \ref{map3/dia}$'$ are restrictions of the
corresponding maps in diagram \ref{map3/dia}.

\begin{Diagram}[htb]{map4/dia}\vskip-6pt
\diag{}{map4.eps}[12pt]{}
\end{Diagram}

Then, we can construct the commutative diagram \ref{map4/dia}. Here, $s_|$, $\pi_{G|}$ and
$p_{G/K|}$ are restrictions of maps in diagram \ref{map3/dia}. In particular, $s_|$ is a
regular branched covering by Proposition \ref{restaction/thm}, while $\pi_{G|}$ is a
(finite) regular branched covering if we assume, as we obviously can, that the conical
restriction $(U', P', V', G' = H'K')$ is special in the sense specified at page
\pageref{smallcharts/rem}. In this case, also $p_{G/K|}$ is regular branched covering by
Proposition \ref{restaction/thm}. Moreover, $t = s' \circ t_1 = s_| \circ t_2 : T \to Q'$
is the connected pullback of $s'$ and $s_|$, while $\widehat r = \pi_G \circ q'_1 \circ t
\circ \widehat s: \widehat P \to \pi_G(P')$ is the minimal regularization of $\pi_G \circ
q'_1 \circ t$ (this is a branched covering by Corollary \ref{finitecovercomp/thm}, having
$q'_1$ and $\pi_{G|}$ finite degree). Finally, the maps $p_{\widehat H}$, $p_{\widehat
H}$, $f_1$ and $f_2$ are define by composition.

By putting $\widehat G = \widehat H \widehat K$, where $\widehat H$ and $\widehat K$ are
the subgroups of the group of deck transformations of $r$ corresponding to $p_{\widehat
H}$ and $p_{\widehat K}$ respectively, we get a new chart $(\widetilde U', \widehat P,
\widetilde V', \widehat G = \widehat H \widehat K)$ on $X$. Looking at diagram
\ref{map4/dia}, it is straightforward to verify that this chart dominates both
$(\widetilde U', \widetilde P', \widetilde V', \widetilde G' = \widetilde H' \widetilde
K')$ and $(\widetilde U', \widetilde P'', \widetilde V', \widetilde G'' = \widetilde H''
\widetilde K'')$ through the maps $f_1$ and $f_2$ respectively, giving us the claimed
equivalence between these two charts.

Since $(\widetilde U_x, \widetilde P_x, \widetilde V_x, \widetilde G_x = \widetilde H_x
\widetilde K_x)$ is a reduction of $(\widetilde U', \widetilde P'', \widetilde V',
\widetilde G'' = \widetilde H'' \widetilde K'')$ by its\break very definition, we can
conclude our argument by showing that the reduction of $(U', P', V', G' = H'K')$ to $(U_y,
P_y, V_y, G_y = H_y K_y)$ lifts to a reduction of $(\widetilde U', \widetilde P',
\widetilde V',\break \widetilde G' = \widetilde H' \widetilde K')$ to $(\widetilde U_y,
\widetilde P_y, \widetilde V_y, \widetilde G_y = \widetilde H_y \widetilde K_y)$.

In order to do that, we merge diagram \ref{map3/dia}$'$ with the analogous diagram
\ref{map3/dia}$_y$ that gives the lifting of $(U_y, P_y, V_y, G_y = H_y K_y)$, by
identifying the corresponding spaces and maps which are the same in the two diagrams
(these are $U' = U_y$, $V' = V_y$, $P'/G' = P_y/G_y$, $\widetilde U' = \widetilde U_y$,
$\widetilde V' = \widetilde V_y$ and the maps between them). In this way we get a unique
commutative diagram, which can be completed, still keeping the commutativity, by the
addition of the following maps in turn: the reduction map $P' \to P_y$; a PL map $Q' \to
Q_y$, which can be constructed by using the universal property of the connected pullback
$q_y: Q_y \to U_y$ (cf. Proposition \ref{connectedpullback/thm}); a PL map $\widetilde P'
\to \widetilde P_y$, whose existence derives from the minimality of the regularization
$r_y: \widetilde P_y \to P_y/G_y$. The last map we added clearly gives the wanted
reduction.

At this point, we know that the lifted charts generate a branchfold structure on $X$
making $f: X \to Y$ into a branchfold covering ($f$ is complete, being a non-degenerate PL
map).

The unicity of such branchfold structure, immediately follows once we prove that for any
branchfold charts $(U, P, V, G = HK)$ on $X$ and $(U', P', V', G' = H'K')$ on $Y$, which
are related by diagram \ref{map1/dia} with the additional assumption that $\psi$ is an
isomorphism, there is a chart equivalence between $(U, P, V, G = HK)$ and a suitable
lifting $(\widetilde V', \widetilde P', \widetilde U', \widetilde G' = \widetilde H'
\widetilde K')$ of $(U', P', V', G' = H'K')$.

So, let us look at diagrams \ref{map1/dia} and \ref{map3/dia}$'$. By Proposition
\ref{restaction/thm}, we have that $U$ is a con\-nected component of $f^{-1}(U')$.
Moreover, up to chart isomorphism, assuming that $\psi$ is an isomorphism is the same as
assuming $V = V'$ and $\psi = \id$. Now, let $q' = p_{H'} \circ q'_1 = f_| \circ q'_2: Q'
\to U'$ be the connected pullback of $p_{H'}$ and $f_|$. Then, the universal property of
this pullback gives us a branched covering $s'': P \to Q'$ such that $q_1 \circ s'' =
\phi$ and $q_2 \circ s'' = p_H$. Finally, we put $r'' = \pi_{G'} \circ \phi$. In this
way\break we get a commutative diagram which is like the diagram \ref{map3/dia}$'$
realizing the lifting of $(U', P', V', G' = H'K')$ with $\widetilde U' = U$. Actually, the
differences are that we have $P$, $p_H$, $p_K$, $s''$ and $r''$ respectively in place of
$\widetilde P'$, $p_{\widetilde H'}$, $p_{\widetilde K'}$, $s'$ and $r'$, and that $r''$
is not necessarily regular. However, by the same argument used in diagram \ref{map4/dia},
we can construct a new branchfold chart $(\widehat U, \widehat P, \widehat V, \widehat G =
\widehat H \widehat K)$ starting from the minimal regularization $\widehat r = r'' \circ
\widehat s: \widehat P \to P/G$ of $r'': P \to P'/G'$. This new chart dominates $(U, P, V,
G = HK)$ through $\widehat s: \widehat P \to P$. On the other hand, it dominates also
$(\widetilde V', \widetilde P', \widetilde U', \widetilde G' = \widetilde H' \widetilde
K')$ through a map $t: \widehat P \to \widetilde P'$, which exists by the minimality of
$r'$ as a regularization of $\pi_{G'} \circ q'_1$. Thus we have the desired equivalence.
\end{proof}

\begin{proposition}\label{quotientbranchfold/thm}
Let $X$ be a branchfold and $f: X \to Y$ be a regular branch\-ed covering, whose deck
transformations are branchfold isomorphisms which preserve local orientations at their
fixed points. Then the branchfold structure of $X$ induces a unique branchfold structure
on $Y$ making $f$ into a PL branchfold covering.
\end{proposition}

\begin{proof}
Given $x \in X$, let $L$ be the group of the deck transformations of $f$ that fix $x$ and
let $(U, P, V, G = HK)$ be any $L$-invariant conical chart of $X$ centered at $x$, such
that the induced action of $L$ on $U$ is conical and $t(U) \cap U = \emptyset$ for any
deck transformation $t$ of $f$ not in $L$. Since $L$ is finite and acts on $X$ by
branchfold isomorphisms, such a chart always exists, an example being the local model of
$X$ at $x$, based on any $L$-invariant sufficiently small open star centered at $x$.

Then, the restriction $f_|: U \to \widebar U$ is a conical regular branched covering of an
open conical neighborhood $\widebar U = f(U)$ of $y = f(x)$ in $Y$, and we can think of it
as the canonical projection $\pi_L: U \to U/L \cong \widebar U$, where the action of $L$
on $U$ is given by restriction and preserves orientations by hypothesis.

As a consequence of the $L$-invariance of the chart $(U, P, V, G = HK)$, for any $l \in L$
the restriction $l: U \to U$ lifts to some PL homeomorphism $\lambda: P \to P$. Let
$\widebar H$ be the group of all such liftings when $l$ varies in $L$. Then, $H \nsubgr
\widebar H$, since $H$ consists of the liftings of the identity. Moreover, using the
$L$-invariance of $U$ once again, we see that $K$, and hence also $G$, is a normal
subgroup in the group $\widebar G$ of orientation preserving PL homeomorphisms of $P$
generated by $\widebar H \cup K$. So, we can write $\widebar G = \widebar H K$.

By Propositions \ref{goodaction/thm} and \ref{goodaction-pm/thm}, all maps in the
following commutative diagram are regular branched coverings, possibly except $p$ and
$\widebar p$ (which are however branched coverings). Then $(\widebar U, P, V, \widebar G =
\widebar H K)$ is a conical branchfold chart on $Y$ centered at $y$, which is related to
the original chart $(U, P, V, G = HK)$ by a commutative diagram like diagram
\ref{map1/dia} (with $\phi = \id_P$ and $\psi = \id_V$).

\begin{Diagram}[htb]{map5/dia}
\diag{}{map5.eps}{}
\end{Diagram}

We observe that the branchfold chart $(\widebar U, P, V, \widebar G = \widebar H K)$ turns
out to be well defined up to chart equivalence, being independent on choices involved in
its construction. In particular, if $(U, P', V, G' = H'K')$ is an $L$-invariant conical
chart centered at $x$ which dominates $(U, P, V, G = HK)$, then the resulting chart
$(\widebar U, P', V, \widebar G' = \widebar H' K')$ dominates $(\widebar U, P, V, \widebar
G = \widebar H K)$. Moreover, since the deck transformations of $f$ are branchfold
isomorphisms it does not matter what $x$ we choose in $f^{-1}(y)$.

On the other hand, the above construction also preserves chart restrictions. Namely, if we
start from any $x' \in U$ and denote by $L'$ the group of the deck transformations of $f$
that fix $x'$ and by $(U', P', V', G' = H'K')$ any $L'$-invariant conical restriction of
$(U, P, V, G = HK)$ centered at $x'$, then we have that the resulting chart $(\widebar U,
P', V, \widebar G' = \widebar H' K')$ is a restriction of $(\widebar U, P, V, \widebar G =
\widebar H K)$.

Thus, when $y$ varies in $Y$ the branchfold charts $(\widebar U, P, V, \widebar G =
\widebar H K)$ constructed above form a branchfold atlas on $Y$. This gives us the desired
branchfold structure on $Y$ making $f$ into a branchfold covering ($f$ is complete, being
a non-degenerate PL map).

In order to prove the unicity of such branchfold structure, consider any chart $(U, P, V,
G = HK)$ of the branchfold $X$ and any branchfold chart $(U', P', V', G' = H'K')$ on the
polyhedron $Y$, which are related as in diagram \ref{map1/dia}. Given $x \in U$ and $y =
f(x) \in U'$, up to chart restricion we can assume that both the charts are conical,
respectively centered at $x$ and $y$, in such the way that $U' \cong U/L$, where $L$ is
the group of the deck transformations of $f$ which fix $x$. Then, let $r: P'' \to P'/G'$
be the minimal regularization of the branched covering $\pi_{G'} \circ \phi: P \to P'/G'$
and let $p_{H''}: P'' \to U$ and $p_{K''}: P'' \to V$ be defined by composition, in such a
way that the following diagram commutes.

\begin{Diagram}[htb]{map6/dia}\vskip-6pt
\diag{}{map6.eps}{}\vskip-6pt
\end{Diagram}

By putting $G'' = H'' K''$, where $H''$ and $K''$ are the subgroups of the group of deck
transformations of $r$ corresponding to $p_{H''}$ and $p_{K''}$ respectively, we get a new
chart $(U, P'', V, G'' = H'' K'')$ of $X$, which dominates the original one. On the other
hand, this chart is $L$-invariant and the above construction applied to it gives us a new
chart $(\widebar U, P'', V, \widebar G'' = \widebar H'' K'')$ on $Y$, which dominates
$(U', P', V', G' = H'K')$ through the map $\phi \circ s: P'' \to P'$ (notice that $V \cong
V'$, being $\psi$ an isomorphism, and $\widebar U = U/L \cong U'$). Hence, we are done
thanks to the arbitrary choice of $y$ and $(U', P', V', G' = H'K')$.
\end{proof}

It is worth emphasizing that Propositions \ref{coverbranchfold/thm} and
\ref{quotientbranchfold/thm} apply in the order to the branched coverings $p_K$ and $p_H$
of diagram \ref{chart2/dia} for any local chart $(U, P, V, G = HK)$, to determine the
branchfold structure of $U$ starting from the smooth one of $V$ (see comments to diagram
\ref{good1/dia} below).

\medskip

Finally, thanks to Propositions \ref{branchfoldcover/thm} and \ref{coverbranchfold/thm},
we can consider pullbacks and regularizations in the category of PL branchfold coverings,
by endowing the new space arising in those constructions of the right branchfold structure
which makes all the new maps into PL branchfold coverings. From now on, we will do that
without any further comment.

\subsection*{Universal branchfold covering\label{universal/sub}}

Let $X$ be a connected branchfold. We recall that the singular locus $\Sigma X$ is a good
subpolyhedron of $X$ and hence $X - \Sigma X$ is connected by Proposition
\ref{goodconnect/thm}.

Given any point $x \in X$, let $(U_x, P, V, G = H K)$ be a conical chart of $X$ centered
at $x$. Referring to the chart diagram \ref{chart2/dia}, we consider the maps $p_|: U_x -
p_H(S_G) \to P/G - B_G$ and $p_{G/K|}: V - p_K(S_G) \to P/G - B_G$, where $S_G$ is the
singular set of the action of $G$ on $P$ and $B_G$ is the branch set of the canonical
projection $\pi_G$. These are connected ordinary coverings by Proposition
\ref{goodconnect/thm} and \ref{covergood/thm}. On the other hand, $p_H(S_G) \subset U_x$
is a subpolyhedron of codimension $\geq 2$ containing the singular set $\Sigma U_x =
\Sigma X \cap U_x$. As a consequence, since $U_x -\Sigma U_x$ is a manifold, by
transversality we have that the homomorphism $i_*: \pi_1(U_x - p_H(S_G)) \to \pi_1(U_x -
\Sigma U_x)$ induced by the inclusion is surjective.

We put $\Gamma_x = i_*((p_{|*})^{-1} (\Im p_{G/K|*})) \nsubgr \pi_1(U_x - \Sigma U_x)$.
From a different point of view, $\Gamma_x$ is the normal subgroup of $\pi_1(U_x - \Sigma
U_x)$ whose elements are represented by the loops $\omega$ in $U_x - p_H(S_G)$ with the
following property: if $\widetilde \omega$ is any path lifting $\omega$ to $P - S_G$
through the ordinary covering $p_{H|}: P - S_G \to U_x - p_H(S_G)$, then its projection
$p_K \circ \widetilde \omega$ in $V - p_K(S_G)$ is a loop. In particular, $\Gamma_x$
contains $i_*(\Im p_{H|*})$ and hence it has finite index in $\pi_1(U_x - \Sigma U_x)$,
being $p_{H|}$ a finite covering.

We notice that the group $\Gamma_x$ depend only on the local model of $X$ at $x$ and not
on the particular conical chart $(U_x, P, V, G = H K)$ we started with. In fact, any such
chart dominates $(U_x, P_x, V_x, G_x = H_x K_x)$, hence it gives raise to the same group
$\Gamma_x$, as it can be easily seen by looking at diagram \ref{domination/dia}.

Once fixed base points $* \in X$ and $*_x \in U_x$ for any $x \in X$, we choose a path
$\alpha_x$ in $X - \Sigma X$ connecting $*$ to $*_x$ and denote by $h_x: \pi_1(U_x -
\Sigma U_x) \to \pi_1(X - \Sigma X)$ the natural homomorphism induced by the map $\omega
\mapsto \widebar \alpha_x \omega \alpha_x$.

At this point, we define $\Gamma_X$ to be the smallest normal subgroup of $\pi_1(X -
\Sigma X)$ containing all the $h_x(\Gamma_x)$'s, which is clearly independent on the
choice of the $\alpha_x$'s. Of course, here it is enough to consider the groups
$h_x(\Gamma_x)$ with $x \in \Sigma X$, being the other ones trivial. Actually, due to the
fact that the local model is constant on the connected components of the strata of
$\Sigma_X$, it would suffice to let $x$ vary on a set of representatives of such
components.

\begin{definition}\label{chargroup/def}
We call $\Gamma_X \nsubgr \pi_1(X - \Sigma X)$ the {\sl characteristic group} of $X$, and
$\Gamma_x \nsubgr \pi_1(U_x - \Sigma U_x)$ the {\sl local characterictic group} of $X$ at
$x \in X$.
\end{definition}

Now, we consider the ordinary regular covering $r: R \to X - \Sigma X$ corresponding to
the normal subgroup $\Gamma_X \nsubgr \pi_1(X - \Sigma X)$, meaning that $r_*(\pi_1(R)) =
\Gamma_X$.

First of all, we observe that, for each $x \in X$ and each connected component $C$ of
$r^{-1}(U_x - \Sigma U_x)$, the group $r_{|*}(\pi_1(C)) \subgr \pi_1(U_x - \Sigma U_x)$
contains the local characteristic group $\Gamma_x$. In fact, choose base points
$\widetilde *$ and $\widetilde *_x$ respectively for $R$ and $C$, such that $r(\widetilde
*) = *$ and $r(\widetilde *_x) = *_x$, and choose $\alpha_x = r \circ \beta_x$ with
$\beta_x$ a path in $R$ from $\widetilde *$ to $\widetilde *_x$. Then, by definition of
$\Gamma_X$, for every loop $\omega$ in $U_x - \Sigma U_x$ representing an element of
$\Gamma_x$, the loop $\widebar \alpha_x \omega \alpha_x$ lifts to a loop in $R$ with
respect to the covering $r$. This lifting has to be of the form $\widebar \beta_x
\widetilde \omega \beta_x$, where $\widetilde \omega$ is a loop in $C$ such that $r_|
\circ \widetilde \omega = \omega$.

In the light of the above observation, since the local characteristic group $\Gamma_x$ has
finite index in $\pi_1(U_x - \Sigma U_x)$, a fortiori the same is true for the group
$r_{|*}(\pi_1(C))$.

At this point, we can easily conclude that $r$ satisfies the monodromy hypothesis of
Proposition \ref{completion/thm}, therefore it can be completed to a regular branched
covering $u_X : \widetilde X \to X$. According to Proposition \ref{branchfoldcover/thm},
we think $\widetilde X$ as a branchfold endowed with the unique structure making $u_X$
into a PL branchfold covering.

\begin{definition}\label{univbranchcover/def}
We call $u_X: \widetilde X \to X$ the {\sl universal branchfold covering} of the connected
branchfold $X$. Moreover, we define the {\sl branchfold fundamental group} as the group
$\pi_1^b(X) = \pi_1(X - \Sigma X)/\Gamma_X$ of deck transformations of $u_X$.
\end{definition}

We emphasize that universal branchfold coverings are PL branchfold coverings and they are
natural in the sense specified by the following proposition.

\begin{proposition}\label{naturalcover/thm}
For any PL branchfold covering $f: X \to Y$, there exists a PL branchfold covering
$\widetilde f: \widetilde X \to \widetilde Y$ such that $f \circ u_X = u_Y \circ
\widetilde f$, and this is unique up to deck transformations of $u_X$ and $u_Y$.
\end{proposition}

\begin{proof}
Let us consider $S = \Sigma Y \cup f(\Sigma X) \subset Y$ and $f^{-1}(S) \subset X$, which
are good subpolyhedra by Propositions \ref{covergood/thm} and \ref{branchfoldcover/thm}.
The restriction $f_|: X - f^{-1}(S) \to Y - S$, as well as the restriction of $u_X$ over
$X - f^{-1}(S)$ and that of $u_Y$ over $Y - S$, are all ordinary coverings.

Any loop in $X - f^{-1}(S)$ representing an element of $\Im u_{X|*}$ is a product of loops
of the form $\widebar \alpha_x \omega_x \alpha_x$, where $\omega_x$ is a loop in $U_x -
f^{-1}(S)$ whose lifting to $P_x$ through $p_{H_x}$ projects to a loop in $V_x$ through
$p_{K_x}$. We put $y = f(x)$, $\alpha_y = f \circ \alpha_x$ and $\omega_y = f \circ
\omega_x$. By inspection on the chart diagram \ref{map1/dia} for to the conical
restriction $f_|: U_x \to U_y$, we see that $\omega_y$ is a loop in $U_y - S$ whose
lifting to $P_y$ through $p_{H_y}$ projects to a loop in $V_y$ through $p_{K_y}$. Hence,
$\widebar \alpha_y \omega_y \alpha_y$ represents an element of $\Im u_{Y|*}$.

Therefore, the homomorphism $f_{|*}: \pi_1(X - f^{-1}(S)) \to \pi_1(Y - S)$ sends $\Im
u_{X|*}$ into $\Im u_{Y|*}$ and hence the restriction $f_|: X - f^{-1}(S) \to Y - S$ can
be lifted to an ordinary covering $\widetilde f_|: \widetilde X - u_X^{-1}(f^{-1}(S)) \to
\widetilde Y - u_Y^{-1}(S)$ through the ordinary coverings $u_{X|}$ and $u_{Y|}$.

Then, the wanted lifting $\widetilde f$ can be obtained as the completion of $\widetilde
f_|$. Proposition \ref{coverbranchfold/thm}, locally applied to $f \circ u_X = u_Y \circ
\widetilde f$, tells us that $\widetilde f$ is actually a branchfold covering. On the
other hand, the uniqueness of such a lifting $\widetilde f$ up to deck transformations
immediately follows from that of $\widetilde f_|$.
\end{proof}

Next proposition characterizes universal branchfold coverings in terms of their universal
property with respect to a certain class of PL branchfold coverings.

\begin{proposition}\label{universalcover/thm}
Let $X$ be a connected branchfold. The universal branchfold covering $u_X: \widetilde X
\to X$ satisfies the following property: for any $\widetilde x \in \widetilde X$, the
local model $(U_x, P_x, V_x, G_x = H_x K_x)$ of $X$ at $x = u_X(\widetilde x)$ lifts to a
conical (possibly non-reduced) branchfold chart $(U_{\widetilde x}, P_x, V_x, \widebar G_x
= \widebar H_x K_x)$ of $\widetilde X$ centered at $\widetilde x$, such that $p_{H_x} =
u_{X|} \circ p_{\widebar H_x}$, with $\widebar H_x \subgr H_x$ and $\widebar H_x \cap K_x
= H_x \cap K_x$. Moreover, for any PL branchfold covering $f: \widehat X \to X$ satisfying
the same property, there exists a PL branchfold covering $g: \widetilde X \to \widehat X$
such that $u_X = f \circ g$.
\end{proposition}

\begin{proof}
Let $x$ and $\widetilde x$ be as in the statement. When defining $u_X$, we proved that
$i_*(\Im p_{H_x|*}) \subgr \Gamma_x \subgr \Im u_{X|*} \subgr \pi_1(U_x - \Sigma U_x)$,
where $u_{X|}: U_{\widetilde x} - \Sigma U_{\widetilde x} \to U_x - \Sigma U_x$ is the
restriction of $u_X$ to the non-singular part of the open star $U_{\widetilde x}$ of
$\widetilde X$ at $\widetilde x$. By the standard theory of ordinary coverings, there
exists a lifting $r_x: P_x - S_{G_x} \to U_{\widetilde x} - \Sigma U_{\widetilde x}$ of $i
\circ p_{H_x|}: P_x - S_{G_x} \to U_x - \Sigma U_x$ through $u_{X|}$. The restriction
$r_{x|} : P_x - S_{G_x} \to r_x(P_x - S_{G_x})$ is a regular ordinary covering, which can
be completed to a regular branched covering $p_{\widebar H_x}: P_x \to U_{\widetilde x}$,
whose deck transformations form a subgroup $\widebar H_x \subgr H_x$. By considering the
special case of diagram \ref{map3/dia}, with $\widetilde P \cong Q \cong P = P_x$,
$\widetilde V = V = V_x$, $\widetilde U = U_{\widetilde x}$, $U = U_x$, $s \cong q_1 \cong
\phi = \id_{P_x}$, $q_2 \cong p_{\widetilde H} = p_{\widebar H_x}$, $p_{\widetilde K} =
p_K =p_{K_x}$ and $p_H = p_{H_x}$, we immediately see that $(U_{\widetilde x}, P_x, V_x,
\widebar G_x = \widebar H_x K_x)$ is a conical branchfold chart of $\widetilde X$ centered
at $\widetilde x$. In order to prove the equality $\widebar H_x \cap K_x = H_x \cap K_x$,
let us consider any $g \in H_x \cap K_x$ and let $\alpha$ be any path in $P_x - S_{G_x}$
from the base point $*_x$ to its image $g(*_x)$. Then, $\omega = p_{H_x} \circ \alpha$ is
a loop in $U_x - \Sigma U_x$ and $p_{K_x} \circ \alpha$ is a loop in $V_x$, hence $\omega
\in \Gamma_x$. The inclusion $\Gamma_x \subgr \Im u_{X|*}$ implies that $\widetilde \omega
= p_{\widebar H_x} \circ \alpha$ is a loop in $U_{\widetilde x} - \Sigma U_{\widetilde
x}$. This in turn allows us to conclude that $g \in \widebar H_x$.

To prove the second part of the proposition, let us consider any PL branchfold covering
$f: \widehat X \to X$ as in the statement. It follows from the stated property that the
restriction $f_|: \widehat X - f^{-1}(\Sigma X) \to X - \Sigma X$ is an ordinary covering,
in such a way that $f$ can be thought as the completion of it.

Let $\widebar \alpha_x \omega_x \alpha_x$ be any generator of $\Gamma_X$, where $\alpha_x$
is a path from the base point $*$ of $X$ to the base point $*_x$ of $U_x$ and $\omega_x$
is a loop in $U_x - \Sigma U_x$. For any lifting $\widetilde \omega_x$ of $\omega_x$ to
$P_x$ through $p_{H_x}$, there exists an element $g \in H_x \cap K_x$ such that
$\widetilde \omega_x(1) = g(\widetilde \omega_x(0))$.\break Now, let $\widehat \alpha_x$
be the lifting of $\alpha_x$ through $f_|$ starting from the base point of $\widehat X -
f^{-1}(\Sigma X)$ and denote by $\widehat x \in f^{-1}(x)$ be the point such that
$\widehat \alpha_x(1) \in U_{\widehat x}$. The equality $\widebar H_x \cap K_x = H_x \cap
K_x$ for the conical restriction of $f$ at $\widehat x$, implies that $\widehat \omega_x =
p_{\widebar H_x} \circ \widetilde \omega_x$ is a loop in $U_{\widehat x}$. Thus, $\widebar
\alpha_x \omega_x \alpha_x$ belongs to $f_{|*}(\pi_1(X - f^{-1}(\Sigma X)))$.

So, we have proved the inclusion $\Gamma_X \subgr \Im f_{|*}$. Recalling from the
construction of $u_X$ that $\Gamma_X$ coincides with $\Im u_{X|*}$ for the restriction
$u_{X|}: \widetilde X - u_X^{-1}(\Sigma X) \to X - \Sigma X$, we can conclude that such
restriction lifts through $f_|$ to an ordinary covering $s: \widetilde X - u_X^{-1}(\Sigma
X) \to \widehat X - f^{-1}(\Sigma X)$. Then, the desired factorization $u_X = f \circ g$
can be obtained by defining $g$ as the completion of $s$ over $\widehat X$.
\end{proof}

Some comments concerning the statement of Proposition \ref{universalcover/thm} is in
order. First of all, we notice that the map $\phi$ in diagram \ref{map1/dia} can be always
assumed to be a branchfold isomorphism for any branchfold covering (to see this, consider
the dominations of the two involved charts induced by the regularization of the branched
covering $\pi_{G'} \circ \phi$). \label{universalcover/rem}
In other words, we can always assume $P = P'$, $K = K'$ and $H \subgr H'$. The significant
point in the local property stated for $u_X$ by Proposition \ref{universalcover/thm} is
that in this case the same assumption can be made keeping the local chart of the range to
be minimal (no domination is required for it). Moreover, the equality $\widebar H_x \cap
K_x = H_x \cap K_x$ tells us that such a local property is invariant up to
dominations/reductions of the local chart of the range.

\medskip

We observe that, as an immediate consequence of Proposition \ref{universalcover/thm},
$u_X$ is an ordinary covering (actually it is the universal covering of $X$, as we will
see shortly) when $X$ is a pure branchfold. Obviously, in this case also $\widetilde X$ is
a pure branchfold.

In general $\widetilde X$ is not a pure branchfold, but in some sense we can say that it
is as pure as possible among the simply connected PL branchfold coverings of $X$ whose
branch set is contained in $\Sigma X$. This fact can be formalized as a universal property
of $u_X$ with respect to the pure branchfold coverings of $X$.

\begin{proposition}\label{mostpurecover/thm}
Let $X$ be a connected branchfold. The universal covering space $\widetilde X$ is simply
connected. Moreover, for any PL branchfold covering $f: \widehat X \to X$ with $\widehat
X$ a simply connected pure branchfold there exists a PL branchfold covering $g: \widehat X
\to \widetilde X$ such that $f = u_X \circ g$.
\end{proposition}

\begin{proof}
The simply connectedness of $\widetilde X$ follows from Proposition
\ref{universalcover/thm}. In fact, if $\widetilde X$ would not be simply connected, the
universal property of $u_X$ stated by that proposition would be contradicted by the
composition $f = u_X \circ v : \widehat X \to X$, where $v: \widehat X \to \widetilde X$
is the ordinary universal covering of $\widetilde X$.

Now, by applying Proposition \ref{naturalcover/thm} to the PL branched covering $f:
\widehat X \to X$ in the statement, we get the factorizing map $g$ as the lifting
$\widetilde f$ of $f$ to the universal branchfold coverings (up to branchfold
isomorphism), being the universal branchfold covering of $\widehat X$ a branchfold
isomorphism by the above observation.
\end{proof}

Finally, as another consequence of Proposition \ref{universalcover/thm}, when $X$ is a
connected orbifold the map $u_X: \widetilde X \to X$ coincides with the universal orbifold
covering and $\pi_1^b(X)$ coincides with the fundamental orbifold group $\pi_1^o(X)$.

Furthermore, in this case the characteristic group $\Gamma_X$ is normally generated by the
powers $\mu_C^{i_C}$, where $C$ varies among the connected components of the
$(m-2)$-stratum of $\Sigma X$, $\mu_C$ is any meridian around $C$ and $i_C$ is the index
of $C$. In fact, in this case $\Gamma_x = \Im p_{H_x|*}$ and the total space $P_x$ of the
local model at $x$ is a disk for any $x \in X$. Then, we can easily conclude that any loop
$\widebar \alpha_x \omega \alpha_x \in h_x(\Gamma_x)$ is homotopic to a product of powers
$\mu_C^{i_C}$ as above, being any loop $\widetilde \omega$ in $P_x - p_{H_x}^{-1} (\Sigma
U_x)$ homotopic to the composition of meridians around $p_{H_x}^{-1} (\Sigma U_x)$.

\subsection*{Good branchfolds\label{goodbranchfolds/sub}}

Thanks to Propositions \ref{coverbranchfold/thm} and \ref{quotientbranchfold/thm}, we can
extend to branchfolds the construction of locally orientable good orbifolds as global
quotients of smooth manifolds by properly discontinuous smooth actions.

Namely, let $P$ be a connected polyhedron and $G = H K$ be a group acting on $P$, such
that $K \nsubgr G$ is normal subgroup, $M = P/K$ is a smooth $m$-manifold and the induced
action of $G/K$ on $M$ is smooth and preserves local orientations, meaning that it is
given by diffeomorphisms which preserve local orientations at their fixed points.

Proposition \ref{coverbranchfold/thm} allows us to endow $P$ with a pure $m$-branchfold
structure, which is uniquely determined by the property of making the canonical projection
$\pi_K: P \to M$ into a PL branchfold covering. The action of $H$ on $P$ leaves this
branchfold structure invariant and preserves local orientations, hence we can apply
Proposition \ref{quotientbranchfold/thm} to obtain a unique $m$-branchfold structure on $X
= P/H$, which makes the canonical projection $\pi_H: P \to X$ into a PL branchfold
covering. On the other hand, the quotient space $P/G = M /(G/K)$ is a locally orientable
good $m$-orbifold and $p: X \to P/G$ is a PL branchfold covering. Hence, we have the
commutative diagram \ref{good1/dia}, consisting of PL branchfold coverings.

\begin{Diagram}[htb]{good1/dia}
\diag{}{good1.eps}{}
\end{Diagram}

Actually, starting from any PL branchfold covering $c: X \to O$ of a good orbifold $O$, we
can produce a diagram like \ref{good1/dia}, in the following way (see diagram
\ref{good2/dia}). Let $r: M \to O$ be a regular orbifold covering with $M$ a smooth
manifold, which always exists being any manifold covering of $O$ virtually regular (cf.
\cite{LS89}). Then, consider the connected pullback $q: Q \to O$ of $r$ and $c$, and the
minimal regularization $q \circ s: P \to O$ of $q$. Denote by $L$ the group of deck
transformations of such regularization, and define $H \subgr L$ as the subgroup
corresponding to the regular covering $\pi_H = q_1 \circ s$, and $K \nsubgr L$ as the
normal subgroup corresponding to the regular covering $\pi_K = q_2 \circ s$. Finally, put
$G = HK \subgr L$ and complete the diagram with the quotient $P/G$ and the coverings $p$,
$t$ and $\pi_{G/K}$. Of course, the new spaces $Q$, $P$ and $P/G$ are endowed with the
branchfold structures which make all the maps into PL branchfold coverings. In particular,
$P/G = M/(G/K)$ turns out to be an orbifold.

\begin{Diagram}[htb]{good2/dia}
\diag{}{good2.eps}{}
\end{Diagram}

We remark that, if $c: X \to O$ is a finite branchfold covering of a very good orbifold,
then all the coverings in the diagram \ref{good2/dia} can be assumed to be finite. Hence,
in this case we get a diagram like \ref{good1/dia}, consisting of finite coverings.

\medskip

In the light of the above considerations, the notions of good and very good branchfold
introduced by the next definition appear as natural extensions of the usual notions of
good and very good orbifold.

\begin{definition}\label{goodbranchfold/def}
A connected branchfold $X$ is called a {\sl good branchfold} when it admits PL branchfold
covering $c: X \to O$ onto a good orbifold $O$. If in addition $c$ has finite degree and
$O$ is a very good orbifold, then we call $X$ a {\sl very good branchfold}.
\end{definition}

This notions of good and very good coincide with the usual ones when referred to
orbifolds. In fact, it is clear from the definition above that a good (resp. very good)
orbifold $X$ is also good (resp. very good) as a branchfold. For the viceversa, it is
enough to observe that the branchfold $P$ in diagram \ref{good2/dia} is a smooth manifold
(that is $\Sigma P = \emptyset$) when $X$ is an orbifold. To see this, we first consider
the universal (orbifold) covering $u_O: \widetilde O \to O$. If $X$ is an orbifold, both
$p$ and $c$ are orbifold cover\-ings, hence $u_O$ lifts to orbifold coverings $\widetilde
O \to M$ and $\widetilde O \to X$. Then, by the universal properties of the pullback and
of the minimal regularization, we get a branchfold covering $\widetilde O \to P$. So, we
can conclude that $P$ is a manifold, being at the same time a pure branchfold (as a
covering of $M$) and an orbifold (as a quotient of $\widetilde O$).

\medskip

It is worth noticing that the universal branchfold covering $\widetilde X$ of a good
branchfold $X$ is a pure branchfold. The reason is that the covering $c: X \to O$ in
Definition \ref {goodbranchfold/def} lifts to a branchfold covering $\widetilde c:
\widetilde X \to \widetilde O$ with $\widetilde O$ a smooth manifold.

\medskip

Finally, by the very definition of branchfold chart, we have that any branchfold is
locally very good, being diagram \ref{chart2/dia} at page \pageref{chart2/dia} a special
case of diagram \ref{good1/dia}, with $U$ in place of $X$ and all the coverings of finite
degree.

\section{Geometric branchfolds\label{geometry/sec}}

In this section we introduce the notion of geometric structure on a branchfold and extend
to geometric branchfolds the well known goodness theorem for geometric orbifolds
(\cite{Thu97}, \cite{MM91}). Then, we study the relation between branchfolds endowed with
geometries modelled on constant curvature spaces and conifolds. In particular, we apply
the geometric goodness theorem to establish what conifolds can be thought as branchfolds.

\medskip

\label{geometry/def}
By a {\sl geometry} we mean a pair $(\G,\M)$, where $\M$ is a simply connected smooth
manifold and $\G$ is a transitive group of diffeomorphisms of $\M$ which is {\sl locally
effective}, meaning that if $g \in \G$ and $g_{|U} = \id_U$ for some open subset $U
\subset \M$ then $g = \id_{\M}$.

\medskip

Following \cite{MM91}, here we use local effectiveness in place of the stronger assumption
of analyticity of the action made by Thurston in his notes \cite{Thu79}.

\medskip

\label{pseudogroup/def}
Given a geometry $(\G,\M)$, we denote by $\widehat \G$ the pseudo-group of all the
diffeomorphisms between non-empty open subsets of $\M$ obtained as restrictions of
elements of $\G$, that is $\widehat \G = \{g_{|}: U \to g(U) \;|\; g \in \G \text{ and } U
\neq \emptyset \text{ open in \M}\}$. Then, the local effectiveness can be reformulated by
saying that each element of $\widehat \G$ extends to an unique elements of $\G$.

\medskip

\label{admissible/def}
A polyhedron $P \subset \M$ is called {\sl $\G$-admissible} if it admits a stratification,
whose strata locally coincide with fixed point sets $\Fix G$, where $G \subset \widehat
\G$ is any finitely generated group acting on an open subset of $\M$.

\medskip

In particular, the singular set $S_G$ of any good action of a finite group $G \subset
\widehat \G$ on an open subset $V \subset \M$ is $\G$-admissible in $\M$. In fact, for
every $x \in S_G$, the fixed point set $\Fix G_x$ of the stabilizer of $x$ is a smooth
submanifold of $V$ contained in $S_G$, and we can define a stratification of $S_G$, by
putting $(S_G)_i = \{x \in S_G \;|\; \dim \Fix G_x \leq i\}$. Such a stratification
satisfies the property required for the $\G$-admissibility of $S_G$.

In the case of a Riemannian geometry, that is when $\M$ is a Riemannian manifold and $\G$
is the group of isometries of $\M$, if $P \subset \M$ is a $\G$-admissible then it can be
stratified by totally geodesic submanifolds of $\M$ (cf. \cite{K72}). The viceversa holds
when $\M$ has constant curvature. Our definition of admissibility is a tentative
reformulation of this metric property in the abstract context of $(\G,\M)$ geometries.

\medskip

\begin{definition}\label{GMchart/def}
A branchfold chart $(U, P, V, G = HK)$ is called a {\sl $(\G,\M)$-chart} when $V$ is
identified with an open subset $V \subset \M$ in such a way that:
\begin{itemize}\parskip0pt\itemsep0pt
\item[(1)] \vskip-\lastskip
the branch set of ${p_K}$ is a $\G$-admissible subpolyhedron of $\M$;
\item[(2)] 
the induced action of $G/K$ on $V$ is given by elements of $\widehat \G$.
\end{itemize}\vskip-\lastskip\vskip-\baselineskip
\end{definition}

Property 1 is aimed to impose a reasonable restriction on the branched covering $\pi_K$,
otherwise any branched covering of a connected open subset $V \subset \M$ would be a
$(\G,\M)$-chart for any geometry $(\G,\M)$. How we will see, this property works well in
the case of constant curvature Riemannian geometries, but we are not sure that it is the
right property to be required in general.

On the other hand, property 2 is a quite natural extension of the analogous one usually
required for an orbifold chart to be geometric. Actually, it implies that the definition
of $(\G,\M)$-chart reduces to the standard one in the orbifold case.

\medskip

Two $(\G,\M)$-charts are called {\sl $(\G,\M)$-isomorphic} (resp. {\sl strongly
$(\G,\M)$-isomor\-phic}) when they are isomorphic (resp. strongly isomorphic) and the
diffeomorphism $V_1 \cong V_2$ in the corresponding definition at page
\pageref{isochart/def} belongs to $\widehat \G$.

Now, any restriction of a $(\G,\M)$-chart is still a $(\G,\M)$-chart. In fact, referring
to Definition \ref{subchart/def}, we have $B_{p_{K'}} = B_{p_K} \cap V'$ and this is
$\G$-admissible, being $\G$-admissibility local property. Concerning the action of
$G'/K'$ on $V'$, this is given by elements of $\widehat \G$, being a restriction of the
action of $G/K$ on $V$.

Analogously, any reduction of a $(\G,\M)$-chart is still a $(\G,\M)$-chart. In this case,
referring to Definition \ref{chartdomin/def}, $B_{p_K}$ is a subpolyhedron of $B_{p_{K'}}$
and it can be obtained by deleting some connected components from the strata of the
stratification giving the admissibility of $B_{p_{K'}}$, hence it is admissible itself.
At the same time, property 2 is trivially preserved by reductions.

On the contrary, it is not difficult to see that a chart dominating a $(\G,\M)$-chart is
not necessarily a $(\G,\M)$-chart. In particular, the common dominating chart giving the
equivalence between two $(\G,\M)$-charts (cf. Definition \ref{chartequiv/def}) is not
necessarily a $(\G,\M)$-chart.

\begin{definition}\label{GMatlas/def}
A branchfold atlas $\cal A = \{(U_i, P_i, V_i, G_i = H_i K_i)\}_{i\,\in\,I}$ on $X$ is
called a {\sl $(\G,\M)$-atlas}, if it consists of $(\G,\M)$-charts which satisfy the
compatibility condition in Definition \ref{atlas/def}, where the strong isomorphisms are
required to be strong $(\G,\M)$-isomorphism (in other words, the diffeomorphisms $V'_i
\cong V' \cong V'_j$ in diagram \ref{compatibility/dia} must belong to $\widehat \G$). A
maximal $(\G,\M)$-atlas on $X$ is called {\sl $(\G,\M)$-structure}.
\end{definition}

\begin{definition}\label{GMbranchfold/def}
By a {\sl geometric branchfold} modelled on the geometry $(\G,\M)$, in short a {\sl
$(\G,\M)$-branchfold}, we mean a pair $X_{\cal S} = (X,\cal S)$, where $X = X_{\cal B}$ is
a branchfold and $\cal S \subset \cal B$ is a $(\G,\M)$-structure on $X$. We will simply
write $X$ instead of $X_{\cal S}$, if no confusion can arise. Moreover, when talking of a
chart (resp. an atlas) of a $(\G,\M)$-branchfold $X = X_{\cal S}$, we will always assume
that it is a $(\G,\M)$-chart (resp. a $(\G,\M)$-atlas) in $\cal S$.
\end{definition}

The same argument we used to see that any branchfold atlas uniquely extends to a
branchfold structure (see page \pageref{atlas->structure}) also works in the geometric
context, to see that any $(\G,\M)$-atlas uniquely estends to a $(\G,\M)$-structure.
This is essentially due to the fact that restricion preserves $(\G,\M)$-charts.

\begin{definition}\label{GMcovering/def}
A branchfold map $f: X \to Y$ between $(\G,\M)$-branchfolds is called a {\sl
$(\G,\M)$-map}, if for every $(\G,\M)$-charts $(U, P, V, G = HK)$ and $(U', P', V',\break
G' = H'K')$ respectively of $X$ and $Y$ as in Definition \ref{branchfoldmap/def}, the map
$\psi: V \to V'$ in diagram \ref{map1/dia} belongs to $\widehat \G$. By {\sl
$(\G,\M)$-covering} (resp. {\sl $(\G,\M)$-isomorphism}) we mean a $(\G,\M)$-map which is a
branchfold covering (resp. isomorphism).
\end{definition}

We notice that the notions of $(\G,\M)$-branchfold and $(\G,\M)$-map (resp. covering)
defined above restrict to the usual ones in the case when referred to orbifolds/manifolds
and maps (resp. coverings) between them.

\medskip

As a consequence of the above observations about $(\G,\M)$-charts, for any point $x \in X$
of a $(\G,\M)$-branchfold $X$, there exist arbitrarily small reduced conical
$(\G,\M)$-charts. We call any of these reduced conical charts a {\sl local
$(\G,\M)$-model} of $X$ at $x$. This is uniquely determined up to $(\G,\M)$-isomorphisms
and conical restrictions (which are not necessarily $(\G,\M)$-isomorphisms).

\medskip

We remark that Proposition \ref{branchfoldcover/thm} does not hold in the geometric
context. Namely, given a branched covering $f: X \to Y$ onto a $(\G,\M)$-branchfold $Y$,
it does not necessarily exist a $(\G,\M)$-structure of $X$ making $f$ into a
$(\G,\M)$-covering. The reason is that, when lifting a $(\G,\M)$-chart of $Y$ through $f$,
we cannot always guarantee that property 1 of Definition \ref{GMchart/def} is preserved.
However, the following proposition says that the universal branchfold covering of a
$(\G,\M)$-branchfold can be thought as a $(\G,\M)$-covering (actually, the same is true
for all the branched coverings having the property described in Proposition
\ref{universalcover/thm}).

\begin{proposition}\label{GMunivcover/thm}
Let $X$ be a connected $(\G,\M)$-branchfold. Then the $(\G,\M)$-structure of $X$ can be
lifted to a unique $(\G,\M)$-structure on $\widetilde X$ making the universal branchfold
covering $u_X: \widetilde X \to X$ into a $(\G,\M)$-covering.
\end{proposition}

\begin{proof}
By Proposition \ref{universalcover/thm}, any local $(\G,\M)$-model $(U_x, P_x, V_x, G_x =
H_x K_x)$ of $X$ at $x = u_X(\widetilde x)$ lifts to a conical $(\G,\M)$-chart
$(U_{\widetilde x}, P_x, V_x, \widebar G_x = \widebar H_x K_x)$ on $\widetilde X$. In
fact, the branched covering $p_{K_x}$ is the same in both the charts, while the action of
$\widebar G_x/K_x$ on $V_x$ is a restriction of that of $G_x/K_x$.

Now, we can reason as in the proof of Proposition \ref{coverbranchfold/thm}, with the
extra requirement that all the chart isomorphisms are $(\G,\M)$-isomorphisms, to show that
these lifted conical $(\G,\M)$-charts form a $(\G,\M)$-atlas on $\widetilde X$ and that
the induced $(\G,\M)$-structure is the only one which makes $u_X$ into a
$(\G,\M)$-covering.
\end{proof}

Next proposition is the geometric version of Proposition \ref{quotientbranchfold/thm}.

\begin{proposition}\label{GMquotient/thm}
Let $X$ be a $(\G,\M)$-branchfold and $f: X \to Y$ be a regular branched covering, whose
deck transformations are $(\G,\M)$-isomorphisms which preserve local orientations at their
fixed points. Then the $(\G,\M)$-structure of $X$ induces a unique $(\G,\M)$-structure on
$Y$ making $f$ into a $(\G,\M)$-covering.
\end{proposition}

\begin{proof}
The same argument of the proof of Proposition \ref{quotientbranchfold/thm} works here,
once we replace charts by $(\G,\M)$-charts and chart isomorphisms by
$(\G,\M)$-isomorphisms.\break In particular, referring to the notations of that proof, if
$(U, P, V, G = HK)$ is a $(\G,\M)$-chart of $X$, then the quotient chart $(\widebar U, P,
V, \widebar G = \widebar H K)$ on $Y$ is a $(\G,\M)$-chart. In fact, the branched covering
$p_K$ is the same in both the charts, while the action of $\widebar G/K$ on $V$ is given
by elements of $\widehat \G$, as it can be easily derived from the fact that $L$ acts of
$U$ by $(\G,\M)$-isomorphisms.
\end{proof}

Before of going on, we briefly recall some constructions and notions concerning
$(\G,\M)$-manifolds, such as those of holonomy and developing map. We refer to
\cite{Thu97} (cf. also \cite{MM91}) for more details.

\medskip

Let $M$ be a connected $(\G,\M)$-manifold. A reduced $(\G,\M)$-chart $(U, P, V, G=HK)$
of $M$ can be thought in the usual way as $(U, \phi = p_K \circ p_H^{-1}: U \to V)$, being
$G$ the trivial group. We will simply denote by $(U,\phi)$ such a chart.

\label{alphaM/def}
Given any path $\alpha: \left[0,1\right] \to M$, we consider a sequence $0 = t_0 \leq t_1
\leq \dots \leq t_k = 1$ such that $\alpha([t_{i-1},t_i])$ is contained in some reduced
$(\G,\M)$-chart $(U_i,\phi_i)$ of $M$, for every $i = 1, \dots, k$. Taking into account
the compatibility condition between $(\G,\M)$-charts, we can assume that $\phi_i$ and
$\phi_{i+1}$ coincide in a neighborhood of $\alpha(t_i)$ (up to composition by elements of
$\G$). The local effectiveness of the action of $\G$ on $\M$ implies that this can be done
in a unique way, once the first chart $(U_1,\phi_1)$ is choosen. Then, we can
define a path $\alpha_\M: [0,1] \to \M$, by putting $\alpha_\M(t) = \phi_i(\alpha(t))$ for
$t \in [t_{i-1},t_i]$. Moreover, we can define a continuous family of local
$(\G,\M)$-models $\{(U'_t, \phi'_t)\}_{t \in [0,1]}$, such that $(U'_t, \phi'_t)$ is the
local model at $\alpha(t)$ induced by $(U_i,\phi_i)$ if $t \in [t_{i-1},t_i]$, in such a
way that $\alpha_\M(t) = \phi'_t(\alpha(t))$ for every $t \in [0,1]$.

By the local effectiveness of the action of $\G$ on $\M$, the path $\alpha_\M$ and the
family $\{(U'_t, \phi'_t)\}_{t \in [0,1]}$ are well defined up to multiplication by
elements of $\G$, depending only on the choice of the local model $(U'_0, \phi'_0)$ at the
starting point $\alpha(0)$.

The {\sl holonomy} $H_M: \pi_1(M) \to \G$ of the $(\G,\M)$-manifold $M$ is the
homomorphism defined as follows, once a local model $(U_\ast,\phi_\ast)$ at the base point
$\ast$ is fixed. For any loop $\omega$ in $(M, \ast)$ we consider the family of local
models $\{(U'_t, \phi'_t)\}_{t \in [0,1]}$ constructed as above, starting from $(U'_0,
\phi'_0) = (U_\ast,\phi_\ast)$. By the compatibility condition between $(\G,\M)$-charts,
there exists a unique $g \in \G$ such that $\phi'_1 = g \circ \phi'_0$ in a neighborhood
of $\ast\,$. Then, we put $H_M([\omega]) = g$.

As it can be easily seen, the holonomy $H_M$ is defined only up to conjugation in
$\G$, depending on the choice of the local model $(U_\ast,\phi_\ast)$ at the base point.

When the holonomy is trivial, we can define a $(\G,\M)$-map $D_M: M \to \M$, by\break
putting $D_M(x) = \alpha_\M(1)$, where $\alpha_\M$ is the path in $\M$ associated by the
above contruction, with $\phi'_0 = \phi_\ast\,$, to any path $\alpha$ in $M$ from $\ast$
to $x$.

This is called a {\sl developing map} for $M$. In this case, different choices of
$(U_\ast, \phi_\ast)$ lead to developing maps which differ by an element of $\G$. Namely,
by the local effectiveness of the action of $\G$ on $\M$, for any other $(\G,\M)$-map
$D'_M: M \to \M$ there exists an element $g \in \G$ such that $D'_M = g \circ D_M$.

In particular, for any connected $(\G,\M)$-manifold $M$ there always exists the developing
map $D_\wtM: \wtM \to \M$ defined on the universal covering $\wtM$ of $M$. Then, for every
path $\alpha: [0,1] \to M$, we can realize $\alpha_\M$, up to multiplication by elements
of $\G$, as $D_\wtM \circ \widetilde \alpha$, where $\widetilde \alpha: [0,1] \to
\wtM$ is a lifting of $\alpha$ through $u_M: \wtM \to M$.

\medskip

Now, we go back to $(\G,\M)$-branchfolds. Given a connected $(\G,\M)$-branchfold $X$, we
consider $X - \Sigma X$ as a connected $(\G,\M)$-manifold with the $(\G,\M)$-structure
induced by the inclusion in $X$.

\begin{definition}\label{holonomy/def}
We denote by $H_X: \pi_1(X - \Sigma X) \to \G$ the holonomy of $X - \Sigma X$ and we call
it the {\sl holonomy} of the $(\G,\M)$-branchfold $X$. Moreover, we denote by $H_x:
\pi_1(U_x - \Sigma U_x) \to \G$ the holonomy of the conical neighborhood $U_x$ of $x \in
X$ and we call it the {\sl local holonomy} of $X$ at $x$.
\end{definition}

Finally, we give a notion of completeness for $(\G,\M)$-branchfolds. We first observe that
the above construction of the path $\alpha_\M$ can be adapted, with a possibly infinite
sequence of $t_i$'s, in order to associate to any half open path $\alpha: \left[0,1\right[
\to M$ an open path $\alpha_\M: \left[0,1\right[ \to \M$ well defined up to multiplication
by elements of $\G$.

\begin{definition}\label{completeness/def}
We say that a $(\G,\M)$-branchfold $X$ is {\sl complete}, when any half open path
$\alpha: \left[0,1\right[ \to X - \Sigma X$, such that $\alpha_\M: \left[0,1\right[ \to
\M$ completes to a rectifiable path $\widebar \alpha_\M: [0,1] \to \M$, admits a
(rectifiable) completion $\widebar \alpha: [0,1] \to X$.
\end{definition}

As observed at page \pageref{rectifiable/def}, in the above definition we can equivalently
adopt the notion of rectifiability in $\M$ as a smooth manifold or as a polyhedron, with
the polyhedral structure given by any smooth triangulation of it.

Furthermore, when $\M$ admits a $\G$-invariant Riemannian metric, the $(\G,\M)$-manifolds
$X - \Sigma X$ and $R \subset \widetilde X$ can be endowed (in a unique way) with
Riemannian metrics which make $D_R: R \to \M$ and $u_{|X}: R \to X - \Sigma X$ into local
isometries. Now, the corresponding geodesic distances can be completed by continuity to
distances on $X$ and $\widetilde X$. So, it makes sense to compare our notion of
completeness with the metric one. Actually, a standard argument shows that they coincide
(cf. \cite{Thu97}).

In particular, if $\M$ admits a $\G$-invariant Riemannian metric, then any compact
$(\G,\M)$-branchfold is complete. The following proposition tells us this is true even
if such $\G$-invariant metric does not exist.

\break

\begin{proposition}
Any compact $(\G,\M)$-branchfold is complete.
\end{proposition}

\begin{proof}
Let $\alpha: \left[0,1\right[ \to X - \Sigma X$ be a half open path such that $\alpha_\M:
\left[0,1\right[ \to \M$ completes to a rectifiable path $\widebar \alpha_\M: [0,1] \to
\M$. Then, $\alpha_\M = D_R \circ \widetilde \alpha$, where $\widetilde \alpha: \left[0,1
\right[ \to R$ is a lifting of $\alpha$ through the ordinary covering $u_{X|}: R \to X -
\Sigma X$. Since $\widebar \alpha_\M$ is rectifiable, $\alpha_\M$ has finite length with
respect to any Riemannian metric on $\M$ (here, we do not need $\G$-invariance). This
metric can be lifted to $R$ through the local diffeomorphism $D_R$, and $\widetilde
\alpha$ has finite length with respect to this lifted metric, having the same length of
$\alpha_\M$.

On the other hand, $\widetilde\alpha(\left[0,1\right[)$ is contained in a compact
subpolyhedron $C \subset \widetilde X$, so it is rectifiable also with respect to the
polyhedral structure of $\widetilde X$. Then, $\alpha$ itself is rectifiable with respect
to the polyhedral structure of $X$, being $u_X$ a PL map.

Now, by the compactness of $X$, the half open path $\alpha$ admits some limit point $x \in
X$, meaning that there exists a sequence $t_n \to 1$ such that $\lim_{n \to \infty}
\alpha(t_n) = x$. To finish the proof, if suffices to observe that such limit point must
be unique, otherwise one could easily conclude that $\alpha$ would not be rectifiable.
\end{proof}

\subsection*{The geometric goodness theorem}

The universal branchfold covering $u_X: \widetilde X \to X$ of a connected
$(\G,\M)$-branchfold $X$ can be thought as a $(\G,\M)$-covering, by putting on $\widetilde
X$ the $(\G,\M)$-branchfold struc\-ture given by Proposition \ref{GMunivcover/thm}. The
restriction of $u_X$ over $X - \Sigma X$ is the ordinary regular covering $r: R \to X -
\Sigma X$ corresponding to the characteristic group $\Gamma_X \nsubgr \pi_1(X - \Sigma
X)$, which we introduced at page \pageref{chargroup/def} just in order to define $u_X$ as
its completion. Then, $R = \widetilde X - u_X^{-1}(\Sigma X)$ can be endowed with the
$(\G,\M)$-manifold structure induced by the inclusion in $\widetilde X$. This makes $r$
into a $(\G,\M)$-covering between $(\G,\M)$-manifolds.

\begin{proposition}\label{holonomy/thm}
Let $X$ a connected $(\G,\M)$-branchfold. Then, $\Gamma_x = \Ker H_x$ for every $x \in X$
and $\Gamma_X \subgr \Ker H_X$. Hence, the holonomy $H_R$ of the $(\G,\M)$-manifold $R =
\widetilde X - u_X^{-1}(\Sigma X)$ is trivial and $R$ admits a developing map $D_R: R \to
\M$.%
\end{proposition}

\begin{proof}
Given $x \in X$ and any conical chart $(U_x, P, V, G = HK)$ centered at $x$, we have
$\Gamma_x = i_*((p_{|*})^{-1}(\Im p_{G/K|*}))$ by definition, where $p_|: U_x - p_H(S_g)
\to P/G - B_G$ and $p_{G/K|}: V - p_K(S_G) \to P/G - B_G$ are the ordinary coverings of
$(\G,\M)$-manifolds given by restrictions of the branched coverings $p$ and $p_{G/K}$ of
the chart, and $i: U_x - p_H(S_G) \to U_x - \Sigma U_x$ is the inclusion.

First of all, we observe that the holonomy $H_{V - p_K(S_G)}$ is trivial, being $V -
p_K(S_G)$ an open subset of $\M$. Then, the equality $H_{V - p_K(S_G)} = H_{P/G - B_G}
\circ p_{G/K|*}$ implies that $\Im p_{G/K|*} \subgr \Ker H_{P/G - B_G}$. On the other
hand, for any loop $\omega$ in $P/G - B_G$, we can construct $\omega_\M$ inside $V -
p_K(S_G) \subset \M$, as a lifting of $\omega$ through the ordinary covering $p_{G/K|}$.
In particular, if $[\omega] \in \Ker H_{P/G - B_G}$ then $\omega_\M$ must be a loop in $V
- p_K(S_G)$, hence $[\omega] \in \Im p_{G/K|*}$. This proves that actually $\Im p_{G/K|*}
= \Ker H_{P/G - B_G}$.

Therefore, we have $\Gamma_x = i_*((p_{|*})^{-1}(\Im p_{G/K|*})) = i_*((p_{|*})^{-1} 
(\Ker H_{P/G - B_G}) = i_*(\Ker H_{U_x - p_H(S_G)}) = \Ker H_x$, where the second equality
derives from $H_{U_x - p_H(S_G)} = H_{P/G - B_G} \circ p_{|*}$ and the last one from the
surjectivity of $i_*$.

At this point, recalling that $\Gamma_X$ is normally generated by the groups
$h_x(\Gamma_x)$ with $x \in X$, the inclusion $\Gamma_X \subgr \Ker H_X$ immediately
follows from the fact that $h_x(\Ker H_x)$ is obviously contained in $\Ker H_X$ for every
$x \in X$.

Now, we can conclude that $H_R$ is trivial, being $H_R = H_X \circ r_*$ and $\Im r_* =
\Gamma_X$, where $r$ is the restriction of $u_X$ over $X - \Sigma X$ (cf. discussion
above).
\end{proof}

\label{holrep/def}
As an immediate consequence of the previous proposition, the holonomy $H_X$ factorizes
through a {\sl holonomy representation} $R_X: \pi_1^b(X) \to \G$ of the branchfold
fundamental group $\pi_1^b(X) = \pi_1(X - \Sigma X)/\Gamma_X$, such that $D_R \circ \gamma
= R_X(\gamma) \circ D_R$ for every $\gamma \in \pi_1^b(X)$. In other words, $D_R$ turns
out to be $R_X$-equivariant.

\begin{proposition}\label{developing/thm}
Let $X$ a connected $(\G,\M)$-branchfold. Then the developing map $D_R: R \to \M$ extends
to an $R_X$-equivariant $(\G,\M)$-map $C_X: \widetilde X \to \M$. Moreover, $C_X$ is a
$(\G,\M)$-covering if and only if $X$ is complete.
\end{proposition}

\begin{proof}
We define the map $C_X$ by local completion of $R_X$. Namely, given any $\widetilde x \in
\widetilde X - R$, we consider the local model $(U_x, P_x ,V_x, G_x = H_x K_x)$ of $X$ at
$x = u_X(\widetilde x)$ and the conical chart $(U_{\widetilde x}, P_x ,V_x, \widebar G_x =
\widebar H_x K_x)$ of $\widetilde X$ centered at $\widetilde x$, as in Proposition
\ref{universalcover/thm}. Then, $D_R \circ p_{\widebar H_x}$ and $p_{K_x}$ both restrict
to a developing map of the $(\G,\M)$-manifold $P_x - S_{G_x}$. Hence, there exists $g \in
\G$ such that $D_R \circ p_{\widebar H_x|} = g \circ p_{K_x|}: P_x - S_{G_x} \to V_x -
p_{K_x}(S_{G_x})$. Therefore $\widebar H_x \subgr K_x$, so that we can identify $P_x/
\widebar G_x$ with $V_x \subset \M$. Under this identification, the completion of $D_{R|}:
U_{\widetilde x} - p_{\widebar H_x}(S_{G_x}) \to \M$ is given by $g \circ \widebar p_x:
U_{\widetilde x} \to \M$, where $\widebar p_x: U_{\widetilde x} \to P_x/\widebar G_x$ is
the covering associated to the chart $(U_{\widetilde x}, P_x ,V_x, \widebar G_x = \widebar
H_x K_x)$ (notice that $p_{\widebar H_x}(S_{G_x})$ is a good subpolyhedron of
$U_{\widetilde x}$, by Proposition \ref{covergood/thm}). The unicity of completions (cf.
Proposition \ref{completion/thm}), guarantees that all this local completions fit together
to give a locally PL map $C_X: \widetilde X \to \M$ that extends $R_X$. Actually, the
above construction also tells us that $C_X$ is a $(\G,\M)$-map.

Now, assume that $X$ is complete. In order to conclude that $C_X$ is a $(\G,\M)$-covering,
we have to prove its completeness with respect to lifting of rectifiable paths. Let
$\alpha: [0,1] \to \M$ be a rectifiable path and $\widetilde \alpha_| : \left[0,1\right[
\to \widetilde X$ be a partial lifting of it through $C_X$. Taking into account that
$u_X^{-1}(\Sigma X)$ is a good subcomplex of $\widetilde X$, we can perturb $\widetilde
\alpha_|$ to a half open path $\beta: \left[0,1\right[ \to R = \widetilde X -
u_X^{-1}(\Sigma X)$ such that $D_R \circ \beta$ is rectifiable and $\lim_{t \to 1}
d(\beta(t), \widetilde\alpha_| (t))) = 0$ for some metric $d$ on $\widetilde X$. Then, for
$u_X \circ \beta: \left[0,1\right[ \to X - \Sigma X$, we have that $(u_X \circ \beta)_\M =
D_R \circ \beta$ is rectifiable. Hence, by the completeness of $X$, there exists an
extension $\gamma: [0,1] \to X$ of $u_X \circ \beta$. Since $u_X$ is a branched covering,
$\gamma$ lifts to $\widetilde \gamma: [0,1] \to \widetilde X$, which extends $\beta$.
Then, there exists $\lim_{t \to 1} \widetilde \alpha_|(t) = \lim_{t \to 1} \beta(t) =
\widetilde \gamma(1) \in \widetilde X$ and we can complete $\widetilde \alpha_|$ to a
lifting $\widetilde \alpha$ of $\alpha$.

Viceversa, assume that $C_X$ is a $(\G,\M)$-covering. Let $\alpha: \left[0,1\right[ \to X
- \Sigma X$ a half open path such that $\alpha_\M$ completes to a rectifiable path
$\widebar \alpha_\M: [0,1] \to \M$. Let $\widetilde \alpha: \left[0,1\right[ \to R$ a
lifting of $\alpha$ through the ordinary covering $r = u_{X|}: R \to X - \Sigma X$. Then,
$C_X \circ \widetilde \alpha = g \circ \alpha_\M$, for a suitable $g \in \G$, is
rectifiable. By the completeness of $C_X$, $\widetilde \alpha$ admits a completion and so
also $\alpha = u_X \circ \widetilde \alpha$ does.
\end{proof}

It is worth remarking that the previous proposition (and more directly its proof) implies
that the universal branchfold covering $\widetilde X$ of a $(\G,\M)$-branchfold $X$ is
always a pure branchfold.

At this point we are ready for the announced geometric goodness theorem.

\begin{theorem}\label{goodness/thm}
Any connected compact $(\G,\M)$-branchfold $X$, whose holonomy group $\H_X =
R_X(\pi_1^b(X)) = H_X(\pi_1(X - \Sigma X)) \subgr \G$ acts properly discontinuously on
$\M$, is a good branchfold. In fact, there exists a PL $(\G,\M)$-covering $p: X \to O_X$
onto the good $(\G,\M)$-orbifold $O_X = \M/\H_X$, such that the following diagram
commutes.
\end{theorem}

\begin{Diagram}[htb]{good3/dia}
\diag{}{good3.eps}{}
\end{Diagram}

\begin{proof}
The existence of the map $p: X \to O_X$, follows from the fact that $C_X$ is
$R_X$-equivariant, being a completion of the $R_X$-equivariant map $D_R$.

In order to prove that $p$ is a $(\G,\M)$-map, let $(U_x, P_x, V_x, G_x = H_x K_x)$ the
local model of $X$ at any point $x \in X$ and $(U_{\widetilde x}, P_x, V_x, \widebar G_x =
\widebar H_x K_x)$ be the conical chart of $\widetilde X$ at any point $\widetilde x \in
\widetilde X$ such that $u_X(\widetilde x) = x$, as in Proposition
\ref{universalcover/thm}. Then, we consider the open subsets $C_X(U_{\widetilde x}) =
g(V_x) \subset \M$ and $p(U_x) = \pi_\H(g(V_x)) \subset O_X$, where $g$ is a suitable
element of $\G$ (cf. proof of Proposition \ref{developing/thm}). By Proposition
\ref{restaction/thm}, the restriction $\pi_{\H|}: g(V_x) \to p(U_x)$ is a regular branched
covering induced by the action of a finite subgroup $L \subgr \H$, hence we can identify
it with the canonical projection $\pi_L$. Thus, we have the following commutative diagram,
telling us that $p$ is $(\G,\M)$-map at $x$.

\begin{Diagram}[htb]{good4/dia}\vskip-9pt
\diag{}{good4.eps}{}\vskip3pt
\end{Diagram}

Finally, the compactness of $X$ allows us to conclude that $p$ is a PL map, hence a PL
$(\G,\M)$-covering.
\end{proof}

Notice that the compactness of $X$ was used in the proof of the theorem only for
concluding that $p$ is a PL map, hence a PL branchfold covering. However, we emphasize
that the completeness of $X$ would suffice for $p$ to be a branchfold covering.

\medskip

We have already observed at the end of the previous section, that all branchfolds are
locally very good. However, the construction of diagram \ref{good3/dia}, once adapted to
the non-compact (and even non-complete) context of branchfold charts, give us some more
information in the case of $(\G,\M)$-branchfold. In particular, it allows us to see how
the local holonomy group of a $(\G,\M)$-branchfold at a point is related to the conical
$(\G,\M)$-charts centered at that point. This relation is stated by the following
proposition, which will be needed in the next subsection.

\begin{proposition}\label{GMconical/thm}
For any conical $(\G,\M)$-chart $(U_x, P, V, G = HK)$ centered at $x$, the local holonomy
group $\H_x = H_x(\pi_1(U_x - \Sigma U_x))$ coincides with $G/K \subgr \G$.%
\end{proposition}

\begin{proof}
First of all, we observe that $\H_x \subgr G/K$. In fact, the homomorphism $i_*: \pi_1(U_x
- p^{-1}(\Sigma (P/G))) \to \pi_1(U_x - \Sigma U_x)$ induced by the inclusion is
surjective and $H_x \circ i_* = H_{P/G} \circ p_{|*}$, where $p_|: U_x - p^{-1}(\Sigma
(P/G)) \to P/G - \Sigma (P/G)$ is the restriction over $P/G - \Sigma (P/G)$ of the
covering $p: U_x \to P/G$ associated to the chart. Hence, $\H_x = \Im H_x \subgr \Im
H_{P/G} = \H_{P/G} = G/K \subgr \G$.

To see that actually $\H_x = G/K$, look at the commutative diagram \ref{good5/dia}. Here,
we have the diagram of the $(\G,\M)$-chart $(U_x, P, V, G = HK)$ and the maps of diagram
\ref{good3/dia}. Since $P$ is simply connected pure branchfold covering of $U_x$,
Proposition \ref{mostpurecover/thm} tells us that there exists a PL branchfold covering
$l: P \to \widetilde U_x$ such that $u_{U_x} \circ l = p_H$. By a suitable choice of
$C_{U_x}$, we can also assume that $C_{U_x} \circ l = p_K$. Moreover, being $\H_x \subgr
G/K$ a finite group, $\pi_{\H_x}: V \to O$ is a finite PL branchfold covering onto an
orbifold $O$ and there exists a PL branchfold covering $t: O \to P/G$. Analogously,
$C_{U_X}/R_{U_x}$ is a PL branchfold covering, since it is a finite map.

\begin{Diagram}[htb]{good5/dia}\vskip6pt
\diag{}{good5.eps}{}\vskip6pt
\end{Diagram}

Now, let $q: Q \to O$ be the pullback of $C_{U_x}/R_{U_x}$ and $\pi_{\H_x}$, with the
associated PL branchfold coverings $q_1$ and $q_2$, and $q \circ r: R \to O$ be the
minimal regularization of $q$. Then, we define the coverings $p_{H'} = q_1 \circ r$ and
$p_{K'} = q_2 \circ r$, with deck transformations groups $H'$ and $K'$ respectively.
Denoting by $L$ the group of deck transformations of $q \circ r$, we have $H' \subgr L$
and $K' \nsubgr L$. Hence, we can consider the subgroup $G' = H' K' \subgr L$ and the
branchfold chart $(U_x, R, V, G' = H'K')$. By the universal property of the pullback,
there exists a PL branchfold covering $m: P \to Q$ which commutes with the other maps.
Furthermore, $q \circ m$ is regular, since so is $t \circ q \circ m = \pi_G$. Hence, by
the universal property of the minimal regularizations, there exists a factorization $m = r
\circ n$, for a PL branchfold covering $n: P \to R$. This, gives a domination of the chart
$(U_x, P, V, G = HK)$ on the chart $(U_x, R, V, G' = H'K')$. As a consequence, this last
chart is a $(\G,\M)$-chart and $\pi_G = \pi_{G'} \circ n: P \to P/G \cong R/G'$. On the
other hand, since $G' \subgr L$ there is a PL branchfold covering $P/G \cong R/G' \to O$,
which commutes with the other maps. Then, we can conclude that $G/K \subgr \H_x$.
\end{proof}

\subsection*{Rational conifolds as geometric branchfolds}

In this subsection we focus on the geometric branchfolds modelled on constant curvature
Riemannian geometries. In particular, we apply the geometric goodness theorem to establish
the relation between such branchfolds and conifold spaces.

\medskip

\label{riemgeom/def}
We denote by $\M_k^m$ the $m$-dimensional model Riemannian space of constant curvature $k$
and by $\G_k^m$ its isometry group, for any real number $k$. So, the geometric branchfolds
we want to consider are the $(\M_k^m,\G_k^m)$-branchfolds.

\medskip

\label{conifold/def}
Before going on, we recall the definition of {\sl conifold of dimension $m$ and curvature
$k$}, in short {\sl $(m,k)$-conifold}. This is given by induction on the dimension $m$ as
follows. The $(1,k)$-conifolds are the circles of any length, independently on $k$. For $m
\geq 2$, an $(m,k)$-conifold $X$ is a complete metric space locally modelled on $k$-cones
over $(m-1,1)$-conifolds. More precisely, for any $x \in X$ there exist $\epsilon > 0$ and
an isometry between the open ball $B(x,\epsilon)$ and the open $k$-cone
$C_{k,\epsilon}(L_x)$ on a connected compact $(m-1,1)$-conifold $L_x$, letting $x$
correspond to the apex of the cone.

\label{k-cone/def}
Here, by the open {\sl $k$-cone} $C_{k,r}(L)$ of a metric space $L$, with $r \leq k/\sqrt
\pi$ if $k > 0$, we mean the open cone $L \times \left[0,r\right[ / L \times \{0\}$
endowed with the metric $d((x_1,t_1),(x_2,t_2)) = d_{\M_k^2}(p_1,p_2)$, for a geodesic
triangle $p_0,p_1,p_2$ in $\M_k^2$, such that $\angle_{p_0} = \min(d_L(x_1,x_2),\pi)$,
$d(p_0,p_1) = t_1$ and $d(p_0,p_2) = t_2$.

\medskip

According to what we have seen in a more general context (when discussing the notion of
completeness at page \pageref{completeness/def}), any $(\G_k^m,\M_k^m)$-branchfold $X$ can
be endowed with a natural metric, whose restriction to $X - \Sigma X$ has constant
curvature $k$. Our first aim is to prove that such a metric makes $X$ into a
$(m,k)$-conifold.

\medskip

\label{metric/def}
We define the {\sl natural metric} on a $(\G_k^m,\M_k^m)$-branchfold $X$, in the following
way. We start by lifting the metric of $\M_k^m$ to $R$ through the developing map $D_R: R
\to \M$. Then, by the $R_X$-invariance of $D_R$ (cf. page \pageref{holrep/def}), the deck
transformations of the ordinary covering $r = u_{X|} : R \to X - \Sigma X$ preserve such
metric, hence we have an induced metric on $X - \Sigma X$ making $r$ into a local
isometry. Finally, taking into account that $\Sigma X \subset X$ is a good subpolyhedron,
we can extend by continuity the geodesic distance on $X - \Sigma X$ to a distance on $X$.

Actually, in the light of Theorem \ref{goodness/thm}, this natural metric on $X$ can
be obtained by putting the quotient metric on $O_X = \M/\H_X$ and then lifting this metric
to $X$ through the branched covering $p: X \to O_X$.

The above constructions of the natural metric on $X$ can be also performed locally at $x
\in X$. More precisely, for a local $(\G_k^m,\M_k^m)$-model $(U_x, P_x, V_x, G_x = H_x
K_x)$, we can either lift the metric of $V_x$ to $P_x$ through $p_{K_x}$ and then consider
the quotient metric on $U_x = P_x/H_x$, or lift to $U_x$ through $p_x$ the metric on
$P_x/G_x \cong V_x/(G_x/K_x)$ induced by the last quotient. In both cases, if $U_x$ is
a sufficiently small convex conical neighborhood of $x$, we get the restriction to $U_x$
of the natural metric of $X$.

Now, let us consider the linearization of $(U_x, P_x, V_x, G_x = H_x K_x)$ at $x$. This is
the local model $(T_x U_x, T_{\widetilde x} P_x, T_{\widebar x} V_x, G_x = H_x K_x)$,
where $T$ denotes the tangent cone, $\widetilde x$ and $\widebar x$ are respectively the
apices of $P_x$ and $V_x$, while the action of $G_x$ on $T_{\widetilde x} P_x$ is the
unique which preserves the radial structure and corresponds to the original one on $P_x$
in a neighborhood of $\widetilde x$ through the exponential map. We emphasize that
property (1) in Definition \ref{GMchart/def} is essential for the existence of such
action, since $\G_k^m$-admissible subpolyhedra of $\M_k^m$ stratify by totally geodesic
submanifolds.

In particular, we can identify $(T_{\widebar x} V_x,\widebar x)$ with $(R^m,0)$, in such
a way that under this identification $G_x/K_x$ acts on it as a subgroup of $\SO(m)$. Then,
we define\break $T^1_{\widetilde x} P_x = (T_{\widetilde x} p_{K_x})^{-1}(S^{m-1}) \subset
T_{\widetilde x} P_x$, $T^1_{\widehat x}(P_x/G_x) = T_{\widebar x}p_{G_x/K_x}(S^{m-1})
\subset T_{\widehat x}(P_x/G_x)$ and $L_x X = T_{\widetilde x}p_{H_x}(T^1_{\widetilde x}
P_x) = (T_x p)^{-1}(T^1_{\widehat x}(P_x/G_x))$. Clearly, $T^1_{\widehat x}(P_x/G_x) =
S^{m-1}/(G_x/K_x)$ is a very good $(\G_1^{m-1},\M_1^{m-1})$-orbifold, hence $L_x X$ is a
very good $(\G_1^{m-1},\M_1^{m-1})$-branchfold.

\begin{proposition}\label{brachfold->conifold/thm}
Any $(\G_k^m,\M_k^m)$-branchfold $X$ with its natural metric is a $(m,k)$-conifold. In
fact, for every point $x \in X$ there exists $\epsilon > 0$, such that the open geodesic
ball $B(x,\epsilon)$ is isometric to the $k$-cone $C_{k,\epsilon}(L_x X)$, by an isometry
letting $x$ correspond to the apex of the cone.
\end{proposition}

\begin{proof}
We proceed by induction on the $m$. If $m = 1$ there is nothing to prove. So, we assume $m
> 1$ and consider any point $x \in X$. Since the restriction of the natural metric of $X$
to $L_x X$ coincides with the natural metric of $L_x X$ itself as a
$(\G_1^{m-1},\M_1^{m-1})$-branchfold, by the induction hypothesis tells us that $L_x X$ is
a $(m-1,1)$-conifold. Now, we look at the following diagram.

\begin{Diagram}[htb]{conifold1/dia}\vskip3pt
\diag{}{conifold1.eps}{}
\end{Diagram}

Here, the vertical arrows on the right side are restrictions of the coverings $p_{H_x}$
and $p_{K_x}$ associated to the $(\G_k^m,\M_k^m)$-chart $(U_x, P_x, V_x, G_x = H_x K_x)$,
while those on the left side are obtained by applying the cone construction
$C_{k,\epsilon}$ to restrictions of the corresponding maps associated to the linearization
of such $(\G_k^m,\M_k^m)$-chart. Moreover, the horizontal arrows are all induced by the
exponential map on the top. Since this is an isometry and the vertical arrows are local
isometries out of the singularities, we can easily conclude that also the other two
horizontal maps are isometries.
\end{proof}

In order to characterize the conifold that can be obtained from branchfold as in the above
proposition, we need the notions of (local) holonomy of a conifold.

First of all, we observe that any $(m,k)$-conifold $X$ is a speudo-manifold of dimension
$m$. Moreover, the {\sl singular locus} $\Sigma X = \{x \in X \;|\; X$ is not an
$m$-manifold at $x\}$ has dimension $\leq m - 2$, hence it is a good subpolyhedron of $X$.
On the other hand, the complement $X - \Sigma X$ is a $(\G_k^m,\M_k^m)$-manifold.

Then, we the holonomy $H_X: \pi_1(X - \Sigma X) \to \G_k^m$ is defined and we call it the
{\sl holonomy} of $X$. Analogously, the holonomy $H_x: \pi_1(B(x,\epsilon) - \Sigma
B(x,\epsilon)) \to \G_k^m$ is defined for $\epsilon > 0$ sufficiently small and we call it
the {\sl local holonomy} of $X$ at $x$.

\medskip

\label{rational/def}
By a {\sl rational conifold} we mean a conifold $X$ such that the local holonomy group
$\H_x = \Im H_x \subgr \G_k^m$ is finite for every $x \in X$. The reason of this
terminology is that, the local holonomy at a codimension $2$ point $x \in \Sigma X$ is
finite if and only if the singular angle of $X$ at $x$ is a rational multiple of $\pi$
radians. In other words, the local model of $X$ at $x$ is like that shown in Figure
\ref{example1/fig} at page \pageref{example1/fig} crossed by $R^{m-2}$.

\begin{theorem}\label{conifoldtobranchfold/thm}
A connected $(m,k)$-conifold $X$ admits a $(\G_k^m,\M_k^m)$-branch\-fold structure as in
Proposition \ref{brachfold->conifold/thm}, if and only if it is a rational conifold.%
\end{theorem}

\begin{proof}
The ``only if'' part is quite trivial. In fact, the local holonomies of $X$ as a conifold
and as a branchfold are the same. Then, the finiteness of $\H_x$, for every $x \in X$,
derives from Proposition \ref{GMconical/thm}.

We prove the ``if'' part by induction on $m$, starting from the trivial case $m = 1$. So,
let us assume that $X$ is a connected $(m,k)$-conifold with $m > 1$. Given any $x \in X$,
there exists $\epsilon > 0$ such that $B(x,\epsilon) \cong C_{k,\epsilon}(L_x)$ for some
connected compact $(m-1,1)$-conifold $L_x$. Moreover, the holonomy group $\H_{L_x}$ of
$L_x$ is finite, since it coincides with the local holonomy group $\H_x$ of $X$ at $x$.
By the induction hypothesis, $L_x$ is a $(\G_1^{m-1},\M_1^{m-1})$-branchfold.
Then, Theorem \ref{goodness/thm} tells us that $L_x$ is a good branchfold. Actually, $L_x$
is very good, being a finite branchfold covering of the very good orbifold
$\G_1^{m-1}/\H_x = S^{m-1}/\H_x$. According to our comment to diagram \ref{good2/dia},
this allows us to construct the following diagram of finite branchfold coverings.

\begin{Diagram}[htb]{conifold2/dia}
\diag{}{conifold2.eps}{}
\end{Diagram}

By applying the $k$-cone construction $C_{k,\epsilon}$ to all spaces and maps in the
diagram, we get a conical $(\G_k^m,\M_k^m)$-chart $(U_x, P_x, V_x, G = HK)$ for $X$ at
$x$, where $U_x = B(x,\epsilon) \cong C_{k,\epsilon}(L_x)$, $P_x = C_{k,\epsilon}(P)$,
$V_x = C_{k,\epsilon}(S^{m-1}) \cong B(0,\epsilon) \subset \M_k^m$ and $P_x/G \cong
C_{k,\epsilon}(O)$ (cf. Proposition \ref{GMconical/thm}).

In order to prove that the conical $(\G_k^m,\M_k^m)$-charts we have just constructed form
a $(\G_k^m,\M_k^m)$-atlas, it suffices to verify that for any two such charts $(U_x, P_x,
V_x, G = HK)$ and $(U_y, P_y, V_y, \widebar G = \widebar H \widebar K)$ with $U_x \subset
U_y$, we have that the former is equivalent to a restriction $(U_x, P'_x, V'_x, G' = H'
K')$ of the latter to $U_x$.

Since the above construction of $(\G_k^m,\M_k^m)$-charts for $X$ originates from
developing maps (cf. diagrams \ref{good2/dia} and \ref{good3/dia}), we can assume, up to
multiplication by elements of $\G$, that $V_x = V'_x \subset V_y$ and $G/K = \H_{U_x}
\subgr \H_{U_y} = \widebar G/\widebar K$, where the equalities follows from Proposition
\ref{GMconical/thm}. On the other hand, any restriction of a $(\G_k^m,\M_k^m)$-chart is
still a $(\G_k^m,\M_k^m)$-chart, so we also have $G'/K' = \H_{U_x}$, by applying once
again Proposition \ref{GMconical/thm}. Then, $P_x/G = P_x'/G'$ and we have the same PL
branchfold covering $p = C_{U_x}/R_{U_x}: U_x \to P_x/G = P_x'/G'$ associated to the
charts we are comparing. Therefore, those charts induce the same branchfold structure on
$U_x$, hence they are equivalent.
\end{proof}

\end{document}